\pdfoutput=1

\documentclass[12pt,a4paper]{article}

\usepackage{amsfonts,amssymb,amsthm,amsmath,latexsym,mathtools} 
\usepackage{eqsection,indent}
\usepackage{subeqnarray, eepicemu, url,cite,bm}
\usepackage{multirow}  
\usepackage{tensor}  
\usepackage{graphicx,graphbox}
\usepackage{float}  

\usepackage{tikz}
\usepackage{color}

\usetikzlibrary{matrix}
\usetikzlibrary{cd}
\usepackage{forest}
\usepackage{subfig}

\date{July 24, 2025}

\DeclareFontFamily{U}{mathx}{}
\DeclareFontShape{U}{mathx}{m}{n}{<-> mathx10}{}
\DeclareSymbolFont{mathx}{U}{mathx}{m}{n}
\DeclareMathAccent{\widecheck}{0}{mathx}{"71}

\begin{document}

\title{\vspace*{-0.5cm}
       Higher-order Stirling cycle and subset triangles: \\[2mm]
    \hspace*{-1.3cm}
  \hbox{\quad Total positivity, continued fractions and real-rootedness}
      }

\author{
      \hspace*{-1cm}
      {\large Bishal Deb${}^{1}$ and Alan D.~Sokal${}^{2,3}$}
   \\[5mm]
      \hspace*{-1.8cm}
      \normalsize
      ${}^1$Yau Mathematical Sciences Center, Tsinghua University, Beijing 100084, China\\
     \hspace*{-1.3cm}
      \normalsize
           ${}^2$Department of Mathematics, University College London,
                    London WC1E 6BT, UK   \\
     \hspace*{-2.9cm}
      \normalsize
           ${}^3$Department of Physics, New York University,
                    New York, NY 10003, USA
       \\[4mm]
     \hspace*{-1cm}
     {\tt bishal@gonitsora.com},
     {\tt sokal@nyu.edu}
}

\maketitle
\thispagestyle{empty}   

\begin{abstract}
Given a lower-triangular matrix of real numbers,
one can ask the following four total-positivity questions:
total positivity of the triangle itself;
total positivity of its row-reversal;
Toeplitz-total positivity of its row sequences
(equivalent to negative-real-rootedness of the row-generating polynomials);
and coefficientwise Hankel-total positivity of
the sequence of row-generating polynomials.
In this paper, 
we introduce two infinite families of lower-triangular matrices
generalising the Stirling cycle and subset triangles,
parametrised by an integer $r \ge 1$;
we call these the $r$th-order Stirling cycle and subset numbers.
We then ask the foregoing four questions for each of these triangles,
leading us to several conjectures.
We then prove some of these conjectures for the case $r=2$.
\end{abstract}

\medskip
\noindent
{\bf Key Words:}
Totally positive matrix, 
Hankel matrix, 
zeros of polynomials, 
higher-order Stirling cycle numbers, 
higher-order Stirling subset numbers,
Stirling permutations,
increasing ternary trees,
increasing ordered trees,
branched continued fractions,
permutations,
set partitions,
phylogenetic trees,
log-concavity.

\bigskip

\noindent
{\bf Mathematics Subject Classification (MSC 2020) codes:}
05A15 (Primary);
05A05, 05A18, 05A19, 05A20, 11B73, 15B48, 30B70 (Secondary).

\vspace*{1cm}

\newtheorem{theorem}{Theorem}[section]
\newtheorem{proposition}[theorem]{Proposition}
\newtheorem{lemma}[theorem]{Lemma}
\newtheorem{corollary}[theorem]{Corollary}
\newtheorem{definition}[theorem]{Definition}
\newtheorem{conjecture}[theorem]{Conjecture}
\newtheorem{question}[theorem]{Question}
\newtheorem{problem}[theorem]{Problem}
\newtheorem{openproblem}[theorem]{Open Problem}
\newtheorem{example}[theorem]{Example}
\newtheorem{remark}[theorem]{Remark}

\renewcommand{\theenumi}{\alph{enumi}}
\renewcommand{\labelenumi}{(\theenumi)}
\def\eop{\hbox{\kern1pt\vrule height6pt width4pt
depth1pt\kern1pt}\medskip}
\def\prf{\par\noindent{\bf Proof.\enspace}\rm}
\def\rmk{\par\medskip\noindent{\bf Remark\enspace}\rm}

\newcommand{\textbfit}[1]{\textbf{\textit{#1}}}

\newcommand{\bigdash}{%
\smallskip\begin{center} \rule{5cm}{0.1mm} \end{center}\smallskip}

\newcommand{\safepar}{ {\protect\hfill\protect\break\hspace*{5mm}} }
\newcommand{\ceil}[1]{\lceil #1 \rceil }

\newcommand{\be}{\begin{equation}}
\newcommand{\ee}{\end{equation}}
\newcommand{\<}{\langle}
\renewcommand{\>}{\rangle}
\newcommand{\widebar}{\overline}
\def\reff#1{(\protect\ref{#1})}
\def\spose#1{\hbox to 0pt{#1\hss}}
\def\ltapprox{\mathrel{\spose{\lower 3pt\hbox{$\mathchar"218$}}
    \raise 2.0pt\hbox{$\mathchar"13C$}}}
\def\gtapprox{\mathrel{\spose{\lower 3pt\hbox{$\mathchar"218$}}
    \raise 2.0pt\hbox{$\mathchar"13E$}}}
\def\textprime{${}^\prime$}
\def\proof{\par\medskip\noindent{\sc Proof.\ }}
\def\firstproof{\par\medskip\noindent{\sc First Proof.\ }}
\def\secondproof{\par\medskip\noindent{\sc Second Proof.\ }}
\def\alternateproof{\par\medskip\noindent{\sc Alternate Proof.\ }}
\def\algebraicproof{\par\medskip\noindent{\sc Algebraic Proof.\ }}
\def\graphicalproof{\par\medskip\noindent{\sc Graphical Proof.\ }}
\def\combinatorialproof{\par\medskip\noindent{\sc Combinatorial Proof.\ }}
\def\proofof#1{\bigskip\noindent{\sc Proof of #1.\ }}
\def\firstproofof#1{\bigskip\noindent{\sc First Proof of #1.\ }}
\def\secondproofof#1{\bigskip\noindent{\sc Second Proof of #1.\ }}
\def\thirdproofof#1{\bigskip\noindent{\sc Third Proof of #1.\ }}
\def\algebraicproofof#1{\bigskip\noindent{\sc Algebraic Proof of #1.\ }}
\def\combinatorialproofof#1{\bigskip\noindent{\sc Combinatorial Proof of #1.\ }}
\def\completionofproofof#1{\bigskip\noindent{\sc Completion of the Proof of #1.\ }}
\def\sketchofproof{\par\medskip\noindent{\sc Sketch of proof.\ }}
\renewcommand{\qed}{ $\square$ \bigskip}
\newcommand{\myendremark}{ $\blacksquare$ \bigskip}
\def\half{ {1 \over 2} }
\def\third{ {1 \over 3} }
\def\twothird{ {2 \over 3} }
\def\smfrac#1#2{{\textstyle{#1\over #2}}}
\def\smhalf{ {\smfrac{1}{2}} }
\newcommand{\real}{\mathop{\rm Re}\nolimits}
\renewcommand{\Re}{\mathop{\rm Re}\nolimits}
\newcommand{\imag}{\mathop{\rm Im}\nolimits}
\renewcommand{\Im}{\mathop{\rm Im}\nolimits}
\newcommand{\sgn}{\mathop{\rm sgn}\nolimits}
\newcommand{\tr}{\mathop{\rm tr}\nolimits}
\newcommand{\supp}{\mathop{\rm supp}\nolimits}
\newcommand{\disc}{\mathop{\rm disc}\nolimits}
\newcommand{\diag}{\mathop{\rm diag}\nolimits}
\newcommand{\tridiag}{\mathop{\rm tridiag}\nolimits}
\newcommand{\AZ}{\mathop{\rm AZ}\nolimits}
\newcommand{\EAZ}{\mathop{\rm EAZ}\nolimits}
\newcommand{\NC}{\mathop{\rm NC}\nolimits}
\newcommand{\PF}{{\rm PF}}
\newcommand{\rk}{\mathop{\rm rk}\nolimits}
\newcommand{\perm}{\mathop{\rm perm}\nolimits}
\def\hboxscript#1{ {\hbox{\scriptsize\em #1}} }
\renewcommand{\emptyset}{\varnothing}
\newcommand{\eqdef}{\stackrel{\rm def}{=}}

\newcommand{\restrict}{\upharpoonright}

\newcommand{\compinv}{{\langle -1 \rangle}}   

\newcommand{\scra}{{\mathcal{A}}}
\newcommand{\scrb}{{\mathcal{B}}}
\newcommand{\scrc}{{\mathcal{C}}}
\newcommand{\scrd}{{\mathcal{D}}}
\newcommand{\scrdtilde}{{\widetilde{\mathcal{D}}}}
\newcommand{\scre}{{\mathcal{E}}}
\newcommand{\scrf}{{\mathcal{F}}}
\newcommand{\scrg}{{\mathcal{G}}}
\newcommand{\scrh}{{\mathcal{H}}}
\newcommand{\scri}{{\mathcal{I}}}
\newcommand{\scrj}{{\mathcal{J}}}
\newcommand{\scrk}{{\mathcal{K}}}
\newcommand{\scrl}{{\mathcal{L}}}
\newcommand{\scrlbar}{{\overline{\mathcal{L}}}}
\newcommand{\scrm}{{\mathcal{M}}}
\newcommand{\scrn}{{\mathcal{N}}}
\newcommand{\scro}{{\mathcal{O}}}
\newcommand{\scroo}{{
  \mathchoice
    {{\scriptstyle\mathcal{O}}}
    {{\scriptstyle\mathcal{O}}}
    {{\scriptscriptstyle\mathcal{O}}}
    {\scalebox{0.6}{$\scriptscriptstyle\mathcal{O}$}}
  }}
\newcommand{\scrp}{{\mathcal{P}}}
\newcommand{\scrq}{{\mathcal{Q}}}
\newcommand{\scrr}{{\mathcal{R}}}
\newcommand{\scrs}{{\mathcal{S}}}
\newcommand{\scrt}{{\mathcal{T}}}
\newcommand{\scrtprime}{\mathcal{T}^{\,'}}
\newcommand{\scrv}{{\mathcal{V}}}
\newcommand{\scrw}{{\mathcal{W}}}
\newcommand{\scrx}{{\mathcal{X}}}
\newcommand{\scrz}{{\mathcal{Z}}}
\newcommand{\SP}{{\mathcal{SP}}}
\newcommand{\ST}{{\mathcal{ST}}}

\newcommand{\scrthat}{{\widehat{\mathcal{T}}}}

\newcommand{\bfa}{{\mathbf{a}}}
\newcommand{\bfb}{{\mathbf{b}}}
\newcommand{\bfc}{{\mathbf{c}}}
\newcommand{\bfd}{{\mathbf{d}}}
\newcommand{\bfe}{{\mathbf{e}}}
\newcommand{\bfh}{{\mathbf{h}}}
\newcommand{\bfj}{{\mathbf{j}}}
\newcommand{\bfi}{{\mathbf{i}}}
\newcommand{\bfk}{{\mathbf{k}}}
\newcommand{\bfl}{{\mathbf{l}}}
\newcommand{\bfL}{{\mathbf{L}}}
\newcommand{\bfm}{{\mathbf{m}}}
\newcommand{\bfn}{{\mathbf{n}}}
\newcommand{\bfp}{{\mathbf{p}}}
\newcommand{\bfr}{{\mathbf{r}}}
\newcommand{\bft}{{\mathbf{t}}}
\newcommand{\bfu}{{\mathbf{u}}}
\newcommand{\bfv}{{\mathbf{v}}}
\newcommand{\bfw}{{\mathbf{w}}}
\newcommand{\bfx}{{\mathbf{x}}}
\newcommand{\bfX}{{\mathbf{X}}}
\newcommand{\bfy}{{\mathbf{y}}}
\newcommand{\bfz}{{\mathbf{z}}}
\renewcommand{\k}{{\mathbf{k}}}
\newcommand{\n}{{\mathbf{n}}}
\newcommand{\vv}{{\mathbf{v}}}
\newcommand{\w}{{\mathbf{w}}}
\newcommand{\x}{{\mathbf{x}}}
\newcommand{\y}{{\mathbf{y}}}
\newcommand{\cc}{{\mathbf{c}}}
\newcommand{\zero}{{\mathbf{0}}}
\newcommand{\one}{{\mathbf{1}}}
\newcommand{\bmm}{{\mathbf{m}}}

\newcommand{\ahat}{{\widehat{a}}}
\newcommand{\Ghat}{{\widehat{G}}}
\newcommand{\Zhat}{{\widehat{Z}}}

\newcommand{\Pcheck}{{\widecheck{P}}}
\newcommand{\pcheck}{{\widecheck{p}}}

\newcommand{\C}{{\mathbb C}}
\newcommand{\D}{{\mathbb D}}
\newcommand{\Z}{{\mathbb Z}}
\newcommand{\N}{{\mathbb N}}
\newcommand{\Q}{{\mathbb Q}}
\newcommand{\PP}{{\mathbb P}}
\newcommand{\R}{{\mathbb R}}
\newcommand{\RR}{{\mathbb R}}
\newcommand{\E}{{\mathbb E}}
\newcommand{\AAA}{{\mathbb A}}

\newcommand{\Sym}{{\mathfrak{S}}}
\newcommand{\SymB}{{\mathfrak{B}}}
\newcommand{\Alt}{{\mathrm{Alt}}}

\newcommand{\germanA}{{\mathfrak{A}}}
\newcommand{\germanB}{{\mathfrak{B}}}
\newcommand{\germanQ}{{\mathfrak{Q}}}
\newcommand{\germanh}{{\mathfrak{h}}}

\newcommand{\myle}{\preceq}
\newcommand{\myge}{\succeq}
\newcommand{\mygt}{\succ}

\newcommand{\B}{{\sf B}}
\newcommand{\OB}{B^{\rm ord}}
\newcommand{\OS}{{\sf OS}}
\newcommand{\OO}{{\sf O}}
\newcommand{\OSP}{{\sf OSP}}
\newcommand{\Eu}{{\sf Eu}}
\newcommand{\ERR}{{\sf ERR}}
\newcommand{\sfB}{{\sf B}}
\newcommand{\sfD}{{\sf D}}
\newcommand{\sfE}{{\sf E}}
\newcommand{\sfG}{{\sf G}}
\newcommand{\sfJ}{{\sf J}}
\newcommand{\sfL}{{\sf L}}
\newcommand{\sfLhat}{{\widehat{{\sf L}}}}
\newcommand{\sfLcheck}{{\widecheck{{\sf L}}}}
\newcommand{\sfLtilde}{{\widetilde{{\sf L}}}}
\newcommand{\sfM}{{\sf M}}
\newcommand{\sfP}{{\sf P}}
\newcommand{\sfQ}{{\sf Q}}
\newcommand{\sfS}{{\sf S}}
\newcommand{\sfT}{{\sf T}}
\newcommand{\sfW}{{\sf W}}
\newcommand{\sfMV}{{\sf MV}}
\newcommand{\AMV}{{\sf AMV}}
\newcommand{\BM}{{\sf BM}}
\newcommand{\emIB}{B^{\rm irr}}
\newcommand{\emIP}{P^{\rm irr}}
\newcommand{\emOB}{B^{\rm ord}}
\newcommand{\emCB}{B^{\rm cyc}}
\newcommand{\emSC}{P^{\rm cyc}}

\newcommand{\cyc}{{\rm cyc}}
\newcommand{\e}{{\rm e}}
\newcommand{\ecyc}{{\rm ecyc}}
\newcommand{\epa}{{\rm epa}}
\newcommand{\iv}{{\rm iv}}
\newcommand{\pa}{{\rm pa}}
\newcommand{\pk}{{\rm p}}
\newcommand{\val}{{\rm v}}
\newcommand{\da}{{\rm da}}
\newcommand{\dd}{{\rm dd}}
\newcommand{\fp}{{\rm fp}}
\newcommand{\pkcyc}{{\rm pcyc}}
\newcommand{\valcyc}{{\rm vcyc}}
\newcommand{\dacyc}{{\rm dacyc}}
\newcommand{\ddcyc}{{\rm ddcyc}}
\newcommand{\pkpa}{{\rm ppa}}
\newcommand{\valpa}{{\rm vpa}}
\newcommand{\dapa}{{\rm dapa}}
\newcommand{\ddpa}{{\rm ddpa}}

\newcommand{\yp}{{y_\pk}}
\newcommand{\yptilde}{{\widetilde{y}_\pk}}
\newcommand{\yphat}{{\widehat{y}_\pk}}
\newcommand{\yv}{{y_\val}}
\newcommand{\yiv}{{y_\iv}}
\newcommand{\yda}{{y_\da}}
\newcommand{\ydatilde}{{\widetilde{y}_\da}}
\newcommand{\ydd}{{y_\dd}}
\newcommand{\yddtilde}{{\widetilde{y}_\dd}}
\newcommand{\yfp}{{y_\fp}}
\newcommand{\zp}{{z_\pk}}
\newcommand{\zv}{{z_\val}}
\newcommand{\zda}{{z_\da}}
\newcommand{\zdd}{{z_\dd}}

\newcommand{\lev}{{\rm lev}}
\newcommand{\stat}{{\rm stat}}
\newcommand{\umiv}{{\rm umiv}}
\newcommand{\umpa}{{\rm umpa}}
\newcommand{\mysteryone}{{\rm mys1}}
\newcommand{\mysterytwo}{{\rm mys2}}
\newcommand{\Asc}{{\rm Asc}}
\newcommand{\asc}{{\rm asc}}
\newcommand{\Des}{{\rm Des}}
\newcommand{\des}{{\rm des}}
\newcommand{\Exc}{{\rm Exc}}
\newcommand{\exc}{{\rm exc}}
\newcommand{\Wex}{{\rm Wex}}
\newcommand{\wex}{{\rm wex}}
\newcommand{\Fix}{{\rm Fix}}
\newcommand{\fix}{{\rm fix}}
\newcommand{\lrmax}{{\rm lrmax}}
\newcommand{\rlmax}{{\rm rlmax}}
\newcommand{\Rec}{{\rm Rec}}
\newcommand{\rec}{{\rm rec}}
\newcommand{\Arec}{{\rm Arec}}
\newcommand{\arec}{{\rm arec}}
\newcommand{\ERec}{{\rm ERec}}
\newcommand{\erec}{{\rm erec}}
\newcommand{\EArec}{{\rm EArec}}
\newcommand{\earec}{{\rm earec}}
\newcommand{\recarec}{{\rm recarec}}
\newcommand{\nonrec}{{\rm nonrec}}
\newcommand{\Cpeak}{{\rm Cpeak}}
\newcommand{\cpeak}{{\rm cpeak}}
\newcommand{\Cval}{{\rm Cval}}
\newcommand{\cval}{{\rm cval}}
\newcommand{\Cdasc}{{\rm Cdasc}}
\newcommand{\cdasc}{{\rm cdasc}}
\newcommand{\Cddes}{{\rm Cddes}}
\newcommand{\cddes}{{\rm cddes}}
\newcommand{\cdrise}{{\rm cdrise}}
\newcommand{\cdfall}{{\rm cdfall}}
\newcommand{\Peak}{{\rm Peak}}
\newcommand{\peak}{{\rm peak}}
\newcommand{\Val}{{\rm Val}}
\newcommand{\Dasc}{{\rm Dasc}}
\newcommand{\dasc}{{\rm dasc}}
\newcommand{\Ddes}{{\rm Ddes}}
\newcommand{\ddes}{{\rm ddes}}
\newcommand{\inv}{{\rm inv}}
\newcommand{\maj}{{\rm maj}}
\newcommand{\rs}{{\rm rs}}
\newcommand{\cross}{{\rm cr}}
\newcommand{\crosshat}{{\widehat{\rm cr}}}
\newcommand{\nest}{{\rm ne}}
\newcommand{\rodd}{{\rm rodd}}
\newcommand{\reven}{{\rm reven}}
\newcommand{\lodd}{{\rm lodd}}
\newcommand{\leven}{{\rm leven}}
\newcommand{\sg}{{\rm sg}}
\newcommand{\bl}{{\rm bl}}
\newcommand{\tran}{{\rm tr}}
\newcommand{\area}{{\rm area}}
\newcommand{\ret}{{\rm ret}}
\newcommand{\peaks}{{\rm peaks}}
\newcommand{\hl}{{\rm hl}}
\newcommand{\sll}{{\rm sl}}
\newcommand{\negg}{{\rm neg}}
\newcommand{\imp}{{\rm imp}}
\newcommand{\osg}{{\rm osg}}
\newcommand{\ons}{{\rm ons}}
\newcommand{\isg}{{\rm isg}}
\newcommand{\ins}{{\rm ins}}
\newcommand{\LL}{{\rm LL}}
\newcommand{\height}{{\rm ht}}
\newcommand{\as}{{\rm as}}

\newcommand{\ba}{{\bm{a}}}
\newcommand{\bahat}{{\widehat{\bm{a}}}}
\newcommand{\sfa}{{{\sf a}}}
\newcommand{\bb}{{\bm{b}}}
\newcommand{\bc}{{\bm{c}}}
\newcommand{\bchat}{{\widehat{\bm{c}}}}
\newcommand{\bd}{{\bm{d}}}
\newcommand{\bee}{{\bm{e}}}
\newcommand{\beh}{{\bm{eh}}}
\newcommand{\bff}{{\bm{f}}}
\newcommand{\bg}{{\bm{g}}}
\newcommand{\bh}{{\bm{h}}}
\newcommand{\bll}{{\bm{\ell}}}
\newcommand{\bp}{{\bm{p}}}
\newcommand{\bq}{{\bm{q}}}
\newcommand{\br}{{\bm{r}}}
\newcommand{\bs}{{\bm{s}}}
\newcommand{\bt}{{\bm{t}}}
\newcommand{\bu}{{\bm{u}}}
\newcommand{\bv}{{\bm{v}}}
\newcommand{\bw}{{\bm{w}}}
\newcommand{\bx}{{\bm{x}}}
\newcommand{\by}{{\bm{y}}}
\newcommand{\bz}{{\bm{z}}}
\newcommand{\bA}{{\bm{A}}}
\newcommand{\bB}{{\bm{B}}}
\newcommand{\bC}{{\bm{C}}}
\newcommand{\bE}{{\bm{E}}}
\newcommand{\bF}{{\bm{F}}}
\newcommand{\bG}{{\bm{G}}}
\newcommand{\bH}{{\bm{H}}}
\newcommand{\bI}{{\bm{I}}}
\newcommand{\bJ}{{\bm{J}}}
\newcommand{\bL}{{\bm{L}}}
\newcommand{\bLhat}{{\widehat{\bm{L}}}}
\newcommand{\bM}{{\bm{M}}}
\newcommand{\bN}{{\bm{N}}}
\newcommand{\bP}{{\bm{P}}}
\newcommand{\bQ}{{\bm{Q}}}
\newcommand{\bR}{{\bm{R}}}
\newcommand{\bS}{{\bm{S}}}
\newcommand{\bT}{{\bm{T}}}
\newcommand{\bW}{{\bm{W}}}
\newcommand{\bX}{{\bm{X}}}
\newcommand{\bY}{{\bm{Y}}}
\newcommand{\bIB}{{\bm{B}^{\rm irr}}}
\newcommand{\bOB}{{\bm{B}^{\rm ord}}}
\newcommand{\bOS}{{\bm{OS}}}
\newcommand{\bERR}{{\bm{ERR}}}
\newcommand{\bSP}{{\bm{SP}}}
\newcommand{\bMV}{{\bm{MV}}}
\newcommand{\bBM}{{\bm{BM}}}
\newcommand{\balpha}{{\bm{\alpha}}}
\newcommand{\balphapre}{{\bm{\alpha}^{\rm pre}}}
\newcommand{\bbeta}{{\bm{\beta}}}
\newcommand{\bgamma}{{\bm{\gamma}}}
\newcommand{\bdelta}{{\bm{\delta}}}
\newcommand{\bkappa}{{\bm{\kappa}}}
\newcommand{\bmu}{{\bm{\mu}}}
\newcommand{\bomega}{{\bm{\omega}}}
\newcommand{\bsigma}{{\bm{\sigma}}}
\newcommand{\btau}{{\bm{\tau}}}
\newcommand{\bphi}{{\bm{\phi}}}
\newcommand{\bphihat}{{\skew{3}\widehat{\vphantom{t}\protect\smash{\bm{\phi}}}}}
\newcommand{\bpsi}{{\bm{\psi}}}
\newcommand{\bxi}{{\bm{\xi}}}
\newcommand{\bzeta}{{\bm{\zeta}}}
\newcommand{\bone}{{\bm{1}}}
\newcommand{\bzero}{{\bm{0}}}

\newcommand{\Cbar}{{\overline{C}}}
\newcommand{\Dbar}{{\overline{D}}}
\newcommand{\dbar}{{\overline{d}}}
\newcommand{\Pbar}{{\bar{P}}}
\newcommand{\pbar}{{\bar{p}}}
\newcommand{\Lbar}{{\bar{L}}}
\newcommand{\Ubar}{{\bar{U}}}
\def\Btilde{{\widetilde{B}}}
\def\Ctilde{{\widetilde{C}}}
\def\Ftilde{{\widetilde{F}}}
\def\Gtilde{{\widetilde{G}}}
\def\Htilde{{\widetilde{H}}}
\def\Lhat{{\widehat{L}}}
\def\Ltilde{{\widetilde{L}}}
\def\Ptilde{{\widetilde{P}}}
\def\ptilde{{\widetilde{p}}}
\def\Chat{{\widehat{C}}}
\def\ctilde{{\widetilde{c}}}
\def\zbar{{\overline{Z}}}
\def\pitilde{{\widetilde{\pi}}}
\def\omegahat{{\widehat{\omega}}}

\def\CAsc{{\rm CAsc}}

\newcommand{\sech}{{\rm sech}}

%
%
\newcommand{\sn}{{\rm sn}}
\newcommand{\cn}{{\rm cn}}
\newcommand{\dn}{{\rm dn}}
\newcommand{\sm}{{\rm sm}}
\newcommand{\cm}{{\rm cm}}

%
%
\newcommand{\zfz}{ {{}_0 \! F_0} }
\newcommand{\zfo}{ {{}_0  F_1} }
\newcommand{\ofz}{ {{}_1 \! F_0} }
\newcommand{\ofo}{ {{}_1 \! F_1} }
\newcommand{\oft}{ {{}_1 \! F_2} }

%
%
\newcommand{\FHyper}[2]{ {\tensor[_{#1 \!}]{F}{_{#2}}\!} }
\newcommand{\FHYPER}[5]{ {\FHyper{#1}{#2} \!\biggl(
   \!\!\begin{array}{c} #3 \\[1mm] #4 \end{array}\! \bigg|\, #5 \! \biggr)} }
\newcommand{\tfo}{ {\FHyper{2}{1}} }
\newcommand{\tfz}{ {\FHyper{2}{0}} }
\newcommand{\threefz}{ {\FHyper{3}{0}} }
\newcommand{\FHYPERbottomzero}[3]{ {\FHyper{#1}{0} \hspace*{-0mm}\biggl(
   \!\!\begin{array}{c} #2 \\[1mm] \hbox{---} \end{array}\! \bigg|\, #3 \! \biggr)} }
\newcommand{\FHYPERtopzero}[3]{ {\FHyper{0}{#1} \hspace*{-0mm}\biggl(
   \!\!\begin{array}{c} \hbox{---} \\[1mm] #2 \end{array}\! \bigg|\, #3 \! \biggr)} }

\newcommand{\phiHyper}[2]{ {\tensor[_{#1}]{\phi}{_{#2}}} }
\newcommand{\psiHyper}[2]{ {\tensor[_{#1}]{\psi}{_{#2}}} }
\newcommand{\PhiHyper}[2]{ {\tensor[_{#1}]{\Phi}{_{#2}}} }
\newcommand{\PsiHyper}[2]{ {\tensor[_{#1}]{\Psi}{_{#2}}} }
\newcommand{\phiHYPER}[6]{ {\phiHyper{#1}{#2} \!\left(
   \!\!\begin{array}{c} #3 \\ #4 \end{array}\! ;\, #5, \, #6 \! \right)\!} }
\newcommand{\psiHYPER}[6]{ {\psiHyper{#1}{#2} \!\left(
   \!\!\begin{array}{c} #3 \\ #4 \end{array}\! ;\, #5, \, #6 \! \right)} }
\newcommand{\PhiHYPER}[5]{ {\PhiHyper{#1}{#2} \!\left(
   \!\!\begin{array}{c} #3 \\ #4 \end{array}\! ;\, #5 \! \right)\!} }
\newcommand{\PsiHYPER}[5]{ {\PsiHyper{#1}{#2} \!\left(
   \!\!\begin{array}{c} #3 \\ #4 \end{array}\! ;\, #5 \! \right)\!} }
\newcommand{\zerophizero}{ {\phiHyper{0}{0}} }
\newcommand{\ophizero}{ {\phiHyper{1}{0}} }
\newcommand{\zphio}{ {\phiHyper{0}{1}} }
\newcommand{\ophio}{ {\phiHyper{1}{1}} }
\newcommand{\tphio}{ {\phiHyper{2}{1}} }
\newcommand{\tphiz}{ {\phiHyper{2}{0}} }
\newcommand{\tPhio}{ {\PhiHyper{2}{1}} }
\newcommand{\opsio}{ {\psiHyper{1}{1}} }

%
%
\newcommand{\stirlingsubset}[2]{\genfrac{\{}{\}}{0pt}{}{#1}{#2}}
\newcommand{\stirlingcycle}[2]{\genfrac{[}{]}{0pt}{}{#1}{#2}}
\newcommand{\stirlingcyclesub}[2]{\widehat{\genfrac{\{}{\}}{0pt}{}{#1}{#2}}}
\newcommand{\assocstirlingsubset}[3]{{\genfrac{\{}{\}}{0pt}{}{#1}{#2}}_{\! \ge #3}}
\newcommand{\bigstirlingcycle}[2]{\left[\genfrac{}{}{0pt}{}{#1}{#2}\right]}
\newcommand{\genstirlingsubset}[4]{{\genfrac{\{}{\}}{0pt}{}{#1}{#2}}_{\! #3,#4}}
\newcommand{\irredstirlingsubset}[2]{{\genfrac{\{}{\}}{0pt}{}{#1}{#2}}^{\!\rm irr}}
\newcommand{\euler}[2]{\genfrac{\langle}{\rangle}{0pt}{}{#1}{#2}}
\newcommand{\eulergen}[3]{{\genfrac{\langle}{\rangle}{0pt}{}{#1}{#2}}_{\! #3}}
\newcommand{\eulersecond}[2]{\left\langle\!\! \euler{#1}{#2} \!\!\right\rangle}
\newcommand{\eulersecondgen}[3]{{\left\langle\!\! \euler{#1}{#2} \!\!\right\rangle}_{\! #3}}
\newcommand{\binomvert}[2]{\genfrac{\vert}{\vert}{0pt}{}{#1}{#2}}
\newcommand{\binomsquare}[2]{\genfrac{[}{]}{0pt}{}{#1}{#2}}
\newcommand{\doublebinom}[2]{\left(\!\! \binom{#1}{#2} \!\!\right)}
\newcommand{\lahnum}[2]{\genfrac{\lfloor}{\rfloor}{0pt}{}{#1}{#2}}
\newcommand{\rlahnum}[3]{{\genfrac{\lfloor}{\rfloor}{0pt}{}{#1}{#2}}_{\! #3}}

%
%
\newcommand{\Lna}{{L_n^{(\alpha)}}}
\newcommand{\scrlna}{{\scrl_n^{(\alpha)}}}
\newcommand{\scrlbarna}{{\scrlbar_n^{(\alpha)}}}
\newcommand{\scrlhatna}{{\widehat{\scrl}_n^{(\alpha)}}}
\newcommand{\scrlhatbarna}{{\overline{\widehat{\scrl}}_n^{(\alpha)}}}
\newcommand{\scrlbarnka}{{\scrlbar_{n,k}^{(\alpha)}}}
\newcommand{\scrlbarnzeroa}{{\scrlbar_{n,0}^{(\alpha)}}}
\newcommand{\Lnb}{{L_n^{[\beta]}}}
\newcommand{\scrlnb}{{\scrl_n^{[\beta]}}}
\newcommand{\scrlbarnb}{{\scrlbar_n^{[\beta]}}}
\newcommand{\Lah}{{\textrm{Lah}}}
\newcommand{\LD}{{\mathbf{LD}}}
\newcommand{\MLD}{{\mathbf{MLD}}}


\newenvironment{sarray}{
             \textfont0=\scriptfont0
             \scriptfont0=\scriptscriptfont0
             \textfont1=\scriptfont1
             \scriptfont1=\scriptscriptfont1
             \textfont2=\scriptfont2
             \scriptfont2=\scriptscriptfont2
             \textfont3=\scriptfont3
             \scriptfont3=\scriptscriptfont3
           \renewcommand{\arraystretch}{0.7}
           \begin{array}{l}}{\end{array}}

\newenvironment{scarray}{
             \textfont0=\scriptfont0
             \scriptfont0=\scriptscriptfont0
             \textfont1=\scriptfont1
             \scriptfont1=\scriptscriptfont1
             \textfont2=\scriptfont2
             \scriptfont2=\scriptscriptfont2
             \textfont3=\scriptfont3
             \scriptfont3=\scriptscriptfont3
           \renewcommand{\arraystretch}{0.7}
           \begin{array}{c}}{\end{array}}


\newcommand*\circled[1]{\tikz[baseline=(char.base)]{
  \node[shape=circle,draw,inner sep=1pt] (char) {#1};}}
\newcommand{\ostar}{{\circledast}}
\newcommand{\ostarN}{{\,\circledast_{\vphantom{\dot{N}}N}\,}}
\newcommand{\ostarPsi}{{\,\circledast_{\vphantom{\dot{\Psi}}\Psi}\,}}
\newcommand{\starN}{{\,\ast_{\vphantom{\dot{N}}N}\,}}
\newcommand{\starpsi}{{\,\ast_{\vphantom{\dot{\bpsi}}\!\bpsi}\,}}
\newcommand{\starone}{{\,\ast_{\vphantom{\dot{1}}1}\,}}
\newcommand{\startwo}{{\,\ast_{\vphantom{\dot{2}}2}\,}}
\newcommand{\starinfty}{{\,\ast_{\vphantom{\dot{\infty}}\infty}\,}}
\newcommand{\starT}{{\,\ast_{\vphantom{\dot{T}}T}\,}}

\newcommand*{\Scale}[2][4]{\scalebox{#1}{$#2$}}

\newcommand*{\Scaletext}[2][4]{\scalebox{#1}{#2}} 

\newcommand{\bolddot}{\boldsymbol{\cdot}}

\clearpage 

\tableofcontents

\clearpage

\section{Introduction and statement of main results}
\label{sec.intro}

A (finite or infinite) matrix of real numbers is called
\textbfit{totally positive} (TP) if all its minors are nonnegative.
A matrix of polynomials with real coefficients is called
\textbfit{coefficientwise totally positive}
if all its minors are polynomials with nonnegative coefficients.
Let $A = (a_{n,k})_{n,k\geq 0}$ be a lower-triangular matrix
of nonnegative real numbers. 
One can then ask the following four questions:
\begin{enumerate}
   \item Is $A$ totally positive?
   \item Is the lower-triangular matrix
$A^{\rm rev} \eqdef (a_{n,n-k})_{n,k\geq 0}$,
obtained by reversing the rows of $A$, totally positive?
(Here $a_{n,n-k} \eqdef 0$ when $n<k$.)
   \item Let $A_n(x) = \sum_{k=0}^n a_{n,k} x^k$ denote
the row-generating polynomial of the $n$-th row of $A$.
Are the polynomials $A_n(x)$ negative-real-rooted?\footnote{
      A polynomial with real coefficients is said to be \textbfit{real-rooted} 
      if it is either identically zero or all its (complex) zeros are real;
      it is said to be \textbfit{negative-real-rooted}
      if it is either identically zero or all its (complex) zeros
      lie in $(-\infty,0]$.
      Since we are here assuming that all the matrix entries $a_{n,k}$
      are nonnegative, clearly $A_n(x)$ cannot have any
      {\em positive}\/ real zeros unless it is identically zero.
      So negative-real-rootedness and real-rootedness are equivalent here.
}
This is equivalent \cite[pp.~395, 399]{Karlin_68}
to asking if the Toeplitz matrix 
of the $n$-th row sequence $(a_{n,k})_{k\geq 0}$
--- that is, the matrix $(a_{n,i-j})_{i,j \ge 0}$
with $a_{n,k} \eqdef 0$ for $k < 0$ ---
is totally positive.
In this situation we say that the $n$-th row sequence $(a_{n,k})_{k\geq 0}$
is \textbfit{Toeplitz-totally positive}.\footnote{
      Sequences whose Toeplitz matrices are totally positive
      are also known as \textbfit{P\'olya frequency sequences}.
      The just-cited characterisation of {\em finitely supported}\/
      P\'olya frequency sequences is a special case of the
      (much more difficult) theorem of
      Aissen, Edrei, Schoenberg and Whitney \cite{Aissen_52,Edrei_52}
      that characterises arbitrary P\'olya frequency sequences
      \cite[p.~412, Theorem~5.3]{Karlin_68}.
}
   \item Is the Hankel matrix of the polynomial sequence $(A_n(x))_{n\geq 0}$
--- that is, the matrix $(A_{i+j}(x))_{i,j \ge 0}$ ---
coefficientwise totally positive in the variable $x$?
In this situation we say that the sequence $(A_n(x))_{n\geq 0}$
is \textbfit{coefficientwise Hankel-totally positive}.
One sufficient but far-from-necessary condition for this to hold
is for the ordinary generating function $\sum_{n=0}^\infty A_n(x) t^n$
to have a Stieltjes-type continued fraction 
(or more generally, a Thron-type continued fraction)
with coefficients that are polynomials with nonnegative coefficients
in the variable $x$.\footnote{
   This result for Stieltjes-type continued fractions (S-fractions)
   was first stated explicitly in \cite{Sokal_flajolet};
   see \cite[Theorem~9.8]{latpath_SRTR} for a published proof
   (in fact of a much more general result).
   However, this result is in fact
   an easy corollary of Flajolet's \cite{Flajolet_80} 
   combinatorial interpretation of Stieltjes-type continued fractions
   together with Viennot's \cite[Chapter~4, Section~3]{Viennot_83}
   combinatorial characterization of the Hankel minors.
   The generalization to Thron-type continued fractions (T-fractions)
   is given (again in much greater generality)
   in \cite[Theorem~9.9]{latpath_SRTR}.
}
\end{enumerate}

For several important combinatorial triangles,
the answers to all four of these questions
have either been proven to be yes or are conjectured to be yes.
Some examples are:
\begin{enumerate}
   \item When $A = \left(\binom{n}{k}\right)_{n,k\geq 0}$
is the matrix of binomial coefficients.
The binomial matrix $A$ is well known to be totally positive,
and there are many proofs of this fact:
see e.g.~\cite[pp.~232--233]{Aigner_18} for a proof using the 
Lindstr\"om--Gessel--Viennot lemma.
Here $A^{\rm rev} = A$ as the rows are palindromic.
The row-generating polynomials $A_n(x) = (1+x)^n$ 
are clearly negative-real-rooted.
And finally, the sequence $\big( (1+x)^n \big)_{n\geq 0}$
is coefficientwise Hankel-totally positive in the variable $x$:
this can easily be seen directly,
or it follows from the continued fraction 
$\sum\limits_{n=0}^\infty (1+x)^n \, t^n = 1/[1-(1+x)t]$.
   \item When $A = \left(\stirlingcycle{n}{k}\right)_{n,k\geq 0}$
is the matrix of Stirling cycle numbers,
which count the permutations of $[n] \eqdef \{1,\ldots,n\}$ with $k$ cycles.
The total positivity of the matrix $A$ and its reversal $A^{\rm rev}$ 
are due to Brenti \cite{Brenti_95} in 1995:
they hold because the entries $\stirlingcycle{n}{k}$ satisfy a
binomial-like recurrence with nonnegative $n$-dependent coefficients.
The row-generating polynomials $A_n(x) = x(x+1)\cdots (x+n-1)$
are clearly negative-real-rooted.
The sequence $\big( A_n(x) \big)_{n\geq 0}$
has a Stieltjes-type continued fraction,
first discovered by Euler \cite{Euler_1760} in 1760:
see e.g. \cite[eqn.~(1.2)]{Sokal_alg_contfrac}
\cite[eqns.~(2.2)/(2.6)]{Sokal-Zeng_masterpoly}.
This implies the coefficientwise Hankel-total positivity.
   \item When $A = \left(\stirlingsubset{n}{k}\right)_{n,k\geq 0}$
is the matrix of Stirling subset numbers,
which count the set partitions of $[n]$ with $k$ blocks.
The total positivity of the matrix $A$ 
is due to Brenti \cite{Brenti_95} in 1995:
it holds because the entries $\stirlingsubset{n}{k}$ satisfy a
binomial-like recurrence with nonnegative $k$-dependent coefficients.
However, the total positivity of the reversal $A^{\rm rev}$ 
was only recently shown by us in \cite{Chen_TPrecurrences}.
The real-rootedness of the row-generating polynomials $A_n(x)$
was proven by Harper \cite{Harper_67} in 1967.
The sequence $\big( A_n(x) \big)_{n\geq 0}$
has a Stieltjes-type continued fraction,
discovered by various people in the second half of the twentieth century:
see e.g. \cite[eqns.~(3.2)/(3.5) and footnote~19]{Sokal-Zeng_masterpoly}.
This implies the coefficientwise Hankel-total positivity.
   \item When $A = \left(\euler{n}{k}\right)_{n,k\geq 0}$,
is the matrix of Eulerian numbers, 
which count the permutations on $[n]$ with $k$ descents (or $k$ excedances).
The row-generating polynomials $A_n(x)$ were shown to be real-rooted 
by Frobenius \cite{Frobenius_1910} in 1910.
The sequence $\big( A_n(x) \big)_{n\geq 0}$
has a Stieltjes-type continued fraction
going back to Stieltjes \cite{Stieltjes_1894}:
see e.g. \cite[eqns.~(2.2)/(2.8) and footnote~4]{Sokal-Zeng_masterpoly}.
This implies the coefficientwise Hankel-total positivity.
In 1996, Brenti \cite{Brenti_96} conjectured
that the matrix $A$ is totally positive;
this conjecture is, alas, still open.\footnote{
   Dyachenko \cite{Dyachenko_private_Eulerian}
   has verified this conjecture for the leading $512\times 512$ submatrix
   of the matrix $\left(\euler{n+1}{k}\right)_{n,k\geq 0}$,
   by computing the bidiagonal factorization
   using Neville elimination \cite{Gasca_92}.
}
Here the rows of $A$ are almost-palindromic
in the sense that $\euler{n}{k} = \euler{n}{n-1-k}$ for $n \ge 1$;
it~easily follows from this that $A^{\rm rev}$ is totally positive
if and only if $A$ is.
In \cite{Chen_TPrecurrences} we have generalized Brenti's conjecture:
we conjecture the {\em coefficientwise}\/ total positivity
of a lower-triangular matrix in six indeterminates
that is defined by a particular recurrence;
this matrix includes the Eulerian matrix,
the Stirling subset matrix and the reversed Stirling subset matrix
as special cases.
\end{enumerate}

\medskip

In this paper, we introduce two infinite families of lower-triangular matrices
generalizing the Stirling cycle and subset triangles,
parametrized by an integer $r \ge 1$:
we call them the
\textbfit{$\bm{r^{th}}$-order Stirling cycle and subset triangles}.
For $r=1$ they reduce to the ordinary Stirling cycle and subset triangles.
The central theme of this paper is to investigate the foregoing four questions
for these triangular arrays.
We have combined an experimental approach
with proofs using both bijective and generating-function methods.
Even though we are able to prove only very few of our conjectures
--- notably those for the case $r=2$ ---
we have proven several interesting combinatorial interpretations
for the $r$th-order Stirling triangles,
and our results suggest further interesting consequences
in the form of combinatorial interpretations and additional conjectures.

The rest of this introduction is structured as follows:
In Section~\ref{sec.def.higherorder}
we define the $r$th-order Stirling cycle and subset numbers
and state some recurrences that these numbers satisfy;
we also set up some basic notation for this paper.
Then in Section~\ref{sec.def.eulerian} we introduce
the second-order Eulerian numbers
and relate them to the second-order Stirling cycle and subset numbers; 
this connection will play a major role in our proofs.
We will also introduce two generalisations:
the higher-order Eulerian numbers and the higher-order quasi-Eulerian numbers.
In Section~\ref{sec.conjecture} we state our main conjectures
on the higher-order Stirling cycle and subset numbers,
and also state our theorems.
In Section~\ref{sec.conj.eulerian} we state some conjectures for the 
higher-order Eulerian and quasi-Eulerian numbers.
In Section~\ref{sec.structure} we present the outline of the remainder
of this paper.

\subsection{Higher-order Stirling cycle and subset numbers}
\label{sec.def.higherorder}

Several generalizations and extensions of the
Stirling cycle and subset numbers have been introduced in the literature.
One such family of generalizations
are the {\em $r$-associated Stirling cycle and subset numbers}\/,
defined as follows:

The \textbfit{$\bm{r}$-associated Stirling cycle number}
$\stirlingcycle{n}{k}_r$
is the number of permutations on $[n]$ with $k$ cycles,
all of which are of size $\ge r$
\cite[pp.~256--257, 295]{Comtet_74}. 
Thus we have $\stirlingcycle{n}{k}_1 = \stirlingcycle{n}{k}$.
Note that $\stirlingcycle{n}{k}_r$
is nonzero only when $n\geq rk$, 
or equivalently $k\leq \lfloor n/r \rfloor$.
See \cite[A132393/A106828/A050211/A050212/A050213]{OEIS}
for $1 \le r \le 5$.

Similarly, the \textbfit{$\bm{r}$-associated Stirling subset number}
$\stirlingsubset{n}{k}_r$
is the number of set partitions of $[n]$ into $k$ blocks,
all of which are of size $\ge r$
\cite[pp.~221--222]{Comtet_74}.
Thus we have $\stirlingsubset{n}{k}_1 = \stirlingsubset{n}{k}$.
Note that $\stirlingsubset{n}{k}_r$
is nonzero only when $n\geq rk$,
or equivalently $k\leq \lfloor n/r \rfloor$.
See \cite[A048993/A008299/A059022/A059023/A059024/A059025]{OEIS}
for $1 \le r \le 6$.

It is not difficult to show, by considering the possibilities for
the element $n$ in the permutation or set partition,
that the $r$-associated Stirling cycle and subset numbers
satisfy the recurrences \cite[pp.~222, 257]{Comtet_74}
\begin{eqnarray}
   \stirlingcycle{n}{k}_r
   & = &
   (r-1)! \, \binom{n-1}{r-1} \, \stirlingcycle{n-r}{k-1}_r
        \:+\: (n-1) \, \stirlingcycle{n-1}{k}_r
   \label{eq.recurrence.rassoc.cycle}  \\[3mm]
   \stirlingsubset{n}{k}_r
   & = &
   \binom{n-1}{r-1} \, \stirlingcycle{n-r}{k-1}_r
       \:+\: k \, \stirlingsubset{n-1}{k}_r
   \label{eq.recurrence.rassoc.subset}
\end{eqnarray}

For $r\geq 2$, the sequences of row-generating polynomials of these triangles
are {\em not}\/ coefficientwise Hankel-totally positive. 
This is because the $n$-th row-generating polynomial is of degree 
$\lfloor n/r\rfloor$ --- not $n$ ---
and it follows from this that even the $2\times 2$ Hankel minors
must have some negative coefficients.

Instead, we modify these triangles by shifting the columns up
until we get a lower-triangular matrix with nonzero diagonal entries.
That is, we define, for each $r \ge 1$,
\begin{eqnarray}
\stirlingcycle{n}{k}^{(r)} & \eqdef & \stirlingcycle{n+(r-1)k}{k}_r   
	\label{eq.stircyc.def.new}  \\[2mm]
\stirlingsubset{n}{k}^{(r)} & \eqdef & \stirlingsubset{n+(r-1)k}{k}_r
	\label{eq.stirsub.def.new}
\end{eqnarray}
We call these the
\textbfit{$\bm{r^{th}}$-order Stirling cycle and subset numbers},
respectively;
and we call the lower-triangular matrices
\begin{subeqnarray}
   C^{(r)} & \eqdef & \left(\stirlingcycle{n}{k}^{(r)}\right)_{n,k\geq 0}
       \\[2mm]
   S^{(r)} & \eqdef & \left(\stirlingsubset{n}{k}^{(r)}\right)_{n,k\geq 0}\;.
\end{subeqnarray}
the \textbfit{$\bm{r^{th}}$-order Stirling cycle and subset triangles}.
Of course we have
\begin{eqnarray}
\stirlingcycle{n}{k}^{(1)} & = &\stirlingcycle{n}{k}_1 \;=\; \stirlingcycle{n}{k}  \\[2mm]
\stirlingsubset{n}{k}^{(1)} & = &\stirlingsubset{n}{k}_1 \;=\; \stirlingsubset{n}{k}
\end{eqnarray}
so that the numbers of order $r=1$
are the usual Stirling cycle and subset numbers, respectively.
For the convenience of the reader and for our use in the rest of the paper, 
we rewrite \reff{eq.stircyc.def.new}/\reff{eq.stirsub.def.new} for $r=2$: 
\begin{eqnarray}
	\stirlingcycle{n}{k}^{(2)} & = & \stirlingcycle{n+k}{k}_2  \\[2mm]
	\stirlingsubset{n}{k}^{(2)} & = & \stirlingsubset{n+k}{k}_2
\end{eqnarray}
The triangle $C^{(2)}$ is \cite[A259456/A269940/A111999]{OEIS}
and its row sums are \cite[A032188]{OEIS}.
The triangle $S^{(2)}$ is \cite[A134991/A269939]{OEIS}
and its row sums are \cite[A000311]{OEIS}. 
The second-order Stirling subset numbers are often referred 
to as the \textbfit{Ward numbers}
\cite{Ward_34,Carlitz_71,Clark_99,Barbero_15,Elvey-Price-Sokal_wardpoly}.
The row sums of $C^{(r)}$ and $S^{(r)}$ for $r \ge 3$
are not currently in the OEIS \cite{OEIS}.
In the Appendix we provide tables of the triangles
$C^{(r)}$ and $S^{(r)}$ for $1 \le r \le 4$
and the first few values of $n$.

From the combinatorial definition one sees that the diagonal elements
of the arrays $C^{(r)}$ and $S^{(r)}$ are
\begin{eqnarray}
\stirlingcycle{n}{n}^{(r)} & = & \stirlingcycle{rn}{n}_r
     \;=\;  {(rn)! \over r^n \, n!}
        \label{eq.stircyc.diagonal}  \\[2mm]
\stirlingsubset{n}{n}^{(r)} & = & \stirlingsubset{rn}{n}_r
     \;=\;  {(rn)! \over (r!)^n \, n!}
        \label{eq.stirsub.diagonal}
\end{eqnarray}
\noindent
\hspace*{-1.5mm}
\cite[A001147/A052502/A060706/A052504]{OEIS}
\cite[A001147/A025035/A025036/A025037/ A025038]{OEIS}.
In particular, for $r=2$ we have
\be
   \stirlingcycle{n}{n}^{(2)}  \;=\;  \stirlingsubset{n}{n}^{(2)}
        \;=\;  (2n-1)!!
   \;.
 \label{eq.stir.diagonal.r=2}
\ee
The numbers $(2n-1)!!$ \cite{Callan_09} \cite[A001147]{OEIS}
will play a central role in this paper:
see Section~\ref{subsec.remark} for an observation about this.

The recurrence relations
\reff{eq.recurrence.rassoc.cycle}/\reff{eq.recurrence.rassoc.subset}
for the $r$-associated Stirling cycle and subset numbers
easily imply, by change of variables,
the following recurrence relations
for the $r$th-order Stirling cycle and subset numbers:

\begin{lemma}
\begin{enumerate}
   \item[(a)] The numbers $\stirlingcycle{n}{k}^{(r)}$
satisfy the recurrence
\be
   \stirlingcycle{n}{k}^{(r)} 
	\;=\;
   (r-1)! \, \binom{n+(r-1)k-1}{r-1} \, \stirlingcycle{n-1}{k-1}^{(r)} 
       \:+\:
       \big(n+(r-1)k-1\big)\, \stirlingcycle{n-1}{k}^{(r)}
   \label{rec.cycle}
\ee
for $n \ge 1$, with initial conditions
$\stirlingcycle{0}{k}^{(r)} = \delta_{k0}$
and $\stirlingcycle{n}{-1}^{(r)} = 0$.
   \item[(b)] The numbers $\stirlingsubset{n}{k}^{(r)}$
satisfy the recurrence
\be
   \stirlingsubset{n}{k}^{(r)}
   \;=\;
   \binom{n+(r-1)k-1}{r-1} \, \stirlingsubset{n-1}{k-1}^{(r)} 
         \:+\:  k\, \stirlingsubset{n-1}{k}^{(r)}
   \label{rec.subset}
\ee
for $n \ge 1$, with initial conditions
$\stirlingsubset{0}{k}^{(r)} = \delta_{k0}$
and $\stirlingsubset{n}{-1}^{(r)} = 0$.
\end{enumerate}
	\label{lem.rec}
\end{lemma}

When $r=2$, the recurrences~\reff{rec.cycle}/\reff{rec.subset} are
\begin{eqnarray}
	\stirlingcycle{n}{k}^{(2)} 
	&=&
	(n+k-1) \, \stirlingcycle{n-1}{k-1}^{(2)}
	\:+\:
	(n+k-1) \,  \stirlingcycle{n-1}{k}^{(2)}
	\label{eq.rec.stirlingcycle.2nd} 	\\[3mm]
	\stirlingsubset{n}{k}^{(2)} 
	&=&
	(n+k-1) \, \stirlingsubset{n-1}{k-1}^{(2)}
	\:+\: 
	k  \ \stirlingsubset{n-1}{k}^{(2)}
	\label{eq.rec.stirlingsubset.2nd} 
\end{eqnarray}
Note that the coefficients here are affine in $n$ and $k$ 
(this is no longer the case for $r\geq 3$).
It is curious that we see here the combination $n+k$,
in contrast with the $n-k$ occurring in
\cite[eqs.~(1.2)/(1.4)/(1.5)]{Chen_TPrecurrences}.

We write $c_{r,n}(x)$ for the row-generating polynomial
of the $n$-th row of the matrix $C^{(r)}$;
we call them the \textbfit{$\bm{r^{th}}$-order Stirling cycle polynomials}.
Similarly, we write $s_{r,n}(x)$ for the row-generating polynomial
of the $n$-th row of the matrix $S^{(r)}$;
we call them the \textbfit{$\bm{r^{th}}$-order Stirling subset polynomials}.
Thus, we have
\begin{eqnarray}
   c_{r,n}(x) & \eqdef & \sum_{k=0}^{n} \stirlingcycle{n}{k}^{(r)} x^k
      \label{eq.cycle.base}  \\[2mm]
   s_{r,n}(x) & \eqdef & \sum_{k=0}^{n} \stirlingsubset{n}{k}^{(r)} x^k
      \label{eq.subset.base}
\end{eqnarray}
Finally, we use $\widecheck{C}^{(r)}$ and $\widecheck{S}^{(r)}$
to denote the two row-reversed matrices 
$\widecheck{C}^{(r)} \eqdef \bigl(C^{(r)}\bigr)^{\rm rev}$
and 
$\widecheck{S}^{(r)} \eqdef \bigl(S^{(r)}\bigr)^{\rm rev}$.

\subsection{Relation with the second-order Eulerian numbers}
\label{sec.def.eulerian}

The second-order ($r=2$) Stirling cycle and subset numbers
are closely related to the second-order Eulerian numbers,
as we now explain.

Let $w = w_1 \cdots w_N$ be a word on a totally ordered alphabet $\AAA$.
A pair $(w_i, w_{i+1})$ with $i \in [N-1]$ is called an
\textbfit{ascent} (resp.~\textbfit{descent}, \textbfit{plateau})
in case $w_i < w_{i+1}$ (resp.~$w_i > w_{i+1}$, $w_i = w_{i+1}$).\footnote{
   Note that we do not allow $i=N$ as a position for an ascent, descent
   or plateau.  Some authors do allow this, by fixing a ``boundary condition''
   for $w_{N+1}$ (usually defining it to be either $<$ every element of $\AAA$
   or $>$ every element of $\AAA$).  It is important, when comparing formulae,
   to check carefully which convention the author is using.
}
The word $w$ is said to be a \textbfit{Stirling word} if,
for every letter $a \in \AAA$,
all letters between two successive occurrences of $a$ are larger than~$a$.
A permutation of a multiset $M$
is said to be a \textbfit{Stirling permutation}
if it is a Stirling word.
The second-order Eulerian number $\eulersecond{n}{k}$
is then defined to be the number of Stirling permutations
of the multiset $\{1,1,2,2,\ldots,n,n\}$ with $k$~ascents
(or $k$~descents: they are equidistributed,
 since reversing the word exchanges ascents with descents).\footnote{
   Here we follow the definition of
   Graham--Knuth--Patashnik \cite[p.~270]{Graham_94},
   which entails that, when $n \ge 1$,
   $\eulersecond{n}{k}$ is nonzero for $0 \le k \le n-1$.
   Many other authors \cite{Gessel_78,Park_94a,Janson_11,Dzhumadildaev_14}
   use the convention that the final index $N = 2n$
   is always the beginning of a descent
   (e.g.~by imposing the boundary condition $w_{N+1} = 0$),
   so that, when $n \ge 1$,
   $\eulersecond{n}{k}$ is nonzero for $1 \le k \le n$.
}
They satisfy the recurrence \cite[eq.~(6.41)]{Graham_94}
\be
   \eulersecond{n}{k}   \;=\;
      (2n-k-1) \, \eulersecond{n-1}{k-1}  \:+\:
      (k+1) \, \eulersecond{n-1}{k}
 \label{def.eulersecond}
\ee
with initial conditions
$\eulersecond{0}{k} = \delta_{k0}$ and $\eulersecond{n}{-1} = 0$.

Now define the \textbfit{second-order Eulerian triangle}
$E^{(2)} = \big( \eulersecond{n}{k} \big)_{n,k \ge 0}$
and the \textbfit{shifted-reversed second-order Eulerian triangle}
$\widecheck{E}^{(2)} = \big( \eulersecond{n}{n-k-1} \big)_{n,k \ge 0}$
[with the modification for $n=0$ to make the row $(1,0,0,\ldots)$].
We will prove:

\begin{proposition}[Relation of second-order Stirling cycle and subset numbers
   with second-order Eulerian numbers]
   \label{prop.secondeuler}
Let $B = \big(\, \binom{n}{k} \,\big)_{n,k \ge 0}$ be the binomial matrix.
Then:
\begin{itemize}
   \item[(a)] $\widecheck{C}^{(2)} \:=\: E^{(2)} \, B$.
   \item[(b)] $\widecheck{S}^{(2)} \:=\: \widecheck{E}^{(2)} \, B$.
\end{itemize}
\end{proposition}

\noindent
Part~(b) seems to have been known for a long time
\cite[eq.~(6)]{Smiley_00} \cite[eqns.~(1a)/(31b)]{Barbero_15};
and part~(a) is implicit in some formulas at
\cite[A106828/A358622]{OEIS}.
We will give several proofs of these equalities:
a simple algebraic proof based on the recurrences,
which however is not very enlightening;
and several (more illuminating) combinatorial proofs based on bijections
to other combinatorial objects.

If one defines the \textbfit{second-order Eulerian polynomials}
\be
   E^{(2)}_n(x)  \;\eqdef\;  \sum_{k=0}^{n} \eulersecond{n}{k} x^k
   \;,
\ee
then the identities of Proposition~\ref{prop.secondeuler}
imply a relation with the second-order Stirling cycle and subset polynomials:
\begin{eqnarray}
   x^n \, c_{2,n}(1/x)  & = &  E^{(2)}_n(1+x)
      \label{eq.cycle.secondeuler}  \\[2mm]
   x^n \, s_{2,n}(1/x)  & = &  (1+x)^{n-1} \, E^{(2)}_n\Big( {1 \over 1+x} \Big)
                \quad\hbox{for $n \ge 1$}
      \label{eq.subset.secondeuler}
\end{eqnarray}
or equivalently\footnote{
   The identity \reff{eq.subset.secondeuler.bis} can also be found in
   \cite[eq.~(1.15b)]{Elvey-Price-Sokal_wardpoly}:
   that paper used a different definition of the
   second-order Eulerian numbers, so that their $E^{[2]}_n(x)$
   equals our $x E^{(2)}_n(x)$ for $n \ge 1$.
}
\begin{eqnarray}
   c_{2,n}(x)  & = &  x^n \, E^{(2)}_n\Big( {1+x \over x} \Big)
      \label{eq.cycle.secondeuler.bis}  \\[2mm]
   s_{2,n}(x)  & = &  x \, (1+x)^{n-1} \, E^{(2)}_n\Big( {x \over 1+x} \Big)
                \quad\hbox{for $n \ge 1$}
      \label{eq.subset.secondeuler.bis}
\end{eqnarray}

Let us remark that one can define, more generally,
the $r$th-order Eulerian numbers $\euler{n}{k}^{\! (r)}$
as the number of Stirling permutations of the multiset $\{1^r,\ldots,n^r\}$
with $k$~ascents.
They satisfy the recurrence
\be
   \euler{n}{k}^{\! (r)}   \;=\;
      \big[ rn - k - (r-1) \big] \, \euler{n-1}{k-1}^{\! (r)}   \:+\:
      (k+1) \,  \euler{n-1}{k}^{\! (r)}
   \;.
 \label{eq.eulerian.recurrence}
\ee
We analogously define $r$-th order Eulerian triangle 
$E^{(r)} = \big(\euler{n}{k}^{\! (r)}\big)_{n,k\geq 0}$
and shifted-reversed $r$-th order Eulerian triangle
$\widecheck{E}^{(r)} = \big(\euler{n}{n-k-1}^{\! (r)}\big)_{n,k\geq 0}$
for $r\geq 1$.

However, as far as we know there is no relation between
the $r$th-order Stirling cycle or subset numbers
and the $r$th-order Eulerian numbers, for $r \ge 3$.
Indeed, this lack of relation is to be expected, since the recurrence
\reff{eq.eulerian.recurrence} is affine in $n$ and $k$ for all~$r$,
while the recurrences \reff{rec.cycle}/\reff{rec.subset}
are not affine in $n$ or $k$ when $r \ge 3$.
Instead, the matrices which do seem to be analogous in this context for 
$r\ge 3$
are $\widecheck{C}^{(r)} \, B^{-1}$ and $\widecheck{S}^{(r)} \, B^{-1}$.
We call their entries the
\textbfit{$\bm{r^{th}}$-order quasi-Eulerian numbers}.
We will further look at these matrices in Sections~\ref{sec.conj.eulerian}
and \ref{sec.interpretation.stirperm.gen}:
we will state some conjectures related to these numbers in Section~\ref{sec.conj.eulerian},
and we will provide a combinatorial interpretation for the 
entries of $\widecheck{C}^{(r)} \, B^{-1}$
in Section~\ref{sec.interpretation.stirperm.gen}.

\subsection{Statement of main conjectures and results}
\label{sec.conjecture}

Having defined the Stirling cycle and subset numbers of all orders $r$,
we are now ready to state our conjectures,
which constitute our first principal contribution in this paper.
We will then state which cases are known and which cases
we will prove in the present paper.

Our conjecture for the Stirling cycle numbers of all orders is the following:

\begin{conjecture}[Conjecture for Stirling cycle numbers of all orders]
   \label{conj.cycle}
The following are true:
\begin{itemize}
\item[(a)] The triangle $C^{(r)}$ is totally positive for all $r\geq 1$.

\item[(b)] The triangle $\widecheck{C}^{(r)}\eqdef \left(C^{(r)}\right)^{\rm rev}$
formed by reversing the rows of $C^{(r)}$
is totally positive for $r=1,2$
and is not totally positive when $r\geq 3$.

\item[(c)] The row-generating polynomials $c_{r,n}(x)$ are negative-real-rooted 
for $r=1,2$ 
and they have non-real complex zeros for $r\geq 3$
and $n\geq 3$.
Also, the row sequences of the matrix $C^{(r)}$ are log-concave
for $1\leq r \leq 5$.

\item[(d)]  The sequence $(c_{r,n}(x))_{n\geq 0}$
is coefficientwise Hankel-totally positive in the variable $x$
for all $r\geq 1$.
\end{itemize}
\end{conjecture}

As mentioned earlier, 
the $r=1$ cases of Conjecture~\ref{conj.cycle} are all known to be true.
In this paper we will prove (c) and (d) for $r=2$.\footnote{
  The log-concavity problems in (c) for $r=3,4,5$ were proposed
  by the first author at the 30th British Combinatorial Conference:
  see \cite{Cameron_24} for the full list of problems.
}

\bigskip

Our conjecture for the Stirling subset numbers of all orders is the following:

\begin{conjecture}[Conjecture for Stirling subset numbers of all orders]
   \label{conj.subset}
The following are true:
\begin{itemize} 
\item[(a)] The triangle $S^{(r)}$ is totally positive for all $r\geq 1$.

\item[(b)] The triangle $\widecheck{S}^{(r)}\eqdef \left(S^{(r)}\right)^{\rm rev}$ 
formed by reversing the rows of $S^{(r)}$
is totally positive for $r=1,2$
and is not totally positive when $r\geq 3$.

\item[(c)] The row-generating polynomials $s_{r,n}(x)$ are negative-real-rooted
for $r=1,2$ and they have non-real complex zeros for $r\geq 3$
and $n\geq 3$.
Also, the row sequences of the matrix $S^{(r)}$
are log-concave for $1\leq r \leq 5$.

\item[(d)]  The sequence $(s_{r,n}(x))_{n\geq 0}$ 
is coefficientwise Hankel-totally positive 
in the variable $x$ for $r=1,2$
and is not coefficientwise Hankel-totally positive for $r\geq 3$.
\end{itemize}
\end{conjecture}

As mentioned earlier,
the $r=1$ cases of Conjecture~\ref{conj.subset} are all known to be true.
For $r=2$, the sequence of polynomials $(s_{2,n}(x))_{n\geq 0}$
was studied by Elvey Price and Sokal \cite{Elvey-Price-Sokal_wardpoly},
who called them the \textbfit{Ward polynomials}:
compare the recurrence~\reff{eq.rec.stirlingsubset.2nd}
with \cite[eq.~(1.6)]{Elvey-Price-Sokal_wardpoly}.
They showed \cite[Theorem~1.1]{Elvey-Price-Sokal_wardpoly}
that these polynomials have a Thron-type continued fraction (T-fraction)
with coefficients $\alpha_n = nx$ and $\delta_n = n-1$.
It then follows from \cite[Theorem~9.9]{latpath_SRTR}
that this sequence of polynomials
is coefficientwise Hankel-totally positive in the variable $x$.
Therefore, part (d) of Conjecture~\ref{conj.subset}
is known to be true for $r=2$.
In this paper, we will prove part (c) for the case $r=2$.

\bigskip

Our main result on Hankel-total positivity is therefore the following:

\begin{theorem}[Hankel-total positivity of the second-order Stirling cycle polynomials]
   \label{thm.second.HankelTP}
The sequence of second-order Stirling cycle polynomials $(c_{2,n}(x))_{n\geq 0}$
is coefficientwise Hankel-totally positive in the variable $x$.
\end{theorem}

This theorem is, in fact, an immediate consequence of the
coefficientwise Hankel-total positivity
of the second-order Eulerian polynomials
\cite[Corollaries~12.35 and 12.36]{latpath_SRTR}
together with the identity \reff{eq.cycle.secondeuler.bis}.
So our main work will be to prove Proposition~\ref{prop.secondeuler}(a):
we will provide an algebraic proof (Section~\ref{sec.proof.prop.secondeuler})
and a combinatorial proof (Section~\ref{sec.hankeltp.2nd.cycle}).
For our combinatorial proof,
we need a new combinatorial interpretation for the numbers
$\stirlingcycle{n}{k}^{(2)}$,
using Stirling permutations, due to David Callan
\cite{Callan_private};
we will do this in Section~\ref{sec.interpretation.stirperm}.

More specifically, \cite[Corollary~12.35]{latpath_SRTR}
showed that the second-order Eulerian polynomial $E^{(2)}_n(x)$
is equal to a specialization of the multivariate Eulerian polynomial
$\scrp_n^{(2)}(x_0,x_1,x_2)$ introduced in \cite[Chapter~12]{latpath_SRTR}:
namely, $E^{(2)}_n(x) = \scrp_n^{(2)}(1,x,x)$.
(We will recall the precise definitions in Section~\ref{sec.r=2.hankel.proof}.)
It then follows immediately from \reff{eq.cycle.secondeuler.bis}
that the second-order Stirling cycle polynomial $c_{2,n}(x)$
is also equal to a specialisation of the multivariate Eulerian polynomial:
\be
   c_{2,n}(x) \;=\; \scrp_n^{(2)}(x,1+x,1+x)\;.
 \label{eq.ternary}
\ee
And from this, Theorem~\ref{thm.second.HankelTP} follows immediately
as a consequence of the coefficientwise Hankel-total positivity of the
multivariate Eulerian polynomials \cite[Corollary~12.3]{latpath_SRTR}.

An analogous argument using 
Proposition~\ref{prop.secondeuler}(b) and \reff{eq.subset.secondeuler.bis},
instead, provides an alternate proof
of the coefficientwise Hankel-total positivity
of the second-order Stirling subset polynomials $s_{2,n}(x)$,
using 2-branched S-fractions instead of classical T-fractions.
This was already observed in
\cite[Note Added, p.~12]{Elvey-Price-Sokal_wardpoly}.
The important remark here is that~\reff{eq.subset.secondeuler.bis}
gives us a way to write the second-order Stirling subset polynomials
$s_{2,n}(x)$ as a specialization of the multivariate Eulerian polynomials,
namely,
\be
   s_{2,n}(x) \;=\; \scrp_n^{(2)}(x,x,1+x)\;.
 \label{eq.ternary.subset}
\ee
See the Remark at the end of Section~\ref{sec.r=2.hankel.proof}.

Our second main contribution in this paper
is to introduce, in Section~\ref{sec.interpretation.other},
three further combinatorial interpretations
for the numbers $\stirlingcycle{n}{k}^{(2)}$;
thus, in this paper we provide a total of five combinatorial interpretations
for the second-order Stirling cycle numbers.
Even though these latter interpretations are not necessary
for our proof of Theorem~\ref{thm.second.HankelTP},
and we are unable at present to see how to use them to 
prove any of our other main conjectures,
we think that they are still interesting in their own right
and provide interesting consequences;
they also suggest new conjectures, as we will explain.
We begin, in Section~\ref{sec.interpretation.ternary},
by introducing our third combinatorial interpretation
for the second-order Stirling cycle numbers,
in terms of increasing ternary trees,
using a well-known bijection to Stirling permutations.
Then, in Section~\ref{sec.interpretation.other.inctree},
we introduce our fourth interpretation in terms of 
certain families of weighted increasing trees with ordered children,
and we provide two different proofs:
one bijective and one using generating functions.
Both proofs suggest new lines of investigation,
and also show that the hypothesis in~\cite[Theorem~1.1(c)]{latpath_lah}
due to P\'etr\'eolle and Sokal, though sufficient, is far from necessary.
In Section~\ref{sec.interpretation.other.phylogenetic}
we introduce our fifth and final interpretation
in terms of certain leaf-labeled trees.
We then show the relation to the multivariate Ward polynomials
introduced by Elvey Price and Sokal
\cite[pp.~9--11]{Elvey-Price-Sokal_wardpoly},
and we show that these multivariate polynomials are a common generalization 
of both the second-order Stirling cycle polynomials $c_{2,n}(x)$
and the second-order Stirling subset polynomials $s_{2,n}(x)$.
This observation leads us very naturally to another triangle of numbers 
that also seems to satisfy our four positivity properties; 
we state this in Conjecture~\ref{conj.phylogenetic.ordered}.

\bigskip

Our main results on real-rootedness are the following:

\begin{theorem}[Location of zeros for the of the second-order Stirling cycle and subset polynomials]
For each $n\geq 1$,
the polynomials $c_{2,n}(x)$ and $s_{2,n}(x)$ have $0$ as a root.
The zeros of the polynomial sequence $\big(c_{2,n}(x)/x \big)_{n\geq 1}$,
and also of the sequence $\big(s_{2,n}(x)/x \big)_{n\geq 1}$,
are simple, interlacing, and lie in the interval $(-1, 0)$.
   \label{thm.stirling.2nd.zeros}
\end{theorem}

Thus, in Theorem~\ref{thm.stirling.2nd.zeros} we prove not only
that the zeros are non-positive, but that in fact they lie
in a specific interval.

As mentioned previously,
it is known \cite[pp.~395, 399]{Karlin_68}
that the Toeplitz-total positivity of the rows of a lower-triangular matrix
is equivalent to the negative-real-rootedness
of its row-generating polynomials.
It therefore follows from Theorem~\ref{thm.stirling.2nd.zeros}
that the lower-triangular Toeplitz matrices
$\big( \stirlingcycle{n}{i-j}^{(2)} \big)_{i,j\geq 0}$
and
$\big( \stirlingsubset{n}{i-j}^{(2)} \big)_{i,j\geq 0}$
are totally positive for every $n \ge 0$.
In particular, the positivity of the contiguous $2\times 2$ minors 
of these Toeplitz matrices implies
that the row sequences of the second-order Stirling cycle and subset matrices
$C^{(2)}$ and $S^{(2)}$ are log-concave.

One can show that Theorem~\ref{thm.stirling.2nd.zeros}
follows very easily from a very general theorem
due to Liu and Wang \cite{Liu_06}.
However, for the convenience of the reader
we will provide a direct proof of this theorem
in Section~\ref{sec.zeros.2nd}.

Even though the zeros of the polynomials $c_{r,n}(x)$ and $s_{r,n}(x)$
do not lie on the real line for $r\geq 3$,
it might be worthwhile to try to
understand how they are distributed in the complex plane.
In Section~\ref{sec.zero.distribution}, 
we will provide some plots for the distribution of zeros of the polynomials
for orders $r\geq 3$.
These plots will suggest that these zeros, when correctly normalised,
seem to accumulate around some limiting curves.
However, at this moment we are unable to provide 
any precise conjectures or proofs describing these limiting curves.

\subsection{Higher-order Eulerian and quasi-Eulerian numbers:
   Statement of conjectures}
\label{sec.conj.eulerian}

A major longstanding open problem in the theory of total positivity is 
Brenti's conjecture on the total positivity of the triangle of Eulerian numbers 
\cite{Brenti_96}.
Brenti's conjecture is in fact the first in a family of conjectures
concerning the higher-order Eulerian numbers
that arose in discussions between us and Alex Dyachenko
\cite{Dyachenko_private_Eulerian}:

\begin{conjecture}[Conjecture for $r$-th order Eulerian numbers \cite{Dyachenko_private_Eulerian}]
   \label{conj.secondeuler}
For any $r\geq 1$, the following are true:
\begin{itemize}
\item[(a)] The $r$-th order Eulerian triangle $E^{(r)}=\big(\euler{n}{k}^{\! (r)}\big)_{n,k\geq 0}$
      is totally positive.
   \item[(b)] The shifted-reversed $r$-th order Eulerian triangle
      $\widecheck{E}^{(r)}=\big(\euler{n}{n-k-1}^{\! (r)}\big)_{n,k\geq 0}$ is totally positive.
\end{itemize}
\end{conjecture}

We remark that (a) and (b) are equivalent for $r=1$,
by virtue of the symmetry $\euler{n}{k} = \euler{n}{n-k-1}$ for $n \ge 1$;
but this is not the case for $r \ge 2$,
where the two conjectures need to be considered separately.
Conjecture~\ref{conj.secondeuler} is unproven even for $r=1$,
where it is Brenti's conjecture from 1996 \cite{Brenti_96}.
Dyachenko \cite{Dyachenko_private_Eulerian}
has verified Conjecture~\ref{conj.secondeuler}(a,b)
for the leading $512\times 512$ submatrix for $r=1$, 
and for the leading $256\times 256$ submatrix for $2 \le r \le 6$,
by computing the bidiagonal factorization
using Neville elimination \cite{Gasca_92}.

Conjecture~\ref{conj.secondeuler} asks our questions~(a) and~(b)
for the higher-order Eulerian triangles.
Questions~(c) and~(d) for these triangles
have already been answered in the affirmative:
the negative-real-rootedness of the row-generating polynomials 
is due to Brenti \cite[Theorem~6.6.3(i)]{Brenti_89},
and their coefficientwise Hankel-total positivity
was proven by P\'etr\'eolle, Sokal and Zhu
\cite[Corollaries~12.35 and 12.36]{latpath_SRTR}
using branched continued fractions.



\bigskip

However, as mentioned previously, as far as we know there is no relation
between the $r$th-order Stirling cycle or subset numbers
and the $r$th-order Eulerian triangles $E^{(r)}$ for $r\geq 3$.
Rather, Proposition~\ref{prop.secondeuler} suggests that
the matrices that are analogous in this context for $r\ge 3$
are $\widecheck{C}^{(r)} \, B^{-1}$ and $\widecheck{S}^{(r)} \, B^{-1}$.
We call these matrices the
\textbfit{$\bm{r^{th}}$-order quasi-Eulerian cycle triangle}
and the
\textbfit{$\bm{r^{th}}$-order quasi-Eulerian subset triangle},
respectively,
and we denote them as
\be
   Q_{\rm C}^{(r)} \:\eqdef\: \widecheck{C}^{(r)} \, B^{-1}
   ,\quad
   Q_{\rm S}^{(r)} \:\eqdef\: \widecheck{S}^{(r)} \, B^{-1}
   \;.
   \label{eq.def.quasieul.triangle}
\ee
We call their entries the
\textbfit{$\bm{r^{th}}$-order quasi-Eulerian cycle numbers}
and the
\textbfit{$\bm{r^{th}}$-order quasi-Eulerian subset numbers},
respectively:
\be
   \eulersecond{n}{k}^{\! (r)}_{\!\rm C}
   \:\eqdef\: (Q_{\rm C}^{(r)})_{n,k}
   ,\quad
   \eulersecond{n}{k}^{\! (r)}_{\!\rm S}
   \:\eqdef\: (Q_{\rm S}^{(r)})_{n,k}
   \;.
\ee
Thus, when $r=2$, we have from Proposition~\ref{prop.secondeuler}
\be
Q_{\rm C}^{(2)}
\;=\;
E^{(2)}\;,
\qquad
Q_{\rm S}^{(2)}
\;=\;
\widecheck{E}^{(2)}
  \;,
\label{eq.quasieul.secondeul}
\ee
so that the two matrices are reversals of each other.
But for $r \neq 2$ this is not the case,
and the four matrices
$Q_{\rm C}^{(r)}$, $Q_{\rm S}^{(r)}$,
$\bigl(Q_{\rm C}^{(r)}\bigr)^{\rm rev}$,
$\bigl(Q_{\rm S}^{(r)}\bigr)^{\rm rev}$
need to be considered separately.

In Appendices~\ref{app.eq.cycle} and \ref{app.eq.subset}, 
we provide tables of the triangles $Q_{\rm C}^{(r)}$
and $Q_{\rm S}^{(r)}$, respectively,
for $1\leq r\leq 5$ and the first few values of $n$.
When $r=1$, the matrices $Q_{\rm C}^{(1)}$ and $Q_{\rm S}^{(1)}$
are not even entrywise positive:
some of their entries are negative integers.
We therefore restrict our discussion henceforth to $r\geq 2$.

For all $r\geq 2$, we will provide a combinatorial interpretation for the
$r$th-order quasi-Eulerian {\em cycle}\/ numbers
in Section~\ref{sec.interpretation.stirperm.gen}
(Proposition~\ref{prop.combi.interpret.rth}).
We leave it as an open problem to construct a combinatorial interpretation
for the $r$th-order quasi-Eulerian {\em subset}\/ numbers.

For $r=2$, Conjecture~\ref{conj.secondeuler} asserts that
the foregoing four triangles are totally positive.
For $r \ge 3$, the triangles $Q_{\rm C}^{(r)}$ and $Q_{\rm S}^{(r)}$
are {\em not}\/ totally positive:
this is an immediate consequence of the definition
\reff{eq.def.quasieul.triangle}
together with the failure of total positivity
for the triangles $\widecheck{C}^{(r)}$ and $\widecheck{S}^{(r)}$
that was noted in Conjectures~\ref{conj.cycle}(b) and \ref{conj.subset}(b)
[see Appendix~\ref{app.sec.report}].
This fact makes the following conjecture somewhat more surprising:

\begin{samepage}
\begin{conjecture}[Conjectures for the reversed higher-order quasi-Eulerian triangles]
\hfill\break
\vspace*{-10mm}
\begin{itemize}
\item[(b)] The matrix $\bigl(Q_{\rm C}^{(r)}\bigr)^{\rm rev}$
   is totally positive for all $r \ge 2$.
\vspace*{-1mm}
\item[(b${}'$)] The matrix $\bigl(Q_{\rm S}^{(r)}\bigr)^{\rm rev}$
   is totally positive for all $r \ge 2$.
\end{itemize}
  \label{conj.quasi-Eulerian.TP}
\end{conjecture}
\end{samepage}

\noindent
We verified the total positivity of these matrices for $3 \le r \le 10$
for the leading $50 \times 50$ submatrix.

Next, for $r\geq 2$ we define the
\textbfit{$\bm{r^{th}}$-order quasi-Eulerian cycle polynomials}
and the 
\textbfit{$\bm{r^{th}}$-order quasi-Eulerian subset polynomials}
to be the generating polynomials
\begin{eqnarray}
Q^{(r)}_{{\rm C},n}(x)
\;\eqdef\;
\sum_{k=0}^{n} \eulersecond{n}{k}^{\! (r)}_{\!\rm C}  x^k\;,
\label{eq.def.quasieulpoly}
\\
Q^{(r)}_{{\rm S},n}(x)
\;\eqdef\;
\sum_{k=0}^{n} \eulersecond{n}{k}^{\! (r)}_{\!\rm S}  x^k\;,
\label{eq.def.quasieulpoly.subset}
\end{eqnarray}
respectively.
Then by definition, we have
\begin{eqnarray}
x^n c_{r,n}(1/x)  \;=\; Q^{(r)}_{{\rm C},n}(1+x)\;,
\qquad
x^n s_{r,n}(1/x)  \;=\; Q^{(r)}_{{\rm S},n}(1+x)\;,
\label{eq.quasieulpoly.hostirling}
\end{eqnarray}
which are~\reff{eq.cycle.secondeuler} 
and~\reff{eq.subset.secondeuler}, respectively, when $r=2$.

For $r \ge 3$, these polynomials are not real-rooted:
that follows from \reff{eq.quasieulpoly.hostirling}
together with the non-real-rootedness of the
polynomials $c_{r,n}(x)$ and $s_{r,n}(x)$
that was noted in Conjectures~\ref{conj.cycle}(c) and \ref{conj.subset}(c)
[see Proposition~\ref{prop.cycle.discriminant}
 and Corollaries~\ref{cor.cycle.nonreal} and \ref{cor.subset.nonreal}
 and the subsequent remarks].

Likewise, it follows from \reff{eq.quasieulpoly.hostirling}
and the failure of coefficientwise Hankel-total positivity for $s_{r,n}(x)$
that was noted in Conjecture~\ref{conj.subset}(d)
that the polynomials $Q^{(r)}_{{\rm S},n}(x)$ are not
coefficientwise Hankel-totally positive.
On the other hand, we do conjecture coefficientwise Hankel-total positivity
for the quasi-Eulerian {\em cycle}\/ polynomials:

\begin{conjecture}[Conjectured Hankel-total positivity of the quasi-Eulerian cycle polynomials]
For all $r \ge 3$, the sequence $(Q^{(r)}_{{\rm C},n}(x))_{n\geq 0}$
of the $r$th-order quasi-Eulerian cycle polynomials
is coefficientwise Hankel-totally positive in the variable $x$.
   \label{conj.quasi-Eulerian.hankel}
\end{conjecture}

\noindent 
From \reff{eq.quasieulpoly.hostirling}
we see that Conjecture~\ref{conj.quasi-Eulerian.hankel}
strengthens Conjecture~\ref{conj.cycle}(d).
We have tested Conjecture~\ref{conj.quasi-Eulerian.hankel}
up to the $9 \times 9$ Hankel matrix for $3 \le r \le 10$.

\subsection{Structure of this paper}
\label{sec.structure}

The plan of this paper is as follows:
In Section~\ref{sec.proof.prop.secondeuler}
we give an algebraic proof of Proposition~\ref{prop.secondeuler}.
In Section~\ref{sec.hankeltp.2nd.cycle},
we prove the coefficientwise Hankel-total positivity
of the second-order Stirling cycle polynomials
(Theorem~\ref{thm.second.HankelTP}),
based on two combinatorial interprerations of the
second-order Stirling cycle numbers.
Then, in Section~\ref{sec.interpretation.other},
we give our third, fourth and fifth combinatorial interpretations for the
second-order Stirling cycle numbers;
these lead to some interesting observations and another conjecture.
Next, in Section~\ref{sec.interpretation.stirperm.gen},
we generalize our second interpretation for 
the second-order Stirling cycle numbers to higher orders;
this gives a second combinatorial interpretation for 
the $r$th-order Stirling cycle numbers,
and also for the $r$th-order quasi-Eulerian cycle numbers,
for all $r\geq 2$.
After this, in Section~\ref{sec.zeros.2nd},
we prove Theorem~\ref{thm.stirling.2nd.zeros}
on the negative-real-rootedness
of the second-order Stirling cycle and subset polynomials.
Finally, in Section~\ref{sec.zero.distribution}
we provide some plots for the zeros of
\hbox{$r$-th} order Stirling cycle and subset polynomials for 
orders $r\in [3,10]$.

In Appendix~\ref{app.sec.tables}, 
we provide tables of the $r$th-order Stirling cycle and subset numbers
for $1\leq r\leq 4$ and the first few values of $n$;
and in Appendix~\ref{app.sec.tables.qe},
we provide tables of the $r$th-order quasi-Eulerian cycle and subset numbers
for $2\leq r\leq 4$ and the first few values of $n$.
Then in Appendix~\ref{app.sec.report},
we provide a report on the systematic verification of our conjectures 
and mention up to which order we have checked our conjectures.

\section{Algebraic proof of Proposition~\ref{prop.secondeuler}}
\label{sec.proof.prop.secondeuler}

In this section we will prove Proposition~\ref{prop.secondeuler}
by using a result from \cite{GKP_recurrence}
concerning lower-triangular arrays that satisfy a
Graham--Knuth--Patashnik (GKP) recurrence.

\begin{lemma} {\bf \protect\cite[Corollary~A.15]{GKP_recurrence}}
   \label{cor6.prop.spivey.corollary5}
Let $\bA = \bigl( A(n,k) \bigr)_{n,k \ge 0}$
be a lower-triangular array defined by the GKP recurrence
\be
  A(n,k)
  \;=\;
  (\alpha n + \beta k + \gamma)    \, A(n-1,k)
  \:+\:
  (\alpha' n + \beta' k + \gamma') \, A(n-1,k-1)
\ee
for $n \ge 1$, with initial conditions $A(0,k) = \delta_{k0}$ and $A(n,-1) = 0$.
Let $B_\xi$ be the weighted binomial matrix with entries
\be
   (B_\xi)(n,k)  \;=\;  \binom{n}{k} \, \xi^{n-k}
   \;.
\ee
Then the matrix product $\bC = \bA \, B_\xi$ satisfies the recurrence
\begin{eqnarray}
   C(n,k)
   & = &
   \bigl[ (\alpha + \xi \alpha') n \,+\,
          (\beta + 2\xi \beta') k  \,+\,
          \gamma + \xi(\beta' +\gamma')
   \bigr] \, C(n-1,k)
           \nonumber \\[2mm]
   & &
   \quad +\;
   (\alpha' n + \beta' k + \gamma') \, C(n-1,k-1)
           \nonumber \\[2mm]
   & &
   \quad +\;
   \xi \, (\beta + \xi\beta') \, (k+1) \, C(n-1,k+1)
 \label{eq.cor6.prop.spivey.corollary5}
\end{eqnarray}
for $n \ge 1$, with initial conditions $C(0,k) = \delta_{k0}$ and $C(n,-1) = 0$.
\end{lemma}

In particular, if $\xi = -\beta/\beta'$ we again obtain a recurrence
of GKP form.

\proofof{Proposition~\ref{prop.secondeuler}}
(a) The recurrence \reff{def.eulersecond} for the second-order Eulerian numbers
is of GKP form with 
$(\alpha,\beta,\gamma,\alpha',\beta',\gamma') = (0,1,1,2,-1,-1)$.
Applying Lemma~\ref{cor6.prop.spivey.corollary5} with $\xi = 1$,
we deduce that $\bC = E^{(2)} \, B_1$ satisfies the GKP recurrence
\be
   C(n,k)  \;=\;  (2n-k-1) \, C(n-1,k)  \:+\: (2n-k-1) \, C(n-1,k-1)
   \;.
\ee
This is exactly what is obtained from the recurrence
\reff{eq.rec.stirlingcycle.2nd} for the second-order Stirling cycle numbers
by replacing $k \to n-k$.

(b) Making the substitution $k \to n-k-1$ in \reff{def.eulersecond},
we see that the shifted-reversed second-order Eulerian triangle
$\widecheck{E}^{(2)}$ satisfies a GKP recurrence with
$(\alpha,\beta,\gamma,\alpha',\beta',\gamma') = (1,1,0,1,-1,0)$.
Applying Lemma~\ref{cor6.prop.spivey.corollary5} with $\xi = 1$,
we deduce that $\bC = \widecheck{E}^{(2)} \, B_1$ satisfies the GKP recurrence
\be
   C(n,k)  \;=\;  (2n-k-1) \, C(n-1,k)  \:+\: (n-k) \, C(n-1,k-1)
   \;.
\ee
This is exactly what is obtained from the recurrence
\reff{eq.rec.stirlingsubset.2nd} for the second-order Stirling subset numbers
by replacing $k \to n-k$.
\qed

This proof is straightforward
(once one has in hand Lemma~\ref{cor6.prop.spivey.corollary5}),
but it is not terribly illuminating.
Much of the rest of this paper is therefore devoted to providing
combinatorial interpretations of the second-order Stirling cycle numbers
that make the identity of Proposition~\ref{prop.secondeuler}(a) transparent.

\section{Combinatorial proof of Proposition~\ref{prop.secondeuler}(a)}
\label{sec.hankeltp.2nd.cycle}

In this section we will provide 
a combinatorial proof of Proposition~\ref{prop.secondeuler}(a).
We begin by recalling the combinatorial definition
of the numbers $\stirlingcycle{n}{k}^{(2)}$
(Section~\ref{sec.interpretation.derange}).
We then proceed in Section~\ref{sec.interpretation.stirperm}
to establish our second combinatorial interpretation of
$\stirlingcycle{n}{k}^{(2)}$:
it counts Stirling permutations of the multiset $\{1,1,2,2,\ldots,n,n\}$
with optional dots at ascents and containing $n-k$ dots.
Finally, in Section~\ref{sec.r=2.hankel.proof},
this interpretation will immediately imply
Proposition~\ref{prop.secondeuler}(a) as a natural consequence;
we obtain Theorem~\ref{thm.second.HankelTP} as a corollary.


\subsection{Interpretation I: Derangements}
\label{sec.interpretation.derange}

Let us recall that
$\stirlingcycle{n}{k}^{(2)} = \stirlingcycle{n+k}{k}_2$
is defined to be the number of permutations of $[n+k]$
containing $k$ cycles, all of which are of size $\ge 2$.
In other words, it is the number of derangements of $[n+k]$
containing $k$ cycles.
We write $\scrd_{N,k}$ for the set of derangements of $[N]$
with $k$ cycles.
Therefore, $\stirlingcycle{n}{k}^{(2)} = |\scrd_{n+k,k}|$.
These numbers satisfy the recurrence~\reff{eq.rec.stirlingcycle.2nd}.

\subsection{Interpretation II: Stirling permutations}
\label{sec.interpretation.stirperm}

We will now provide a second interpretation for the numbers
$\stirlingcycle{n}{k}^{(2)}$:
they count Stirling permutations with optional dots in ascents.
This interpretation
can be guessed from the entries in \cite[A032188]{OEIS},
and the proof that follows is due to David Callan \cite{Callan_private}.

Let us recall that a \textbfit{Stirling permutation of size $\bm{n}$}
is a permutation of the multiset $\{1,1,2,2,\ldots, n,n\}$
with the {\em Stirling property}\/:
for every $i \in [n]$,
all entries between the two occurrences of $i$ exceed $i$.
Equivalently, a Stirling permutation of size~$n$
is permutation of the multiset $\{1,1,2,2,\ldots, n,n\}$
that avoids the pattern $212$.
We recall~\cite{Gessel_78} that there are $(2n-1)!!$
Stirling permutations of size $n$.
An \textbfit{ascent-marked Stirling permutation of size $\bm{n}$}
is now defined to be a Stirling permutation in which
a dot may or may not be inserted in each ascent
(we recall that an {\em ascent}\/ is a pair $a_i, a_{i+1}$
 such that $a_i < a_{i+1}$;
 the dot, if any, is inserted between $a_i$ and $a_{i+1}$).
Let $\scrv_n$ denote the set of 
ascent-marked Stirling permutations of size $n$;
and let $\scrv_{n,k}$ denote the set of 
ascent-marked Stirling permutations of size $n$ with $k$ dots.
Thus, 
\begin{subeqnarray}
	\scrv_{2,0} & = & \{1122, 1221, 2211\}\\
	\scrv_{2,1} & = & \{11\bolddot 22, 1\bolddot 221\}.
\end{subeqnarray}
Let $v_{n,k} = \left| \scrv_{n,k} \right|$.
We will prove:

\begin{lemma}
   \label{lem.perm.cyc.asc.marked}
The numbers $v_{n,k}$ satisfy the recurrence
\be
   v_{n,k} \;=\; (2n - k - 1) \, \big(v_{n-1,k} + v_{n-1,k-1} \big)
     \label{eq.rec.vnk}
\ee
with initial conditions $v_{0,k} = \delta_{k0}$ and $v_{n,-1}=0$.
Therefore $\stirlingcycle{n}{k}^{(2)} = v_{n,n-k}$.
\end{lemma}

\par\medskip\noindent{\sc Proof} \cite{Callan_private}.
Let $\sigma \in \scrv_{n,k}$ be an ascent-marked Stirling permutation 
of size $n$ with $k$~dots.
Deleting the two $n$'s (which are necessarily adjacent)
together with the dot (if any) preceding them
gives an element of $\scrv_{n-1,k}$ if there was no dot before $nn$,
and an element of $\scrv_{n-1,k-1}$ if there was a dot before $nn$.

Conversely, to get an element $\sigma \in \scrv_{n,k}$
with $nn$ not preceded by a dot from a one-size-smaller
ascent-marked Stirling permutation, $nn$ can be inserted
into an element of $\scrv_{n-1,k}$
in any of the $2n - 1$ spaces except those occupied by a dot;
thus, there are $2n-1-k$ choices. 
To get a $\sigma \in \scrv_{n,k}$ with $nn$ preceded by a dot, 
$\bolddot nn$ can be inserted
into an element of $\scrv_{n-1,k-1}$
in any of the $2n - 1$ spaces except the first space
and those occupied by a dot;
again the number of choices is $2n-1-k$.
This proves the recurrence~\reff{eq.rec.vnk}.

{}From \reff{eq.rec.vnk} it follows that the numbers
$u_{n,k} = v_{n,n-k}$ satisfy the recurrence
\be
   u_{n,k} \;=\; (n+k-1) (u_{n-1,k-1} + u_{n-1,k})
   \;,
\ee
which is the same as the recurrence~\reff{eq.rec.stirlingcycle.2nd}
satisfied by the numbers $\stirlingcycle{n}{k}^{(2)}$.
Since these numbers also satisfy the same initial conditions,
we conclude that $\stirlingcycle{n}{k}^{(2)} = u_{n,k}$.
\qed

\subsection{Proof of Proposition~\ref{prop.secondeuler}(a)
and Theorem~\ref{thm.second.HankelTP}}
\label{sec.r=2.hankel.proof}

We are now ready to prove Proposition~\ref{prop.secondeuler}(a).

\proofof{Proposition~\ref{prop.secondeuler}(a)}
To prove Proposition~\ref{prop.secondeuler}(a),
we need to show for every $n,k\geq 0$ that 
\be
\stirlingcycle{n}{n-k}^{(2)} 
\;=\;
\sum_{i=k}^n \eulersecond{n}{i}\binom{i}{k}
\label{eq.secondeuler.binom.cycle}
\ee
Using Lemma~\ref{lem.perm.cyc.asc.marked}, we know that the left-hand side
is the number of ascent-marked Stirling permutations of size $n$ with $k$ dots.
This is clearly equal to the right-hand side because
we can pick such an ascent-marked Stirling permutation by first 
picking a Stirling permutation of size $n$ with $i$ ascents 
(counted by $\eulersecond{n}{i}$)
and then insert dots in $k$ of the $i$ ascents
(counted by $\binom{i}{k}$).
\qed

Theorem~\ref{thm.second.HankelTP} will now be an immediate consequence,
but we first recall some definitions and results from
\cite[Chapter~12]{latpath_SRTR}.
For a fixed integer $m\geq 1$
and indeterminates $\bfx = (x_0,x_1,\ldots,x_m)$,
we let $\scrq_n^{(m)}(\bfx)$ be the generating polynomial for
increasing $(m+1)$-ary trees on the vertex set $[n] \cup \{0\}$
in which each $i$-edge gets a weight $x_i$.
They are homogeneous polynomials of degree $n$,
and symmetric under permutations of $x_0,\ldots,x_m$.
We then define $\scrp_n^{(m)}(\bfx)$ by
\begin{subeqnarray}
        \scrp_0^{(m)}(\bfx) &=& 1\\
        \scrp_n^{(m)}(\bfx) &=&  x_0 \scrq_{n-1}^{(m)}(\bfx)
        \qquad \hbox{for $n\geq 1$}
\end{subeqnarray}
Therefore, $\scrp_n^{(m)}(\bfx)$ is the generating polynomial for increasing
$(m+1)$-ary trees on the vertex set $[n] \cup \{0\}$
such that the only edge (if any) emanating from the root $0$ is a $0$-edge
and in which each $i$-edge gets a weight $x_i$.
These polynomials were called the multivariate Eulerian polynomials
in \cite{latpath_SRTR}.
They are homogeneous polynomials of degree $n$,
and symmetric under permutations of $x_1,\ldots,x_m$;
but $x_0$ plays a special role.

In \cite[Theorem~9.8]{latpath_SRTR},
P\'etr\'eolle, Sokal and Zhu
showed that the sequence of $m$-Stieltjes--Rogers polynomials
$S_n^{(m)}(\balpha)$
(see \cite[Definition~2.2]{latpath_SRTR} for the definition)
is coefficientwise Hankel-totally positive in the indeterminates
$\balpha = (\alpha_i)_{i \ge m}$.
Then in \cite[eq.~(12.5) and Theorem~12.1(a)]{latpath_SRTR}
they proved that the multivariate Eulerian polynomials $\scrp^{(m)}_n(\bfx)$
can be obtained as a specialization of the $m$-Stieltjes--Rogers polynomials:
\be
   \alpha_{m+j+(m+1)k}  \;=\;  (k+1) x_j
   \;.
\ee
Since these coefficients $\alpha_i$ are polynomials (in fact monomials)
with nonnegative coefficients in the indeterminates $\bfx$,
it follows that the sequence $\big( \scrp_n^{(m)}(\bfx) \big)_{n\geq 0}$
is coefficientwise Hankel-totally positive in the indeterminates $\bfx$.
We can summarise these results as follows:

\begin{theorem}[\!\!\!{\cite[eq.~(12.5) and Theorems~12.1(a) and~9.8]{latpath_SRTR}}]
The polynomials $\scrp_n^{(m)}(\bfx)$ are equal to the
$m$-Stieltjes--Rogers polynomials $S_n^{(m)}(\balpha)$ where the weights
$\balpha = (\alpha_i)_{i \ge m}$ are given by
\be
        \alpha_{m+j+(m+1)k}  \;=\;  (k+1) x_j.
\ee
Thus, the sequence of polynomials $\big( \scrp_n^{(m)}(\bfx) \big)_{n\geq 0}$
is coefficientwise Hankel-totally positive with respect to the variables
$x_0,x_1,\ldots, x_m$.
\label{thm.multivariateEulerian}
\end{theorem}

\proofof{Theorem~\ref{thm.second.HankelTP}}
The second-order Eulerian polynomial $E^{(2)}_n(x)$
is equal to a specialization of the multivariate Eulerian polynomial: namely,
$E^{(2)}_n(x) = \scrp_n^{(2)}(1,x,x)$.
As the polynomials $\scrp_n^{(2)}$ are homogeneous in all three variables,
equation~\reff{eq.cycle.secondeuler.bis}, which is a consequence of 
Proposition~\ref{prop.secondeuler}(a), gives us equation~\reff{eq.ternary}.
The theorem then follows from Theorem~\ref{thm.multivariateEulerian}.
\qed

\medskip

{\bf Remark.}
We can now prove the analogous equation~\reff{eq.ternary.subset}
for the second-order Stirling subset polynomials.
Using \reff{eq.subset.secondeuler.bis} together with
$E^{(2)}_n(x) = \scrp_n^{(2)}(1,x,x)$
and the homogeneity and symmetry of the multivariate Eulerian polynomials,
we have for $n \ge 1$
\begin{subeqnarray}
   s_{2,n}(x)  & = &  x \, (1+x)^{n-1} \, E^{(2)}_n\Big( {x \over 1+x} \Big)
       \\[2mm]
   & = &
   {x \over 1+x} \: \scrp_n^{(2)}(1+x,x,x)
       \\[2mm]
   & = &
   x \: \scrq_{n-1}^{(2)}(1+x,x,x)
       \\[2mm]
   & = &
   x \: \scrq_{n-1}^{(2)}(x,x,1+x)
       \\[2mm]
   & = &
   \scrp_n^{(2)}(x,x,1+x)
     \;.
\end{subeqnarray}
And of course for $n=0$ we have $s_{2,0}(x) = \scrp_0^{(2)}(x,x,1+x) = 1$.
\myendremark



\section{Second-order Stirling cycle numbers: \hfill\break
         Three further combinatorial interpretations}
\label{sec.interpretation.other}

In this section we will provide three more combinatorial interpretations
for the second-order Stirling cycle numbers.
In Section~\ref{sec.interpretation.ternary}
we provide an interpretation using increasing ternary trees.
In Section~\ref{sec.interpretation.other.inctree}
we introduce an interpretation involving increasing trees 
with ordered children.
In Section~\ref{sec.interpretation.other.phylogenetic}
we introduce an interpretation involving a class of leaf-labeled trees
with cyclically-ordered children.

\subsection{Interpretation III: Increasing ternary trees}
\label{sec.interpretation.ternary}

We begin by recalling \cite[p.~573]{Stanley_12}
the recursive definition of an \textbfit{$\bm{m}$-ary tree}
($1 \le m < \infty$):
it is either empty or else consists of a root together with
an ordered list of $m$ subtrees $(T_0,\ldots,T_{m-1})$,
each of which is an $m$-ary tree (which may be empty).
We draw an edge from each vertex to the root of each of its nonempty subtrees;
an edge from a vertex to the root of its $i$th subtree ($0 \le i \le m-1$)
will be called an \textbfit{$\bm{i}$-edge}.
When $m=2$, we refer to a \textbfit{binary tree}
and to its \textbfit{left} and \textbfit{right} edges.
When $m=3$, we refer to a \textbfit{ternary tree}
and to its \textbfit{left}, \textbfit{middle} and \textbfit{right} edges.

A rooted tree on a totally ordered vertex set
is called \textbfit{increasing}
if each vertex is higher-numbered than its parent;
equivalently, the labels increase along each path downward from the root.
It is well known \cite[p.~45, Proposition~1.5.3(a)]{Stanley_12}
that the number of increasing binary trees on the vertex set $[n]$ is $n!$,
and more generally that the number of increasing $m$-ary trees
on the vertex set $[n]$ is the multifactorial $F_n^{(m-1)}$
\cite[p.~30, Example~1]{Bergeron_92},
where
$F_n^{(r)}  \;\eqdef\;  \prod_{j=0}^{n-1} (1+jr)$.
In particular, the number of increasing ternary trees on the vertex set $[n]$
is $(2n-1)!!$.

 We write $\scrt_n$ for the set of increasing ternary trees
 on the vertex set $[n]$.
 (In particular, $\scrt_0$ has a single element, namely the empty tree.)
 For $n \ge 1$, we write $\scrt_{n,k}$ for the set of increasing ternary trees
 on the vertex set $[n]$ with $k$ left edges.

 In what follows we shall also use a variant:
 We write $\scrtprime_n$ for the set of increasing ternary trees
 on the vertex set $[n] \cup \{0\}$
 for which the only edge (if any) emanating from the root $0$ is a left edge.
 (Thus, $\scrtprime_0$ has a single element, namely the tree with only
  the root vertex 0; and in all trees $T \in \scrtprime_n$ with $n \ge 1$,
  the root vertex 0 has a left child 1 and no other children.)
 We write $\scrtprime_{n,k}$ for the set of increasing ternary trees
 $T \in \scrtprime_n$ that have $k$ left edges.

For $n\geq 1$, there is an obvious bijection between $\scrt_n$ and $\scrtprime_n$,
 which also maps $\scrt_{n,k}$ to $\scrtprime_{n,k+1}$: namely,
 given $T \in \scrt_n$ we define $T' \in \scrtprime_n$
 by introducing a new root vertex 0 and making $T$ its left subtree;
 and given $T' \in \scrtprime_n$, we define $T \in \scrt_n$ by deleting the 
 root vertex 0 and making the vertex 1 
 the new root.
 For $n=0$, we augment this bijection by sending the empty tree $T\in \scrt_{0,0}$
 to the tree $T'\in \scrt_{0,0}$ that has only the root vertex $0$.


Since Stirling permutations 
of size $n$
and the set of increasing ternary trees $\scrt_n$
both have cardinality $(2n-1)!!$,
there clearly must exist a bijection between these two sets.
Indeed, such a bijection is well known \cite{Park_94a,Janson_11}.\footnote{
   This bijection was first mentioned by Park in \cite{Park_94a},
   who attributed it to Gessel.
   The details of this bijection were first explicitly stated
   by Janson, Kuba and Panholzer in \cite{Janson_11}.
}
We recall this construction, and the properties that we will need,
in the following lemma:

\begin{lemma}
   \label{lemma.stirling.ternary}
There is a bijection $\sigma \mapsto T(\sigma)$
where $\sigma$ is a Stirling permutations of size $n$
and $T(\sigma)\in \scrt_n$ is an increasing ternary tree on the vertex set $[n]$.
This bijection has the property that, for each letter $a \in [n]$,
the first (resp.~second) occurrence of the letter $a$ in $\sigma$
is the beginning of an ascent (i.e.~$a = a_i < a_{i+1}$)
if and only if the vertex $a$ in $T(\sigma)$
has a middle (resp.~right) child.
\end{lemma}

\proof
We begin by recalling the bijection $\sigma \mapsto T(\sigma)$,
following \cite{Janson_11}.
Let $\sigma = a_1\cdots a_{2n}$ be any word on the alphabet
$\mathbb{Z}_{\geq 1}$ such that
\begin{itemize}
   \item[(a)] each letter, if it occurs in $\sigma$, occurs exactly twice;
      and
   \item[(b)] $\sigma$ has the Stirling property, 
i.e.\ all entries in $\sigma$ between the two occurrences of~$a$ exceed $a$.
\end{itemize}
We then recursively define a ternary tree $T(\sigma)$ as follows: 
\begin{itemize}
	\item $T(\sigma)$ is the empty tree if and only if $\sigma$ is the empty word.
	\item Otherwise, let $a$ be the smallest letter of 
		$\sigma$ and let 
		$\sigma = \sigma_0 \, a \, \sigma_1 \, a \,\sigma_2$.
		Then $T(\sigma)$ is the ternary tree rooted at $a$
		with $T(\sigma_0)$, $T(\sigma_1)$ and $T(\sigma_2)$ 
		as the left, middle and right subtrees of $a$, respectively.
\end{itemize}
It is easy to see that $\sigma\mapsto T(\sigma)$ is a bijection
from Stirling permutations of size $n$
to the set $\scrt_n$ of increasing ternary trees on $[n]$.

From the description of this bijection,
it is clear that the first occurrence of the letter $a$ in $\sigma$
is the first element of an ascent pair
(i.e.~$a = a_i < a_{i+1}$)
if and only if $\sigma_1$ is not the empty word,
which is true if and only if 
the vertex labeled $a$ has a middle child in the ternary tree $T(\sigma)$.
Similarly, the second $a$ in $\sigma$ is the first element of an ascent pair
if and only if $\sigma_2$ is not the empty word,
which is true if and only if the vertex labeled
$a$ has a right child in the ternary tree $T(\sigma)$.
\qed

Thus, the ascents in $\sigma$ are in correspondence 
with the middle and right edges of $T(\sigma)$.

Now define an \textbfit{edge-marked increasing ternary tree}
to be an increasing ternary tree in which each middle or right edge
may be marked or unmarked.  (Left edges are always unmarked.)
Lemma~\ref{lemma.stirling.ternary} immediately implies:

\begin{corollary}
   \label{cor.stirling.ternary}
There is a bijection $\tau \mapsto T(\tau)$
where $\tau\in \scrv_n$ is an ascent-marked Stirling permutation
of size $n$
and $T(\tau)$ is an edge-marked increasing ternary tree on the vertex set $[n]$,
such that, for each letter $a \in [n]$,
the first (resp.~second) occurrence of the letter $a$ in $\tau$
is followed by a dot
if and only if the vertex $a$ in $T(\tau)$
has a middle (resp.~right) child connected by a marked edge.
\end{corollary}

Putting together Lemma~\ref{lem.perm.cyc.asc.marked}
and Corollary~\ref{cor.stirling.ternary}, we conclude:

\begin{proposition}
   \label{prop.secondorder.cycle.ternary}
The number of edge-marked increasing ternary trees on the vertex set $[n]$
with $n-k$ marked edges (hence $k-1$ unmarked edges)
is the second-order Stirling cycle number $\stirlingcycle{n}{k}^{(2)}$.
\end{proposition}

Now we use the bijection $\scrt_n \mapsto \scrtprime_n$
and to write out a different version of Proposition~\ref{prop.secondorder.cycle.ternary}.
First, we let $\scrthat_n$ denote the trees $T\in\scrtprime_n$ 
along with optional edge markings in which each middle or right edge
may be marked or unmarked.  (Left edges are always unmarked.)
Let $\scrthat_{n,k}\subseteq \scrthat_n$ denote the subset of these trees on the vertex set 
$[n] \cup \{0\}$ which have $k$ unmarked edges.
We then have the following rewriting of Proposition~\ref{prop.secondorder.cycle.ternary}:

\begin{proposition}
   \label{prop.secondorder.cycle.ternary.bis}
The number of edge-marked increasing ternary trees on the vertex set
$[n] \cup \{0\}$ with $k$ marked edges --- that is, $|\scrthat_{n,k}|$ ---
is the second-order Stirling cycle number $\stirlingcycle{n}{k}^{(2)}$.
\end{proposition}

\subsection{Interpretation IV: Increasing ordered trees}
\label{sec.interpretation.other.inctree}


An \textbfit{ordered tree} is a rooted tree
in which the children of each vertex are linearly ordered.
An \textbfit{increasing ordered tree} is an ordered tree
on a totally ordered vertex set (usually some set of integers)
such that each vertex is higher-numbered than its parent;
equivalently, the labels increase along each path downward from the root.
It~is well known \cite[Corollary~1(iv)]{Bergeron_92} \cite{Janson_08}
that the number of increasing ordered trees on $n+1$ vertices is $(2n-1)!!$.

A vertex in a rooted tree is called a \textbfit{leaf}
if it has no children;
it is called an \textbfit{internal vertex} (or a \textbfit{non-leaf})
if it has at least one child.

Recall that $\scrtprime_n$ 
denotes the set of increasing ternary trees on the vertex set $[n] \cup \{0\}$
such that the only edge (if any) emanating from the root $0$ is a left edge;
and $\scrtprime_{n,k}$ denotes the set of trees in $\scrtprime_n$
with $k$ left edges.
Now let $\scri_n$ denote the set of increasing ordered trees
on the vertex set $[n] \cup \{0\}$,
and let $\scri_{n,k}$ denote the set of trees in $\scri_n$
with $k$ internal vertices.\footnote{
   Note that when $n=0$, the unique element of $\scri_0$,
   namely the tree with only the root vertex 0,
   belongs to $\scri_{0,0}$, because in this case the root is a leaf.
}
We have already seen that $|\scrtprime_n| = |\scri_n| = (2n-1)!!$.
We will refine this equality as follows:

\begin{proposition}
   \label{thm.bij.planar.ternary}
The numbers $|\scri_{n,k}|$ are given by
\be
	|\scri_{n,k}| \;=\; |\scrtprime_{n,k}| \;=\; \eulersecond{n}{n-k}\;.
\ee
\end{proposition}

The second equality is a direct consequence of Lemma~\ref{lemma.stirling.ternary}.
For the first equality, we will give two proofs:
a bijective proof in Section~\ref{sec.interpretation.other.inctree.bij},
and a proof using generating functions
in Section~\ref{sec.interpretation.other.inctree.gf}.
Having established Proposition~\ref{thm.bij.planar.ternary},
we will finally be able to introduce {\em vertex-marked increasing ordered trees} in 
Section~\ref{subsec.vertex.marked},
which will be our fourth combinatorial interpretation for the numbers $\stirlingcycle{n}{k}^{(2)}$.

\subsubsection{Bijective proof of Proposition~\ref{thm.bij.planar.ternary}}
\label{sec.interpretation.other.inctree.bij}

For $n=0$ the claim of Proposition~\ref{thm.bij.planar.ternary} is trivial,
since $\scrtprime_0 = \scrtprime_{0,0} = \scri_0 = \scri_{0,0}$
has a unique element, namely the tree consisting only of the root vertex 0.
So from now on we fix $n \ge 1$.

We will prove Proposition~\ref{thm.bij.planar.ternary}
by constructing a bijection $\Phi_n \colon \scrtprime_n \to \scri_n$
that maps $\scrtprime_{n,k}$ onto $\scri_{n,k}$.
We will first construct the map $\Phi_n \colon \scrtprime_n \to \scri_n$;
we will then construct another map $\Psi_n \colon \scri_n \to \scrtprime_n$;
finally, we will show that the maps $\Phi_n$ and $\Psi_n$ are inverses
of each other.
Before doing all this, we first introduce some terminology concerning
trees $T\in \scrtprime_n$ and $S\in \scri_n$.
 
Let $T\in \scrtprime_n$ and let $j\in [n] \cup \{0\}$ be a vertex of $T$.
Consider the subtree of $T$ rooted at the vertex $j$ that uses
all the middle and right edges but no left edge.
We call this tree the \textbfit{middle-right subtree of $\bm{j}$ in $\bm{T}$}.

\begin{figure}[t]
\centering
\subfloat[\centering An increasing ternary tree $T\in \scrtprime_{10}$]{{
\newcommand{\nodea}{\node[draw,circle] (a) {$0$}
;}\newcommand{\nodeb}{\node[draw,circle] (b) {$1$}
;}\newcommand{\nodec}{\node[draw,circle] (c) {$2$}
;}\newcommand{\noded}{\node[draw,circle] (d) {$3$}
;}\newcommand{\nodee}{\node[draw,circle] (e) {$4$}
;}\newcommand{\nodef}{\node[draw,circle] (f) {$5$}
;}\newcommand{\nodeg}{\node[draw,circle] (g) {$6$}
;}\newcommand{\nodeh}{\node[draw,circle] (h) {$7$}
;}\newcommand{\nodei}{\node[draw,circle] (i) {$8$}
;}\newcommand{\nodej}{\node[draw,circle] (j) {$9$}
;}\newcommand{\nodeba}{\node[draw,circle] (ba) {$10$}
;}
\begin{tikzpicture}[auto]
\matrix[column sep=.1cm, row sep=0.5cm,ampersand replacement=\&]{
        \& \&         \&         \&         \&         \&  \nodea \\
        \& \&         \&         \&  \nodeb \&         \&         \\
        \& \& \noded  \&         \&  \nodec \&         \&  \nodef \\
        \nodeg  \& {\node {};}{\node {};}{\node {};} \&         \&  \nodee \&         \&  \nodeba \& \\
\& \nodeh \&         \&  \nodei \&  \nodej \&          \& \\
};
	\draw[ultra thick]        (a) -- (b) node [midway, left]  {$x$};
	\draw[ultra thick]        (b) -- (c) node [midway,left] {$y$};
	\draw[ultra thick]        (b) -- (d) node [midway,left] {$x$};
	\draw[ultra thick]       (b) -- (f) node [midway,right] {$y$};
	\draw[ultra thick]        (c) -- (e) node [midway,left] {$x$};
	\draw[ultra thick]       (c) -- (ba) node [midway,right] {$y$};
	\draw[ultra thick]       (d) -- (g) node [midway,left] {$x$};
	\draw[ultra thick]       (e) -- (i) node [midway,left] {$y$};
	\draw[ultra thick]       (e) -- (j) node [midway,right] {$y$};
	\draw[ultra thick]       (g) -- (h) node [midway,left] {$y$};
\end{tikzpicture}
	    }}%
	\subfloat[\centering An increasing ordered tree $S\in \scri_{10}$ ]{{
 \newcommand{\nodea}{\node[draw,circle] (a) {$0$}
;}\newcommand{\nodeb}{\node[draw,circle] (b) {$2$}
;}\newcommand{\nodec}{\node[draw,circle] (c) {$8$}
;}\newcommand{\noded}{\node[draw,circle] (d) {$4$}
;}\newcommand{\nodee}{\node[draw,circle] (e) {$9$}
;}\newcommand{\nodef}{\node[draw,circle] (f) {$10$}
;}\newcommand{\nodeg}{\node[draw,circle] (g) {$1$}
;}\newcommand{\nodeh}{\node[draw,circle] (h) {$3$}
;}\newcommand{\nodei}{\node[draw,circle] (i) {$6$}
;}\newcommand{\nodej}{\node[draw,circle] (j) {$7$}
;}\newcommand{\nodeba}{\node[draw,circle] (ba) {$5$}
;}
\begin{tikzpicture}[auto]
\matrix[column sep=.1cm, row sep=0.5cm,ampersand replacement=\&]{
         \&         \&         \&         \& \nodea  \&         \&         \&         \&         \\ 
         \& \nodeb  \&         \& \nodef  \&         \&         \& \nodeg  \&         \& \nodeba \\ 
 \nodec  \& \noded  \& \nodee  \&         \&         \&         \& \nodeh  \&         \&         \\ 
         \&         \&         \&         \&         \& \nodei  \&         \& \nodej  \&         \\
};

\path[ultra thick] (b) edge (c) edge (d) edge (e)
	(h) edge (i) edge (j)
	(g) edge (h)
	(a) edge (b) edge (f) edge (g) edge (ba);
\node[left of=a,node distance=20pt]{$xy^3$};
\node[left of=b,node distance=20pt]{$xy^2$};
\node[left of=g,node distance=20pt]{$x$};
\node[left of=h,node distance=20pt]{$xy$};
\end{tikzpicture}
	    }}%
    \caption{
	    (a) This tree $T\in \scrtprime_{10}$ 
	    is given edge weights 
	    $x,y,y$ for the left, middle and right edges, respectively. 
    \\[2mm]
    (b) This tree $S\in \scri_{10}$ 
    has $\phi_0=1$ and $\phi_i=xy^{i-1}$ 
    for every $i\geq 1$, 
    where $\phi_i$ is the weight of a vertex with $i$ children. 
    The leaves get weight $\phi_0 = 1$.\\[2mm] 
    We have $S = \Phi_{10}(T)$.
    Both trees have weight $x^4 y^6$.
       \\[4mm]
}
    \label{fig.trees}%
\end{figure}
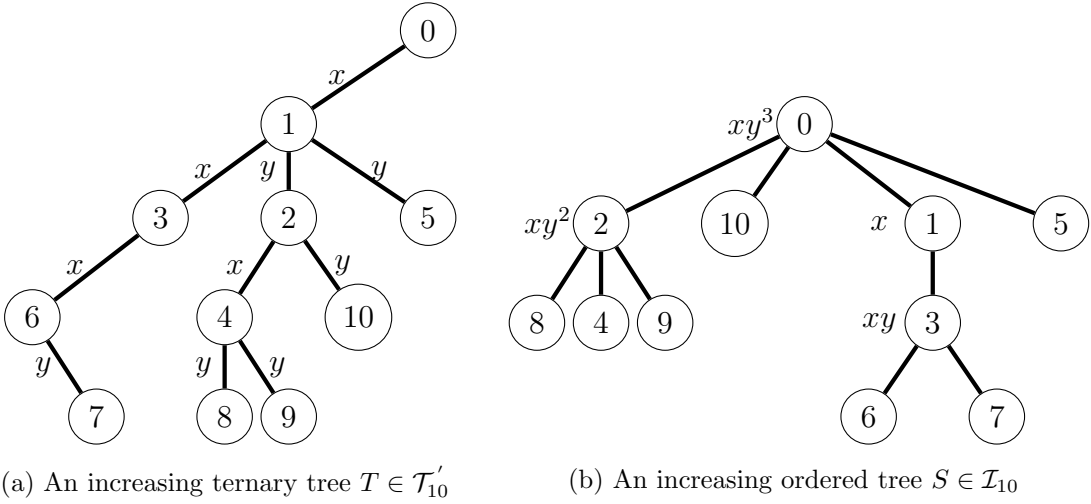

Now let $i$ be a vertex of $T$.
If $i$ has a left child $j$,
we define ${\bf left}\bm{_T(i)}$
to be the middle-right subtree of $j$ in $T$;
if $i$ has no left child, we define ${\rm left}_T(i)$
to be the empty tree.
For example, in Figure~\ref{fig.trees}(a),
the tree ${\rm left}_T(1)$ is the subtree consisting of the vertex $3$,
the tree ${\rm left}_T(2)$ is the subtree
induced by the vertices $4,8$ and $9$,
and the tree ${\rm left}_T(4)$ is the empty tree.

Now let $j$ be any vertex of $T$ other than the root,
and consider the unique path in $T$ from $j$ to the root 0.
Clearly this path must contain at least one left edge,
since the root has only a left edge.
Let ${\rm lanc}(j)$ be the first vertex on this path
that is reached via a left edge;
we refer to ${\rm lanc}(j)$ as the \textbfit{left ancestor} of $j$.
Since the tree $T$ is increasing, we have $0 \le {\rm lanc}(j) < j$.
For example, in Figure~\ref{fig.trees}(a),
${\rm lanc}(1) = 0$, ${\rm lanc}(2) = 0$ and ${\rm lanc}(3) = 1$.

The following lemma is a direct consequence of the definitions of ${\rm left}_T(i)$
and ${\rm lanc}(j)$ and we state it without proof:

\begin{lemma}
   \label{lemma.left_ancestor}
Let $T \in \scrtprime_n$ and let $j$ be a vertex of $T$ different from the root $0$.
Then $j$ is a vertex of ${\rm left}_T({\rm lanc}(j))$
and for every $i\neq {\rm lanc}(j)$, $j$ is not a vertex of ${\rm left}_T(i)$.

\end{lemma}

Next, let $S\in \scri_n$ and let $i\in [n] \cup \{0\}$ be a vertex of $S$.
Then we define the \textbfit{child word of $\bm{i}$},
denoted by ${\rm ch}_S(i)$,
to be the word whose $j$-th letter is the $j$-th child of $i$.
Thus, ${\rm ch}_S(i)$ is the empty word when $i$ is a leaf.
For example, in Figure~\ref{fig.trees}(b), we have
${\rm ch}_S(1) = 3$, ${\rm ch}_S(2) = 849$,
and ${\rm ch}_S(4) = \emptyset$.

We need one more ingredient before describing the map $\Phi_n$.
This is the well-known correspondence between increasing binary trees
on a totally-ordered finite vertex set $A$,
and words $w$ on the alphabet $A$ with no repeated letters:
see e.g.~\cite[pp.~44--45]{Stanley_12}.
However, in our version,
instead of using binary trees with left and right edges,
we will use ternary trees with middle and right edges but no left edges;
the two formulations are of course isomorphic.
The correspondence therefore goes as follows:

If the word $w = \emptyset$, then ${\rm Tree}(w) = \emptyset$.
Otherwise, write $w = u\cdot i \cdot  v$
where $i$ is the smallest letter of $w$;
we then recursively define ${\rm Tree}(w)$ to be the ternary tree with root $i$
whose middle subtree is ${\rm Tree}(u)$ and right subtree is ${\rm Tree}(v)$.

On the other hand, given an increasing ternary tree $T$
on the vertex set $A$ made up solely of middle and right edges,
we define ${\rm Word}(T)$ as follows:
If $T = \emptyset$, then \hbox{${\rm Word}(T)=\emptyset$.}
Otherwise, if $T$ has root $i$ with middle subtree $T'$ 
and right subtree $T''$ (either or both of which may be empty), 
we define ${\rm Word}(T)= {\rm Word}(T')\cdot i \cdot {\rm Word}(T'')$.

From the construction of the maps
$w\mapsto {\rm Tree}(w)$ and $T\mapsto {\rm Word}(T)$,
we clearly have:

\begin{lemma}
\hfill\break\noindent
\vspace*{-7mm}
\begin{itemize}
	\item[(a)] Let $w$ be a word on the alphabet $A$
		with no repeated letters.
		Then\\ $w = {\rm Word}({\rm Tree}(w))$.
	\item[(b)] Let $T$ be an increasing ternary tree
		on the vertex set $A$,
		made up solely of middle and right edges.
		Then $T = {\rm Tree}({\rm Word}(T))$. 
\end{itemize}
	\label{lem.word.tree.inverses}
\end{lemma}

\bigskip

\noindent
\textbfit{Construction of the map $\boldsymbol{\Phi_n:\scrtprime_n\to \scri_n}$.}\\
Given a tree  $T\in \scrtprime_n$,
we will now construct a tree $\Phi_n(T)\in  \scri_n$.
To do this, we begin with a tree containing only the root vertex 0.
We then iterate $j$ from $0$ through $n$
and insert the children of $j$ in $\Phi_n(T)$ as follows:
We consider the subtree ${\rm left}_T(j)$ in $T$;
it consists solely of middle and right edges.
The children of $j$ in $\Phi_n(T)$ are then defined to be
the letters of ${\rm Word}\big({\rm left}_T(j)\big)$,
which are inserted in order from left to right.
That is,
\be
   {\rm ch}_{\Phi_n(T)}(j)
   \;=\;
   {\rm Word}\big({\rm left}_{T}(j)\big)\;.
\label{eq.ch.word.left.phi}
\ee

\noindent
This defines an ordered tree $\Phi_n(T)$ on the vertex set $[n] \cup \{0\}$.
See Figure~\ref{fig.trees} for an example of a tree $T\in \scrtprime_{10}$
and $\Phi_{10}(T)\in\scri_{10}$.

To show that $\Phi_n(T)\in \scri_n$,
we need to show that the tree $\Phi_n(T)$ is increasing,
i.e.~that every vertex $j\in [1,n]$ has a parent in $\Phi_n(T)$
with label strictly smaller than $j$.
We do this by identifying the parent as follows:

\begin{lemma}
   \label{lem.parentofj}
Let $j\in[1,n]$ be a non-root vertex in $T$,
and let ${\rm lanc}(j)$ be its left ancestor.
Then ${\rm lanc}(j)$ is the parent of $j$ in $\Phi_n(T)$.
\end{lemma}

\proof
Let $i$ be the parent of $j$ in $\Phi_n(T)$.
Then, $j$ is a letter in the child word ${\rm ch}_{\Phi_n(T)}(i)$.
From~\reff{eq.ch.word.left.phi},
$j$ must be a vertex of ${\rm left}_T(i).$
Thus, from Lemma~\ref{lemma.left_ancestor} it follows that $i = {\rm lanc}(j)$.
\qed

%

Since ${\rm lanc}(j) < j$, we have shown:

\begin{corollary}
For $T\in \scrtprime_n$, we have $\Phi_n(T)\in \scri_n$.
\end{corollary}

\bigskip

\noindent
{\bf Construction of the map $\Psi_n: \scri_n \to \scrtprime_n$.}\\
Given a tree $S\in \scri_n$,
we will now construct a tree $\Psi_n(S)\in \scrtprime_n$ as follows:
We start with a tree rooted at $0$.
We then iterate through $i\in [n] \cup \{0\}$ and
look at ${\rm ch}_S(i)$, the child word of $i$ in $S$.
If ${\rm ch}_S(i) = \emptyset$, then $i$ will have no left child in $\Psi_n(S)$.
Otherwise, we create a left edge emanating from $i$ in
$\Psi_n(S)$ and attach to it the ternary tree
${\rm Tree}\big({\rm ch}_S(i)\big)$
consisting solely of middle and right edges.
Thus, we have
\be
{\rm left}_{\Psi_n(S)}(i)
\;=\;
{\rm Tree}\big({\rm ch}_S(i)\big)\;.
\label{eq.left.tree.ch.psi}
\ee


\medskip

When $n=0$, we clearly have $\Psi_0(S) \in \scrtprime_0$.
For $n\geq 1$, we will now show that $\Psi_n(S)\in \scrtprime_n$.
In the construction of $\Psi_n(S)$, a vertex $i$ gets a left child $j$
only when ${\rm ch}_S(i)\neq \emptyset$ and 
when $j$ is the root of ${\rm Tree}\big({\rm ch}_S(i)\big)$,
or in other words $j$ is the smallest letter in ${\rm ch}_S(i)$
which is the same as saying that $j$ is 
the smallest-numbered child of the vertex $i$ in $S$.
Thus, it is clear that every vertex $i$ in $\Psi_n(S)$ can contain at most
one left child $j$ and that $j>i$.
Next, notice that a vertex $i$ in $\Psi_n(S)$ 
gets a middle (resp. right) child $j$ only if $i$ and $j$ are both 
letters of the child word ${\rm ch}_S\big(i^\flat\big)$, 
for some vertex~$i^\flat$ in $T$,
and $j$ is the middle (resp. right) child of $i$ in 
\hbox{${\rm Tree}\big({\rm ch}_S\big(i^\flat\big)\big)$}.
Thus, it is clear that every vertex $i$ in $\Psi_n(S)$ can contain at most
one middle child $j$ and at most one right child $j$; furthermore, $j>i$ since 
\hbox{${\rm Tree}\big({\rm ch}_S\big(i^\flat\big)\big)$} is an increasing tree.
Also, $0$ cannot have any middle or right children in $\Psi_n(S)$
as it does not have any siblings in $S$.
This shows that:

\begin{lemma}
 Given $S\in \scri_n$, we have $\Psi_n(T)\in \scrtprime_n$.
\end{lemma}

We are now ready to prove Proposition~\ref{thm.bij.planar.ternary}.

\proofof{Proposition~\ref{thm.bij.planar.ternary} using bijection}
We will first show that the maps $\Phi_n$ and $\Psi_n$ are inverses of each other.
Then we will argue that the number of left edges of a tree $T \in \scrtprime_n$
 equals the number of internal vertices of the tree $\Phi_n(T)\in \scri_n$.

From equations~\reff{eq.ch.word.left.phi} and \reff{eq.left.tree.ch.psi}
and Lemma~\ref{lem.word.tree.inverses},
we observe the following:
\begin{itemize}
	\item[(a)] For a tree $T\in \scrtprime_n$ and a vertex $j \in [n] \cup \{0\}$ of $T$,
		we have 
		\begin{subeqnarray}
		{\rm left}_{\Psi_n(\Phi_n(T))}(j)
		&= &
		{\rm Tree}({{\rm ch}_{\Phi_n(T)}(j)})	
		\\
		&=&
		{\rm Tree}({\rm Word}({\rm left}_{T}(j)))\\
		&=&
		{\rm left}_{T}(j)\;.
		\label{eq.left.psi.phi}
		\end{subeqnarray}
Thus, for the trees $T, \Psi_n(\Phi_n(T))$ every vertex $j\in [n]\cup \{0\}$
satisfies\\
${\rm left}_{\Psi_n(\Phi_n(T))}(j) = {\rm left}_{T}(j)$.
We now need the following lemma:

\begin{lemma}
Let $T_1, T_2 \in \scrtprime_n$.
If ${\rm left}_{T_1}(j) = {\rm left}_{T_2}(j)$ for all $j \in [n] \cup \{0\}$,
then $T_1 = T_2$.
\label{lem.lefttree.same}
\end{lemma}

\proof	
From Lemma~\ref{lemma.left_ancestor}, recall that every vertex $i\neq 0$ 
belongs to exactly one of the subtrees ${\rm left}_{T_1}(j)$
(resp. ${\rm left}_{T_2}(j)$)
in tree $T_1$ (resp. $T_2$).
Since ${\rm left}_{T_1}(j) = {\rm left}_{T_2}(j)$,
a vertex $i$ is the middle (resp. right) child of a vertex 
$i^\flat$ in tree $T_1$
if and only if $i$ is the middle (resp. right) child of
$i^\flat$ in tree $T_2$.
Also, a vertex $i$ is the left child of a vertex 
$i^\flat$ in tree $T_1$ (resp. $T_2$)
if and only if it is the root of the tree 
${\rm left}_{T_1}(i^\flat)$
(resp. ${\rm left}_{T_2}(i^\flat)$).
Thus, 
we also get that $i$ is the left child of
$i^\flat$ in tree $T_1$
if and only if
$i$ is the left child of $i^\flat$ in tree $T_2$.
Thus, $T_1=T_2$.
\qed

Lemma~\ref{lem.lefttree.same} and equation~\reff{eq.left.psi.phi} 
together imply $\Psi_n(\Phi_n(T))=T$.

	\item[(b)] For a tree $S\in \scri_n$ and a vertex $i \in [n] \cup \{0\}$ of $S$,
                we have
		\begin{subeqnarray}
                {\rm ch}_{\Phi_n(\Psi_n(S))}(i)
		&=&
		{\rm Word}({\rm left}_{\Psi_n(S)}(i))
		\\
		&=&
		{\rm Word}({\rm Tree}({\rm ch}_{S}(i)))
		\\
		&=&
		{\rm ch}_{S}(i)
		\label{eq.ch.phi.psi}
		\end{subeqnarray}
		Equation~\reff{eq.ch.phi.psi} gives us
		$\Phi_n(\Psi_n(S)) = S$.
\end{itemize}
Thus, $\Phi_n$ and $\Psi_n$ are inverses of each other.


Finally, notice that a vertex $j$ has a left edge in $T\in \scrtprime_n$
if and only if $j$ has a non-empty child word in the tree $\Phi_n(T)$,
which is true if and only if $i$ is an internal vertex of $\Phi_n(T)$.
Thus, we have $|\scrtprime_{n,k}| = |\scri_{n,k}|$.
\qed

%

%

\medskip

{\bf Remark.}
A different bijection between increasing ternary trees and increasing ordered trees
was given by Janson, Kuba and Panholzer~\cite{Janson_11}.
Our bijection is new and is different from theirs.
We thank Markus Kuba for his comments on this.
\myendremark

\subsubsection{Generating-function proof of 
      Proposition~\ref{thm.bij.planar.ternary}}
\label{sec.interpretation.other.inctree.gf}

We will now use generating functions to provide an alternate proof of 
Proposition~\ref{thm.bij.planar.ternary};
this approach also situates it as a special case of a more general result,
to be enunciated in Proposition~\ref{prop.planar.ternary.gf}.
We begin by introducing a framework for enumerating
increasing ordered trees with very general vertex weights,
which was first introduced by
Bergeron, Flajolet and Salvy \cite{Bergeron_92}
and was further applied in \cite{latpath_lah}.

Let $\bphi = (\phi_i)_{i\geq 0}$ be indeterminates,
and consider increasing ordered trees on $n$ vertices
(here we will use the labels $0,\ldots,n-1$)
in which each vertex with $i$~children gets a weight $\phi_i$;
the weight of the tree, denoted ${\rm wt}(T)$,
is the product of the weights of its vertices.
Let $F(t;\bphi)$ be the weighted exponential generating function of 
such trees, not including the empty tree:
\be
F(t;\bphi) 
\;\eqdef\; 
\sum_{n=1}^\infty \left(\sum_{S\in \scri_{n-1}} {\rm wt}(S) \right) \dfrac{t^n}{n!} \;.
\label{eq.def.F}
\ee
Also let $\Phi(u)$ be the ordinary generating function of $\bphi$:
\be
\Phi(u) \;\eqdef\; \sum_{n=0}^\infty \phi_i u^i\;.
\ee
We call $\Phi(u)$ the \textbfit{degree function}.
It is not difficult to show,
by classifying trees according to the degree of the root,
that $F(t;\bphi)$ satisfies the differential equation 
\cite[eq.~(8)]{Bergeron_92}
\be
\dfrac{d}{dt} F(t;\bphi) \;=\; \Phi( F(t;\bphi))\;.
\label{eq.BLL.trees}
\ee

Notice that even though we defined the generating function
$F(t;\bphi)$ using increasing {\em ordered}\/ trees,
this set-up is very general and we can use different 
degree functions to enumerate different families of weighted trees.
For example, if we substitute $\phi_i = \widetilde{\phi}_i/i!$, 
we may think of $F(t; \bphi)$ as the exponential generating function
enumerating increasing {\em unordered}\/ trees
where each vertex with $i$ children gets weight $\widetilde{\phi}_i$.

It turns out that the connection between
increasing ternary trees and increasing ordered trees
that is asserted in Proposition~\ref{thm.bij.planar.ternary}
is a special case of a more general connection between
increasing trees with degree functions $\Phi$ and $\Psi$
that are related in a particular way, as follows:

\begin{proposition}
	\label{prop.planar.ternary.gf}
Let $\widehat{\Phi}(u)$ be a formal power series with zero constant term,
and let $\scrx(s)$ be a formal power series satisfying
\be
\widehat{\Phi}'(u) \;=\; \scrx(\widehat{\Phi}(u))\;.
\label{eq.planar.ternary.gf.b.1}
\ee
Now define
\begin{subeqnarray}
	\Phi(u) & \eqdef & 1\,+\, x\:\widehat{\Phi}(u)
              \slabel{eq.planar.ternary.gf.a.0}  \\[1mm]
	\Psi(u) &\eqdef& (1 \,+\, x u)\:\scrx(u)
	\slabel{eq.planar.ternary.gf.a.1}
\end{subeqnarray}
Let $\scra(t) \eqdef F(t;\bphi)$ with degree function $\Phi(u)$,
and let $\scrb(t) \eqdef F(t;\bpsi)$ with degree function $\Psi(u)$.
Then the following identities hold:
\begin{subeqnarray}
	\scrb(t)& =& \widehat{\Phi}(\scra(t))
      \slabel{eq.lem.exhibit}  \\[2mm]
	\scra'(t) &=& 1 +  x\scrb(t)
      \slabel{eq.planar.ternary.gf}
\end{subeqnarray}
\end{proposition}

\proof
The definition $\scrb(t) = F(t;\bpsi)$ implies the differential equation
\be
   \scrb'(t) \;=\; \Psi(\scrb(t)) \;=\; (1 \,+\, x \scrb(t))\:\scrx(\scrb(t))
\ee
by \reff{eq.BLL.trees} and \reff{eq.planar.ternary.gf.a.1}.
On the other hand, the definition $\scra(t) = F(t;\bphi)$
implies the differential equation
\begin{subeqnarray}
 \dfrac{d}{dt} \widehat{\Phi}(\scra(t)) 
	&=&
 \widehat{\Phi}'(\scra(t)) \: \scra'(t)
	\\
	&=&
 \scrx(\widehat{\Phi}(\scra(t))) \: \Phi(\scra(t))
	\\[2mm]
	&=&
 \scrx(\widehat{\Phi}(\scra(t))) \: \big[ 1 +  x \widehat{\Phi}(\scra(t)) \big]
\end{subeqnarray}
by the chain rule, \reff{eq.planar.ternary.gf.b.1},
\reff{eq.BLL.trees} and \reff{eq.planar.ternary.gf.a.0}.
So $\scrb(t)$ and $\widehat{\Phi}(\scra(t))$ 
satisfy the same differential equation with the same initial condition
(namely, they vanish at $t=0$), hence they are equal.
This proves \reff{eq.lem.exhibit}.

Then
\be
   \scra'(t)  \;=\;  \Phi(\scra(t))  \;=\;  1 + x \widehat{\Phi}(\scra(t))
       \;=\;  1 +  x\scrb(t)
\ee
by \reff{eq.BLL.trees}, \reff{eq.planar.ternary.gf.a.0}
and \reff{eq.lem.exhibit}.
This proves \reff{eq.planar.ternary.gf}.
\qed

\proofof{Proposition~\ref{thm.bij.planar.ternary} using generating functions}
We will now prove Proposition~\ref{thm.bij.planar.ternary}
by using Proposition~\ref{prop.planar.ternary.gf}.
First, consider the exponential generating function
\hbox{$\scra(t) 
\eqdef 
F\big(t; (\delta_{i0}  + (1-\delta_{i0}) xy^{i-1}))_{i\geq 0} \big)$}:
thus, $\scra(t)$ is the exponential generating function 
for increasing ordered trees on $n$ labeled vertices
where each leaf is given weight $\phi_0 = 1$ 
and each vertex with $i\geq 1$ children is given weight $\phi_i = x y^{i-1}$.
Therefore, the degree function is 
\be
\Phi(u) \;=\; 1\,+\,\dfrac{x u}{1-y u}  \;
\label{eq.degfunc.scra}
\ee
 and
 $\scra(t)$ satisfies the differential equation
 \be
 \scra' 
 \;=\;
 1 + \dfrac{ x \scra}{1 - y\scra}\;.
 \label{eq.scra.diff}
 \ee

In the weight $\phi_i = (x/y) y^i$ for internal vertices,
we can transfer a weight $y$ to each of the $i$ children;
since every vertex other than the root is a child,
we end up with a weight $y$ for leaves and $x$ for internal vertices,
with an overall factor $1/y$ because the root is not a child.
Therefore, from~\reff{eq.def.F} we have
\be
   \scra(t) 
   \;=\;
   \sum_{n=1}^\infty \left(\sum_{S\in \scri_{n-1}}
      x^{{\rm intvertex}(S)} y^{{\rm leaf}(S)-1} \right)
      \dfrac{t^n}{n!}
\label{eq.gf.proof.interpretation.1}
\ee
where ${\rm intvertex}(S)$ and ${\rm leaf}(S)$
denote the number of internal vertices and leaves of
of the tree $S\in \scri_{n-1}$, respectively.

We use the notation from Proposition~\ref{prop.planar.ternary.gf} to obtain
\begin{subeqnarray}
	\widehat{\Phi}(u) &=& \dfrac{u}{1-y u}  \\[1mm]
	\scrx(s) &=& (1+ys)^2  \\[2mm]
	\Psi(u) & = & (1+x u)(1+yu)^2
  \label{eq.Phihat.X.Psi}
\end{subeqnarray}
We then define \hbox{$\scrb(t) \eqdef F(t; \bpsi)$}
with degree function $\Psi(u)  =  (1+x u)(1+yu)^2$.
We can interpret $\scrb(t)$ as
the exponential generating function for increasing ternary trees
on the vertex set $[n]$,
where each left edge gets a weight $x$
and each middle or right edge gets a weight $y$,
and the weight of a tree is the product of the weights of its edges.
Passing to trees
on the vertex set $[n] \cup \{0\}$
in which the only edge (if any) emanating from the root is a left edge,
we get an exponential generating function
\be
   1 +  x\, \scrb(t) 
   \;=\;
\sum_{n=0}^\infty
	\left(\sum_{T\in \scrtprime_n} x^{{\rm l}(T)} y^{{\rm m}(T)+{\rm r}(T)}\right) \dfrac{t^n}{n!}
\label{eq.gf.proof.interpretation.2}
\ee
where ${\rm l}(T), {\rm m}(T), {\rm r}(T)$
denote the number of left, middle and right edges of a tree $T\in \scrt_n$,
respectively.
%
%
%
Finally, we use Proposition~\ref{prop.planar.ternary.gf}
to obtain 
\be
        \scra'(t)
        \;=\; 1 +  x\, \scrb(t)
\label{eq.gf.proof.interpretation.123}
\ee
%
%
Combining equations~\reff{eq.gf.proof.interpretation.1},
\reff{eq.gf.proof.interpretation.2}
and~\reff{eq.gf.proof.interpretation.123}
finishes the proof.
\qed



\medskip

We end this section with a few remarks
concerning some interesting special cases of
Proposition~\ref{prop.planar.ternary.gf}:

\medskip

{\bf Remarks.} 1. If we take $\scrx(s) = 1+ys$
and solve \reff{eq.planar.ternary.gf.b.1},
we get $\widehat{\Phi}(u) \;=\; (e^{yu}-1)/y$ and hence
\be
   \Phi(u)  \;=\;  1 \:+\: {x \over y} \, (e^{yu}-1)
   \;,
\ee
so that $\phi_0 = 1$ and $\phi_i = x y^{i-1}/i!$ for $i \ge 1$.
It follows that $\scra(t)$ is the exponential generating function
for increasing {\em unordered}\/ trees on the vertex set $[n-1]\cup\{0\}$,
with a weight~1 for each leaf
and a weight $x y^{i-1}$ for each vertex with $i \ge 1$ children;
or equivalently (as we saw in the foregoing proof),
with a weight $y$ for leaves and $x$ for internal vertices
and an overall factor $1/y$.
On the other hand,
\be
   \Psi(u)  \;=\;  (1+xu)(1+yu)
   \;,
\ee
so that $\scrb(t)$ is the exponential generating function for
increasing binary trees trees on the vertex set $[n]$,
with a weight $x$ for each left edge and $y$ for each right edge.
Specializing to $x=y=1$,
this shows 
\be
\scrb(t)|_{x=y=1} \;=\; \exp(\scra(t)|_{x=y=1}) - 1 \;=\; \scra'(t)|_{x=y=1}-1,
\ee
i.e., unordered increasing trees on $n+1$ vertices are equinumerous
with increasing binary trees on $n$ vertices,
a fact that is well known and can be proven bijectively
\cite[pp.~44--46]{Stanley_12}.
Furthermore, this also tells us that the left edges in the binary trees
are in correspondence with the internal vertices 
in the corresponding unordered increasing tree.

\medskip

2. If we take $\scrx(s) = (1+ys)^k$ with $k \neq 1$
and solve \reff{eq.planar.ternary.gf.b.1},
we obtain 
\be
\widehat{\Phi}(u) \;=\; \dfrac{[1-(k-1)y u]^{-1/(k-1)}-1}{y}\;.
\label{eq.genk.diff.solve}
\ee
In particular, for integer $k\geq 2$,
the degree function $\Phi(u) = 1 + x\widehat{\Phi}(u)$
obtained from~\reff{eq.genk.diff.solve}
describes increasing ordered trees in which each 
leaf gets a weight $\phi_0 = 1$
and each vertex with $i \ge 1$ children gets a weight
$\phi_i = \dfrac{1}{i!}x y^{i-1} \prod\limits_{j=1}^{i-1} [1 + j(k-1)]$
(in particular, $\phi_1 = x$).
When $k=2$ this reduces to the situation
addressed in Proposition~\ref{thm.bij.planar.ternary}
[cf.~\reff{eq.Phihat.X.Psi}].
Such trees on $n+1$ vertices 
are in bijection with $(k+1)$-ary trees on $n+1$ vertices 
in which the root can have only a first child.

\medskip

3. If we take $\scrx(s) = e^{y s}$ and solve \reff{eq.planar.ternary.gf.b.1},
we obtain 
\be
  \widehat{\Phi}(u) \;=\; -\, {1 \over y} \, \log(1-yu)
\ee
and hence
\be
   \Phi(u)  \;=\;  1 \:-\, {x \over y} \, \log(1-yu)
   \;,
\ee
so that $\phi_0 = 1$ and $\phi_i = x y^{i-1}/i$ for $i \ge 1$.
It follows that $\scra(t)$ is the exponential generating function
for increasing {\em cyclically-ordered}\/ trees on the vertex set $\{0\}\cup[n-1]$,
with a weight 1 for each leaf
and a weight $x y^{i-1}$ for each vertex with $i\ge 1$ children;
or equivalently,
with a weight $y$ for leaves and $x$ for internal vertices
and an overall factor $1/y$.
On the other hand,
\be
   \Psi(u)  \;=\;  (1+xu) \, e^{yu}
   \;,
\ee
so that $\scrb(t)$ is the exponential generating function for
increasing trees on vertex set $[n]$
in which every vertex has zero or one left child (with weight $x$)
and an arbitrary number of unordered right children (with weight $y$).
Therefore, $1 + x\scrb(t)$ is the exponential generating function for
increasing trees of the same type on vertex set $\{0\}\cup [n]$,
but in which the root $0$ can only have a left child.
It follows that these trees
are in bijection with cyclically-ordered trees on the vertex set $\{0\}\cup [n]$.


\medskip

 \begin{problem} Find bijections that explain Remarks~2 and 3.
 \end{problem}

4. On specializing to $y=1+x$ in the degree function~\reff{eq.degfunc.scra},
we obtain the weights $\phi_0 = 1$ and $\phi_i = x(1+x)^{i-1}$ for $i \ge 1$.
This sequence is not coefficientwise Toeplitz-totally positive in the 
variable $x$
--- or even Toeplitz-totally positive for any fixed real number $x > 0$ ---
since its generating function
\be
   \Phi(u)  \;=\;  1 \,+\, \sum_{i=1}^\infty x (1+x)^{i-1} u^i
     \;=\;
     {1-u \over 1 \,-\, (1+x)u}
\ee
has a zero on the positive real axis at $u = 1$,
which cannot happen for a Toeplitz-totally positive sequence
\cite[p.~412, Theorem~5.3]{Karlin_68}.
However, by combining 
Proposition~\ref{thm.bij.planar.ternary}
with equation~\reff{eq.ternary},
we deduce that the generating polynomials of 
increasing ordered trees in $\scri_n$ with the weights $\bphi$
are the second-order Stirling cycle polynomials $c_{2,n}(x)$;
and from Theorem~\ref{thm.second.HankelTP}
we see that the sequence $c_{2,n}(x)$ is 
coefficientwise Hankel-totally positive in the variable $x$.

This example has an important consequence for the understanding of
\cite[Theorem~1.1]{latpath_lah}
concerning the Hankel-total positivity of the generic Lah polynomials.
We see that, although this sequence $\bphi$ is not Toeplitz-totally positive,
the resulting sequence of generating polynomials
in \cite[Theorem~1.1(c)]{latpath_lah}
is still coefficientwise Hankel-totally positive.
Thus, the condition of $\bphi$ being Toeplitz-totally positive
is only a sufficient condition and is far from necessary for
the conclusion of \cite[Theorem~1.1(c)]{latpath_lah} to hold.
In \cite{Lah_experiment},
we take $\Phi(u)$ is to be a polynomial in $u$ of degree $3$ or $4$,
and we study the exact necessary and sufficient conditions under which 
the conclusion of \cite[Theorem~1.1(c)]{latpath_lah} holds true.
\myendremark

\subsubsection{Vertex-marked increasing ordered trees}
\label{subsec.vertex.marked}

Having established Proposition~\ref{thm.bij.planar.ternary},
we are now ready to state our fourth interpretation for the second-order Stirling cycle numbers 
$\stirlingcycle{n}{k}^{(2)}$.

We define a \textbfit{vertex-marked increasing ordered tree}
to be an increasing ordered tree in which each leaf can be
marked or unmarked; the internal vertices are always unmarked.
For $n=0$, we declare that the vertex $0$ is always marked.
We write $\widehat{\scri}_{n,k}$ for the set of
vertex-marked increasing ordered trees on the vertex set $[n] \cup \{0\}$
with $k$ unmarked vertices.
In this section we will show that the 
the cardinality of the set $\widehat{\scri}_{n,k}$ is 
$\stirlingcycle{n}{k}^{(2)}$.

Recall also from Section~\ref{sec.interpretation.ternary} that $\widehat{\scrt}_{n,k}$ denotes
the set of trees in $T\in\scrtprime_n$
along with optional edge markings in which each edge can be either marked or unmarked
such that left edges are always unmarked,
and which contain $k$ unmarked edges.

Before providing the proofs, let us observe that using the bijection 
and the correspondence between left edges and internal vertices
in Proposition~\ref{thm.bij.planar.ternary},
and we get the following corollary:

\begin{corollary}
        The sets $\widehat{\scrt}_{n,k}$ and $\widehat{\scri}_{n,k}$
        are in bijective correspondence.
        Thus, the cardinality of $\widehat{\scri}_{n,k}$ is given by
        \be
        |\widehat{\scri}_{n,k}| \;=\; \stirlingcycle{n}{k}^{(2)}\;.
	\label{eq.eq.cor.scrt.scri.bijection}
        \ee
	\label{cor.scrt.scri.bijection}
\end{corollary}

\proof
First notice that for $n=0$, we have $|\scrthat_{0,0}| = |\widehat{\scri}_{0,0}|=1$.
Henceforth, we assume $n\geq 1$.
For a tree $S\in \scri_{n}$,
let ${\rm intvertex}(S)$ denote the number of internal vertices of $S$
and let ${\rm leaf}(S)$ denote the number of leaves of $S$.
Also, for a tree $T\in \scrtprime_{n}$,
let ${\rm l}(T)$, ${\rm m}(T)$, ${\rm r}(T)$ 
denote the number of left, middle and right edges of $T$, respectively.
Then we clearly obtain the following
\begin{subeqnarray}
\sum_{k=0}^{n} |\widehat{\scri}_{n,k}|\:  x^k
	&=&
\sum_{S\in \scri_{n}} x^{{\rm intvertex}(S)} \: (1+x)^{{\rm leaf}(S)}
	\slabel{eq.cor.scrt.scri.bijection.a}\\
	&=&
\sum_{T\in \scrtprime_{n}} x^{{\rm l}(T)} \: (1+x)^{{\rm m}(T)} \: (1+x)^{{\rm r}(T)}
\qquad \qquad  \slabel{eq.cor.scrt.scri.bijection.b}\\
	&=&
\sum_{k=0}^{n} |\widehat{\scrt}_{n,k}|\:  x^k
      \slabel{eq.cor.scrt.scri.bijection.c}
	\label{eq.cor.scrt.scri.bijection}
\end{subeqnarray}
where~\reff{eq.cor.scrt.scri.bijection.a} follows from the definition of the sets 
$\scri_n,\widehat{\scri}_{n,k}$,
\reff{eq.cor.scrt.scri.bijection.b} follows from Proposition~\ref{thm.bij.planar.ternary},
and~\reff{eq.cor.scrt.scri.bijection.c} follows from the definition of 
the sets~$\scrtprime_{n}, \widehat{\scrt}_{n,k}$.
This completes the proof.
\qed

We thus have
our fourth combinatorial interpretation for the numbers
$\stirlingcycle{n}{k}^{(2)}$.

\subsection{Interpretation V: Cyclically-ordered phylogenetic trees}
\label{sec.interpretation.other.phylogenetic}

Here we introduce our fifth and final combinatorial interpretation of the
second-order Stirling cycle numbers.
This interpretation will lead us to a common generalization  
of the second-order Stirling numbers of both kinds
and also lead us to another conjecture
similar to Conjectures~\ref{conj.cycle} and~\ref{conj.subset}.


A \textbfit{phylogenetic tree} is a leaf-labeled rooted tree
(i.e.~the leaves are labeled and the non-leaf vertices are unlabeled)
in which all non-leaf vertices have at least two children.
We refer to them as an \textbfit{unordered phylogenetic tree},
an \textbfit{ordered phylogenetic tree},
or a \textbfit{cyclically-ordered phylogenetic tree}
according to whether the children at each vertex are
unordered, linearly ordered or cyclically ordered.
Cyclically-ordered phylogenetic trees have also been called
\textbfit{series-reduced mobiles}
\cite[A032188]{OEIS} \cite{Seqfan_18092011}.
Let $\scrm_{n,k}$ denote the set of all cyclically-ordered phylogenetic trees
with $n+1$ leaves and $k$ internal (i.e.~non-leaf) vertices.


%
%

We will now establish the following proposition:
\begin{proposition}
The cardinality of the set $\scrm_{n,k}$ of 
cyclically-ordered phylogenetic trees 
with
$n+1$ leaves and $k$ internal vertices is given by
\be
	\left|\scrm_{n,k} \right| \;=\; \stirlingcycle{n}{k}^{(2)}.
\ee
\label{prop.seriesredmobile}
\end{proposition}

We will provide two different proofs of Proposition~\ref{prop.seriesredmobile}.
Our first proof will show that the numbers $\left|\scrm_{n,k} \right|$
satisfy the recurrence relation~\reff{eq.rec.stirlingcycle.2nd}.
Our second proof will use generating functions
for the multivariate Ward polynomials
introduced by Elvey Price and Sokal in \cite{Elvey-Price-Sokal_wardpoly}.

\proofof{Proposition~\ref{prop.seriesredmobile} using recurrence relations}
We will show that the numbers  $\left|\scrm_{n,k} \right|$ satisfy
\be
   |\scrm_{n,k}|
   \;=\;
   (n+k-1) \, |\scrm_{n-1,k-1}| \:+\: (n+k-1) \, |\scrm_{n-1,k}|
   \;.
\ee
To do this, we will consider two cases whether or not the leaf labeled $n+1$
has exactly one sibling.

If  $n+1$ has at least two siblings, 
we can remove it, and the resulting tree belongs to $\scrm_{n-1,k}$.
Conversely, to insert $n+1$ into a tree belonging to $\scrm_{n-1,k}$,
it suffices to choose the vertex that would be its immediately preceding
(in the cyclic order) sibling in the resulting tree.
There are $n+k-1$ such choices.
This contributes the term $(n+k-1) \, |\scrm_{n-1,k}|$. 

If $n+1$ has exactly one sibling, 
then the tree obtained by removing it is no longer a phylogenetic tree,
as the parent of $n+1$ would now have only one child. 
Instead, we proceed as follows:
\begin{itemize}
   \item
If the parent of $n+1$ is not the root,
we contract the edge between the parent and grandparent of $n+1$
and then remove $n+1$ from the resulting tree;
this gives us a tree in $\scrm_{n-1,k-1}$.
Conversely, to insert $n+1$ back into this tree,
we need to select an edge,
and elongate it by adding an extra internal vertex 
and then insert $n+1$ as the other child of this new internal vertex. 
As there are $n+k-2$ edges,
we get a contribution of $(n+k-2)|\scrm_{n-1,k-1}|$. 
   \item
If the parent of $n+1$ is the root, 
we remove the root and the vertex $n+1$
to obtain a tree rooted at the sibling of $n+1$.
Conversely, we can re-obtain the original tree
by creating a new root whose children are the old root and vertex $n+1$.
This contributes $|\scrm_{n-1,k-1}|$.
\end{itemize}
Putting these two contributions together gives $(n+k-1) \, |\scrm_{n-1,k-1}|$.
\qed

For our second proof of Proposition~\ref{prop.seriesredmobile},
we require the multivariate Ward polynomials 
introduced in \cite[pp.~9--11]{Elvey-Price-Sokal_wardpoly}.
As mentioned in Section~\ref{sec.conjecture},  
our second-order Stirling subset polynomials $s_{2,n}(x)$
are usually referred to as the Ward polynomials,
and the numbers $\stirlingsubset{n}{k}^{(2)}$ are the Ward numbers.
It is known \cite[p.~5]{Elvey-Price-Sokal_wardpoly}
that the Ward number $\stirlingsubset{n}{k}^{(2)}$ 
counts the number of phylogenetic trees (with unordered children)
on $n+1$ leaves and $k$ internal vertices.
Thus, the polynomial $s_{2,n}(x)$ is the generating polynomial 
of phylogenetic trees on $n+1$ leaves in which 
each internal vertex gets a weight $x$.
Now let $\bfx = (x_i)_{i \ge 1}$ be an infinite collection of indeterminates.
The multivariate Ward polynomials $\bW_n(\bfx) = \bW_n(x_1,\ldots,x_n)$
are a multivariate generalization of the polynomials $s_{2,n}(x)$:
they are the generating polynomial for phylogenetic trees
on $n+1$ labeled leaves
in which each internal vertex with $i$ ($\ge 2$) children
gets a weight~$x_{i-1}$.
Thus $s_{2,n}(x) = \bW_n(x,\ldots,x)$. 
The polynomial $\bW_n(\bfx)$ is quasi-homogeneous of degree~$n$ 
when $x_i$ is given weight~$i$.
The first few $\bW_n$ are \cite[A134685]{OEIS}
\begin{subeqnarray}
   \bW_0  & = &  1  \\
   \bW_1  & = &  x_1  \\
   \bW_2  & = &  3 x_1^2 + x_2  \\
   \bW_3  & = &  15 x_1^3 + 10 x_1 x_2 + x_3  \\
   \bW_4  & = &  105 x_1^4 + 105 x_1^2 x_2 + 15 x_1 x_3 + 10 x_2^2 + x_4 
\label{eq.multiward.0-4}
\end{subeqnarray}
It is not difficult to show that the exponential generating function 
\be
   \scrw(t;\bfx)
   \;\eqdef\;
   \sum_{n=0}^\infty \bW_n(\bfx) \, {t^{n+1} \over (n+1)!}
   \;=\;
   t \,+\, \sum_{n=2}^\infty \bW_{n-1}(\bfx) \, {t^n \over n!}
 \label{def.scrw.phylogenetic}
\ee
satisfies the functional equation~\cite[eq.~(1.18)]{Elvey-Price-Sokal_wardpoly}
\be
   \scrw(t;\bfx)
   \;=\;
   t \,+\, \sum_{j=2}^\infty x_{j-1} \, {\scrw(t;\bfx)^j \over j!}
   \;,
 \label{eq.functeqn.scrw}
\ee
where the term $j$ in the sum corresponds to the case
in which the root has $j$ children.

\proofof{Proposition~\ref{prop.seriesredmobile} using generating functions}
With the specialization $x_i =  i! \, x$,
the polynomials $\bW_n(x, 2 x, \ldots, n!\, x)$ are the generating polynomials 
of cyclically-ordered phylogenetic trees with weight $x$ given 
to each internal vertex and a weight $1$ to each leaf.
We therefore want to show that
$\bW_n(x, 2 x, \ldots, n!\, x) = c_{2,n}(x)$.
We call these polynomials
$W_n(x)\eqdef \bW_n(x, 2 x, \ldots, n!\, x)$,
and we denote the corresponding exponential generating function
by $\scrw(t;x)$.
Thus, the functional equation~\reff{eq.functeqn.scrw} is now
\be
\scrw(t;x)
   \;=\;
   t \,+\, \sum_{j=2}^\infty  \scrw(t;x)^j \, \dfrac{x}{j} 
   \;.
\ee
This can be rewritten as
\be
   \scrw(t;x) 
   \;=\;
   t  \:-\:  x \, \log[1-\scrw(t;x)]  \:-\: x \, \scrw(t;x)
   \;.
\ee
Differentiating both sides with respect to $t$ and solving for $\scrw'$,
we obtain
\be
\scrw' \;=\; \dfrac{1}{1 - x\: \dfrac{\scrw}{1-\scrw}}
\;=\; 
1+ \dfrac{x\:\scrw}{1-(1+x)\:\scrw}
\;.
\ee
However, this is exactly the same as 
the differential equation~\reff{eq.scra.diff} with $y=1+x$ 
which is satisfied by the exponential generating function
for increasing ordered trees on $[n] \cup \{0\}$
where each leaf is given a  weight $1$ 
and each vertex with $i\geq 1$ children is given a weight
$x(1+x)^{i-1}$.
Thus, the result follows by combining
Proposition~\ref{thm.bij.planar.ternary} and equation~\reff{eq.ternary}.
\qed

%
%

{\bf Final remark and conjecture.}
As mentioned earlier, it was noticed in \cite{Elvey-Price-Sokal_wardpoly}
that the number of phylogenetic trees
with {\em unordered}\/ children having $n+1$ leaves and $k$ internal vertices
is the second-order Stirling subset number $\stirlingsubset{n}{k}^{(2)}$.
In Proposition~\ref{prop.seriesredmobile}
we established that the number of phylogenetic trees
with {\em cyclically}\/-ordered children
having $n+1$ leaves and $k$ internal vertices
is the second-order Stirling cycle number $\stirlingcycle{n}{k}^{(2)}$.
The generating polynomials in both cases can be obtained as 
specializations of the multivariate Ward polynomials:
the former by setting $x_i = x$,
the latter by setting $x_i = i! \, x$.
After these two observations, 
it is very natural to inquire about the number
of phylogenetic trees with {\em linearly}\/-ordered children
having $n+1$ leaves and $k$ internal vertices.
Let us denote this number by $\stirlingcyclesub{n}{k}$.
Their generating polynomials are obtained by setting $x_i = (i+1)! \, x$
in the multivariate Ward polynomials.
The first few of these numbers are

\bigskip

%
%
\begin{table}[H]
\centering
\small
\begin{tabular}{c|rrrrrrr|r}
$n \setminus k$ & 0 & 1 & 2 & 3 & 4 & 5 & 6 & Row sums \\
\hline
0 & 1 &  &  &  &  &  &  & 1  \\
1 & 0 & 2 &  &  &  &  &  & 2  \\
2 & 0 & 6 & 12 &  &  &  &  & 18  \\
3 & 0 & 24 & 120 & 120 &  &  &  & 264  \\
4 & 0 & 120 & 1080 & 2520 & 1680 &  &  & 5400  \\
5 & 0 & 720 & 10080 & 40320 & 60480 & 30240 &  & 141840  \\
6 & 0 & 5040 & 100800 & 604800 & 1512000 & 1663200 & 665280 & 4551120  \\
\end{tabular}
\end{table}

\bigskip

\noindent
and they seem to match \cite[A357367]{OEIS}, namely\footnote{
   This explicit formula is not yet shown at \cite[A357367]{OEIS}.
}
\be
   \stirlingcyclesub{n}{k} \;=\; \dfrac{(n+k)!}{k!} \, \binom{n-1}{k-1}  \;.
  \label{eq.prop.stirlingcyclesub.phy}
\ee
%
%
This is, in fact, true:

\begin{proposition}
The cardinality of the set of all phylogenetic trees
with linearly ordered children
having $n+1$ leaves and $k$ internal vertices
is given by \reff{eq.prop.stirlingcyclesub.phy}.
\end{proposition}

\proof
We show that the numbers $\stirlingcyclesub{n}{k}$ satisfy the recurrence
\be
\stirlingcyclesub{n}{k} 
\;=\; 
(2n+2k-2) \, \stirlingcyclesub{n-1}{k-1}
\:+\:
(n+2k-1)  \, \stirlingcyclesub{n-1}{k}\;.
\ee
Our argument is a modification of the proof of
Proposition~\ref{prop.seriesredmobile} using recurrence relations,
and we leave the details to the reader.
Equation~\reff{eq.prop.stirlingcyclesub.phy} then follows,
as it is easy to see that the numbers $\big[(n+k)!/k!\big] \binom{n-1}{k-1}$
also satisfy the same recurrence.
\qed

Let $d_n(x)$ denote the generating polynomials
\be
d_n(x) \;\eqdef \; \sum_{k=0}^{n} \stirlingcyclesub{n}{k} x^k \;,
\ee
and let $D$ be the lower-triangular matrix 
$D = \left(\stirlingcyclesub{n}{k}\right)_{n,k\geq 0}$.

Similar to Conjectures~\ref{conj.cycle} and~\ref{conj.subset},
we have the following conjecture:
\begin{conjecture}[Conjecture for ordered phylogenetic trees]
\label{conj.phylogenetic.ordered}
The following are true:
\begin{enumerate}
\item The triangle $D$ is totally positive.

\item The triangle $D^{\rm rev}$ formed by reversing the rows of 
	$D$ is totally positive.


\item  The sequence $(d_n(x))_{n\geq 0}$ is 
	coefficientwise Hankel-totally positive in the variable $x$.

\end{enumerate}
\end{conjecture}

The polynomials $d_n(x)$ are also negative-real-rooted,
and in this case we are able to prove it:

\begin{theorem}
   \label{thm.phylogenetic.ordered.zeros}
For every $n\geq 1$, the polynomials $d_{n}(x)$ have $0$ as a root;
and the zeros of the polynomial sequence $\big(d_{n}(x)/x \big)_{n\geq 1}$
are simple, interlacing, and lie in the interval $(-1, 0)$.
\end{theorem}

The proof of Theorem~\ref{thm.phylogenetic.ordered.zeros}
is similar to that of Theorem~\ref{thm.stirling.2nd.zeros}
(which we will see in Section~\ref{sec.zeros.2nd})
and we omit it.
\myendremark

We leave the further study of the multivariate Ward polynomials
as an open problem.
In particular, it would be interesting to find sufficient conditions
for Hankel-total positivity,
analogous to \cite[Theorem~1.1]{latpath_lah}.

\subsection{A final remark}  \label{subsec.remark}

The diagonal entries of the second-order Stirling cycle matrix $C^{(2)}$
are ${(2n-1)!!}$ [cf.~\reff{eq.stir.diagonal.r=2}],
and the reader will observe that the first four of our
combinatorial interpretations of the second-order Stirling cycle numbers
are based on taking well-known \cite{Callan_09} objects
of cardinality $(2n-1)!!$ and then ``augmenting'' them:
\begin{itemize}
   \item Perfect matchings of $[2n]$ $\longrightarrow$
         Derangements of $[n+k]$ containing $k$ cycles
   \item Stirling permutations of $\{1,1,2,2,\ldots,n,n\}$  $\longrightarrow$
         Ascent-marked Stirling permutations of $\{1,1,2,2,\ldots,n,n\}$
   \item Increasing ternary trees on $[n]$  $\longrightarrow$
         Edge-marked increasing ternary trees on $[n]$
   \item Increasing ordered trees on $[n] \cup \{0\}$  $\longrightarrow$
         Edge-marked increasing ordered trees on $[n] \cup \{0\}$
\end{itemize}

In fact, also the second-order Stirling subset matrix $S^{(2)}$
has diagonal entries $(2n-1)!!$ [cf.~again \reff{eq.stir.diagonal.r=2}],
and the construction of Elvey Price and Sokal \cite{Elvey-Price-Sokal_wardpoly}
was also based on a similar ``augmentation'':
\begin{itemize}
   \item Perfect matchings of $[2n]$ $\longrightarrow$
   Augmented and super-augmented perfect matchings of $[2n]$
\end{itemize}

We wonder whether the $r$-th order Stirling cycle and subset polynomials
for $r \ge 3$ can be handled by a similar approach:
namely, find a class of combinatorial objects enumerated by
the numbers \reff{eq.stircyc.diagonal} or \reff{eq.stirsub.diagonal},
and then ``augment'' them to handle the matrix elements with $k < n$.
We will use this approach to construct an analogue of Stirling permutations for $r\ge 3$
in Section~\ref{sec.interpretation.stirperm.gen};
this will give us a new interpretaton for the $r$-th order Stirling cycle numbers.

\section{Interpretation for $\bm{r}$-th order Stirling cycle numbers}
\label{sec.interpretation.stirperm.gen}

We will now provide an interpretation for the 
$r$-th order Stirling cycle numbers that generalizes the interpretation
given for $r=2$ in Section~\ref{sec.interpretation.stirperm}.
We will also provide an interpretation for the
$r$th-order quasi-Eulerian cycle numbers
$\eulersecond{n}{k}^{\! (r)}_{\!\rm C}$.

To begin, let us first provide an interpretation for the numbers 
$\stirlingcycle{n}{n}^{(r)} = (rn)!/(r^n n!)$ for $r\geq 2$.
To do this, let us first consider the multiset 
\be
M_{n}^{(r)} \eqdef \{1^{a_1}, 2^{a_2}, \ldots, (r-1)n^{a_{(r-1)n}}\}
\ee
where $a_i$ is the multiplicity of the letter $i$ and it is defined by
\be
a_i
\;=\;
\begin{cases}
2 \qquad \text{if $i\equiv 0 \pmod{r-1}$}\\
1 \qquad \text{otherwise}
\end{cases}
\ee
Thus, for $r=2$ we have 
$M_n^{(2)} = \{1,1,2,2,3,3,\ldots, n,n\}$
and for $r=3$ we have
$M_n^{(3)} = \{1,2,2,3,4,4,\ldots,2n-1 ,2n,2n\}$.
Let $\Sym_{n}^{(r)}$ denote the set of Stirling permutations of the 
multiset $M_{n}^{(r)}$.
It immediately follows from \cite[eq.~(1)]{Janson_11} that
\be
\stirlingcycle{n}{n}^{(r)} \;=\; \left| \Sym_{n}^{(r)} \right| \,.
\ee
The set $\Sym_{n}^{(r)}$ is our first ingredient.
Notice that each $\sigma\in \Sym_{n}^{(r)}$
is a string of length $rn$ consisting of $(r-1)n$ distinct letters.


Our next ingredient is a new combinatorial statistic on 
Stirling permutations in $\Sym_{n}^{(r)}$.
Given $\sigma\in \Sym_n^{(r)}$,
we say that $i\in [n]$ is a \textbfit{consecutive ascent} in $\sigma$
if the following two conditions hold:
\begin{itemize}
	\item if the following $r-1$ letters:
		$(r-1)i+1, (r-1)i + 2, \ldots, (r-1)i +(r-2)$
		and the first occurence of $(r-1)(i+1)$, 
		occur consecutively in $\sigma$.

	\item the letter $(r-1)i+1$ is preceeded by a lower-numbered letter in $\sigma$.
\end{itemize}
An alternate way of saying this is that the string $a_1a_2\ldots a_{r-1}$
where $a_j =(r-1)i+j$ occurs as a substring in $\sigma$ 
and $a_1$ is preceeded by a lower-numbered letter.
Thus, if $i \in [n]$ is a consecutive ascent in $\Sym_n^{(r)}$,
then the $r-1$ letters,
$(r-1)i+1, (r-1)i + 2, \ldots, (r-1)i +(r-2)$                 
and the first occurence of $(r-1)(i+1)$,
are the ends of ascents (i.e. $a$ in $\sigma_i<\sigma_{i+1}=a$).
Let $\CAsc(\sigma) \subseteq [n]$ denote the set of consecutive ascents
in~$\sigma$.
We then define $\Sym_{n,k}^{(r)}$ to be the set of
Stirling permutations in $\Sym_{n}^{(r)}$ with $k$ consecutive ascents:
\be
\Sym_{n,k}^{(r)}
\;\eqdef\;
\{\sigma \in \Sym_{n}^{(r)} \:\colon \CAsc(\sigma) = k\}\;.
\ee
Notice that when $r=2$,
consecutive ascents are the same as ascents in Stirling permutations,
so $|\Sym_{n,k}^{(2)}| = \eulersecond{n}{k}$.


A \textbfit{consecutive-ascent-marked Stirling permutation 
of order $\bm{r}$ and size $\bm{n}$ with $\bm{k}$ dots}
is defined to be a pair $(\sigma, D)$ where 
$\sigma\in \Sym_n^{(r)}$ and $D\subseteq \CAsc(\sigma) \subseteq [n]$
and $|D| = k$.
We let $\scrv_{n,k}^{(r)}$ denote the set of
consecutive-ascent-marked Stirling permutation
of order $r$ and size $n$ with $k$ dots.
Thus, when $r=2$ we have the interpretation in 
Section~\ref{sec.interpretation.stirperm}.
For $r=3$, we have $\scrv_{2,0}^{(3)}$ is equivalent to $\Sym_2^{(3)}$
(as we always have $D = \emptyset$)
which has cardinality $40$ and hence we refrain from listing out.
For $k=1$, we have 
\be
	\scrv_{2,1}^{(3)} \; = \;  \{1 22 \bolddot 3 \bolddot 44, 
				1 2 \bolddot 3\bolddot 4 4 2,
				1\bolddot 3\bolddot 4 4 2 2\}.
\ee
Here we have inserted dots to the $r-1 = 2$ ascents corresponding
to each of the $i\in D$.

Let $v_{n,k}^{(r)} = | \scrv_{n,k}^{(r)} |$
and $b_{n,k}^{(r)} = | \Sym_{n,k}^{(r)} |$;
then
\be
v_{n,k}^{(r)}
\;=\;
\sum_{i=k}^{n} b_{n,i}^{(r)} \: \binom{i}{k}\;.
\label{eq.vnkr.eulersecond}
\ee
Notice that for $r=2$, the identity $a_{n,k}^{(2)} =\eulersecond{n}{k}$ together with
Proposition~\ref{prop.secondeuler}(a), which gives us~\reff{eq.secondeuler.binom.cycle}, 
implies $v_{n,k}^{(2)} = \stirlingcycle{n}{n-k}^{(2)}$.

We will now establish the following lemma,
which generalises Lemma~\ref{lem.perm.cyc.asc.marked}:

\begin{lemma}
   \label{lem.perm.cyc.asc.marked.gen}
The numbers $v_{n,k}^{(r)}$ satisfy the recurrence
\begin{eqnarray}
v_{n,k}^{(r)} 
	&=& 
	\;(r-1)!\binom{r(n-1) \,-\, (r-1)(k-1)}{r-1}\, v_{n-1,k}^{(r)}
	\nonumber\\[3mm]
	&&
	+\;\; 
	[r(n-1) \,-\, (r-1)(k-1)]\, v_{n-1,k-1}^{(r)} 
     \label{eq.rec.vnkr}
\end{eqnarray}
with initial conditions $v_{0,k}^{r} = \delta_{k0}$ and $v_{n,-1}^{r}=0$.
\end{lemma}


\proof 
Let $(\sigma,D) \in \scrv_{n,k}^{(r)}$ be a
consecutive-ascent-marked Stirling permutation
of order $r$ and size $n$ with $k$~dots.
Throughout this proof, 
we set $a_j = (r-1)(n-1)+j$ where $1\leq j \leq r-1$
and $a_{r} = (r-1)n$.

First, notice that on deleting the $r-2$ letters 
$a_1, \ldots, a_{r-2}$
along with the two letters $a_{r-1}, a_r$ in $\sigma$
together with $n\in D$, if it exists,
gives an element of $\scrv_{n-1,k}^{(r)}$ if $n\not \in D$,
and an element of $\scrv_{n-1,k}^{(r)}$ if $n\in D$.

Conversely, to get an element $(\sigma,D) \in \scrv_{n,k}^{(r)}$ 
with $n\not \in D$ from a one-size-smaller 
ascent-marked Stirling permutation, 
we insert the $r$ letters 
$a_1, \ldots, a_{r-2},\,\, a_{r-1},a_{r}$.
For convenience, let $l = r(n-1)-(r-1)k$.
Notice that there are $l+1$ choices to insert the letter $a_1$ 
(these consist of any of the $r(n-1) +1$ spaces 
except for the $(r-1)k$ corresponding to dots),
$l+2$ to insert the letter $a_2$,
and so on until we finally get  
$l+ (r-1)$ choices for inserting $a_{r-1}\,a_{r}$.
Thus, there are a total of $(r-1)! \binom{l+r-1}{r-1}$ choices.

In the other case, if $(\sigma,D) \in \scrv_{n,k}^{(r)}$ with $n\in D$,
then the string 
$a_1 a_2 \ldots a_{r-1}$
is a substring of $\sigma$.
To get such a pair $(\sigma, D)$, 
the string $a_1 a_2 \ldots a_{r-1} a_{r}$
can be inserted into an element of $\scrv_{n-1,k-1}^{(r)}$
in any of the $r(n-1)+1$ spaces except the first space
and those corresponding to a dot; 
thus, the number of choices is $r(n-1) - (r-1)(k-1)$.

This proves the recurrence~\reff{eq.rec.vnkr}.
\qed

Finally, we use Lemma~\ref{lem.perm.cyc.asc.marked.gen}
to obtain
our combinatorial interpretation for the $r$th-order Stirling cycle numbers and 
for the $r$th-order quasi-Eulerian cycle numbers;
we state this in the following  proposition:

\begin{proposition}
The entries of the matrix $Q_{\rm C}^{(r)} = \widecheck{C}^{(r)} B^{-1}$
of $r$th-order quasi-Eulerian cycle numbers
have the combinatorial interpretation
\be
\eulersecond{n}{k}^{\!(r)}_{\!\rm C} \;=\; |\Sym_{n,k}^{(r)}|
   \;.
\ee
Also, the $r$th-order Stirling cycle number $\stirlingcycle{n}{k}^{(r)}$
counts the number of
consecutive-ascent-marked Stirling permutations of order $r$ and size $n$ 
with $n-k$ dots, 
i.e., $\stirlingcycle{n}{k}^{(r)} = |\scrv_{n,n-k}^{(r)}|$.
\label{prop.combi.interpret.rth}
\end{proposition}

\proof
We now need to show that
$\stirlingcycle{n}{k}^{(r)} = v_{n,n-k}^{(r)}$ and
$\eulersecond{n}{k}^{(r)}_{\!\rm C} = b_{n,k}^{(r)}$.

From \reff{eq.rec.vnkr} it follows that the numbers
$u_{n,k}^{(r)} = v_{n,n-k}^{(r)}$ satisfies
the recurrence~\reff{rec.cycle}
which is satisfied by the numbers $\stirlingcycle{n}{k}^{(r)}$.
Since these numbers also satisfy the same initial conditions,
we infer that $\stirlingcycle{n}{k}^{(r)} = u_{n,k}^{(r)}$.
Finally, from the definition of $\eulersecond{n}{k}^{(r)}_{\!\rm C}$
and~\reff{eq.vnkr.eulersecond}, it follows that~$\eulersecond{n}{k}^{(r)}_{\!\rm C} = b_{n,k}^{(r)}$.
\qed

\section{Location of zeros for the second-order Stirling cycle and subset polynomials}
\label{sec.zeros.2nd}

We will now prove Theorem~\ref{thm.stirling.2nd.zeros}
in this section.

For $r,n\geq 1$, it is clear from~\reff{eq.cycle.base}/\reff{eq.subset.base}
that $c_{r,n}(0) = s_{r,n}(0) = 0$.
Therefore, for $n\geq 0$ and $r=2$ we define the polynomials
$\widehat{c}_{n}(x)$ and $\widehat{s}_{n}(x)$ by
\begin{eqnarray}
	\widehat{c}_{n}(x) &\eqdef& \dfrac{c_{2,n+1}(x)}{x} 
	\slabel{eq.def.chat}  \\[2mm]
	\widehat{s}_{n}(x) &\eqdef& \dfrac{s_{2,n+1}(x)}{x}
	\slabel{eq.def.shat}
	\label{eq.def.hat}
\end{eqnarray}
We will show that the polynomials $\widehat{c}_{n}(x)$ and $\widehat{s}_{n}(x)$
have simple and interlacing zeros
in the interval $(-1,0)$.
To prove this,
we will first obtain recurrences for the polynomials $\widehat{c}_{n}(x)$
and $\widehat{s}_{n}(x)$
and then evaluate them at $x=0$ and at $x=-1$.

\begin{lemma}
	\label{lem.polyrec}
\hfill\break\noindent
\vspace*{-7mm}
\begin{itemize}
\item[(a)]
The polynomials $\widehat{c}_{n}(x)$ satisfy the recurrence relation
	\be
		\widehat{c}_{n}(x)
		\;=\;
		\big[ (n+2) x + n+1 \big] \,\widehat{c}_{n-1}(x)
		\,+\,
		x(x+1)\, \widehat{c}_{n-1}'(x)
		\label{eq.polyrec.cycle}
	\ee
\item[(b)]
The polynomials $\widehat{s}_{n}(x)$ satisfy the recurrence relation
        \be
		\widehat{s}_{n}(x)
                \;=\; 
		\big[ (n+2) x + 1 \big] \, \widehat{s}_{n-1}(x)
                \,+\,
		x(x+1)\,\widehat{s}_{n-1}'(x)
		\label{eq.polyrec.subset}
        \ee
\end{itemize}
%
\end{lemma}
\proof
As a consequence of the recurrences
\reff{eq.rec.stirlingcycle.2nd}/\reff{eq.rec.stirlingsubset.2nd},
the polynomials $c_{2,n}(x)$ and $s_{2,n}(x)$ satisfy the recurrences
\begin{eqnarray}
	c_{2,n}(x)
	       &=&
               [n x + (n-1)] \,c_{2,n-1}(x)
               \,+\,
               x(x+1)\, c_{2,n-1}'(x)  \\[2mm]
	s_{2,n}(x)
	&=& 
                n x\, s_{2,n-1}(x)
                \,+\,
                x(x+1)\, s_{2,n-1}'(x)
\end{eqnarray}
from which the recurrences~\reff{eq.polyrec.cycle}/\reff{eq.polyrec.subset}
follow.
\qed

%
%
%


\begin{lemma}
For all $n\geq 1$, the polynomials $\widehat{c}_n(x)$ and $\widehat{s}_n(x)$ 
have the following calues at $x=0$ and at $x=-1$:
	\begin{subeqnarray}
		\widehat{c}_n(0) & = & (n+1)!\\[1mm]
		\widehat{s}_n(0) & = & 1\\[1mm]
		\widehat{c}_n(-1) & = & (-1)^n\\[1mm]
		\widehat{s}_n(-1) & = & (-1)^n\,(n+1)!
	\label{eq.eval}
	\end{subeqnarray}
\label{lem.eval}
\end{lemma}
\proof
These follow by substituting $x=0$ and $x=-1$
in the recurrences~\reff{eq.polyrec.cycle} and~\reff{eq.polyrec.subset}
and then using induction.
\qed

We are now ready to prove Theorem~\ref{thm.stirling.2nd.zeros}.

%
%

\proofof{Theorem~\ref{thm.stirling.2nd.zeros}}
We prove it for the polynomial sequence $\widehat{c}_{2,n}(x)$.
The proof for the polynomial sequence $\widehat{s}_{2,n}(x)$ 
is mostly the same and is left to the reader.

We first rewrite equation~\reff{eq.polyrec.cycle} as
\be
\widehat{c}_n(x) 
\;=\; 
x\:\big((n+2) \widehat{c}_{n-1}(x) +  (x+1)\widehat{c}_{n-1}'(x)\big)
\,+\,
(n+1)\: \widehat{c}_{n-1}(x)
\ee

Now, we use induction on $n$. 
For the case $n=0$, 
the polynomial $p_0(x) = 1$ 
has no roots and hence the theorem is vacuously true. 
The case $n=1$ can be checked to be true as well. 
Let us assume that the zeros of $\widehat{c}_{n-1}(x)$ 
are distinct and are contained in the interval $(-1,0)$. 

Now let $p_n(x) \eqdef
        {\rm GCD}\left(\widehat{c}_{n-1}(x), 
        (x+1) \widehat{c}_{n-1}'(x)\right)$.
By Rolle's theorem,
the zeros of $\widehat{c}_{n-1}'(x)$
interlace those of $\widehat{c}_{n-1}(x)$
and thus, lie in the interval $(-1,0)$ and are simple.
Thus we get $p_n(x) = 1$.
Thus, there are no common zeros between $\widehat{c}_{n-1}(x)$,
$\widehat{c}_{n-1}'(x)$
and $\widehat{c}_{n}(x)$.

Now let $-1<\lambda_{n-1}^{(n-1)}<\ldots<\lambda_1^{(n-1)}<0$
be the zeros of $\widehat{c}_{n-1}(x)$.
We then have $\widehat{c}_{n-1}'(\lambda_j^{(n-1)})>0$
if $j$ is odd,
and
$\widehat{c}_{n-1}'(\lambda_j^{(n-1)})<0$ 
if $j$ is even.
Hence, 
\be
{\rm sgn}\Big( \widehat{c}_n(\lambda_j^{(n-1)})\Big)
\;=\;
(-1)^{j}\;. 
\ee
Furthermore, using Lemma~\ref{lem.eval},
we also have $\widehat{c}_n(-1) = (-1)^n$ 
and $\widehat{c}_n(0)=(n+1)!$. 
Thus, the proof of this theorem follows by using 
the intermediate value theorem.
\qed

Let us finally conclude this section by mentioning
that even though we have included a proof 
of Theorem~\ref{thm.stirling.2nd.zeros}
for the sake of completeness, 
one can obtain this result as an immediate consequence of
a seminal result due to Liu and Wang~\cite{Liu_06}.

\section{Location of zeros for order $\bm{r \ge 3}$}
\label{sec.zero.distribution}

\newcommand{\RC}{\mathcal{RC}}
\newcommand{\RS}{\mathcal{RS}}

In this section we study the zeros of the
$r$th-order Stirling cycle and subset polynomials for $r \ge 3$.

\subsection[Existence of nonreal roots for $r \ge 3$]{Existence of nonreal roots for $\bm{r \ge 3}$}

Conjectures~\ref{conj.cycle}(c) and \ref{conj.subset}(c) assert that,
for $r\geq 3$, the $r$th-order Stirling cycle and subset polynomials
are not negative-real-rooted.
In this subsection we will prove this conjecture
for the cycle case for all $r \ge 3$,
and for the subset case for all $r \ge 4$.
The subset case with $r=3$ remains open.

We already know that $c_{r,n}(0)=s_{r,n}(0)=0$ for \mbox{$n\geq 1$}.
One sufficient condition to show that a polynomial $p(x)$ has 
non-real complex zeros is to show that the polynomial $\dfrac{d^m}{dx^m} p(x)$
for some $m\leq{\rm deg}(p)$  has non-real complex zeros.
In this section, we will show that the polynomials
$\dfrac{d^{n-3}}{dx^{n-3}}(c_{r,n}(x)/x)$
and $\dfrac{d^{n-3}}{dx^{n-3}}(s_{r,n}(x)/x)$,
which are quadratic polynomials,
have a negative discriminant when $r,n\geq 3$ in the cycle case and
$r\geq 4$, $n\geq 3$ in the subset case.
This will prove that in these cases
the polynomials $c_{r,n}(x)$ and $s_{r,n}(x)$ must have some complex
non-real zeros.
We will then provide plots for the distribution of these zeros
which suggest that these zeros accumulate around some
mysterious limiting curves.

\bigskip

We first notice that
\be
\dfrac{d^{n-3}}{dx^{n-3}}(c_{r,n}(x)/x)
\;=\;
\dfrac{(n-1)!}{2} \stirlingcycle{n}{n}^{(r)} x^2
\,+\,
(n-2)! \stirlingcycle{n}{n-1}^{(r)} x
\,+\,
(n-3)! \stirlingcycle{n}{n-2}
\label{eq.complex.quadratic}
\ee
We have already seen in \reff{eq.stircyc.diagonal} that 
$\stirlingcycle{n}{n}^{(r)} = (rn)!/(r^n n!)$.
We now obtain the following exact expressions
either by using the 
combinatorial descriptions or by the recurrences \reff{rec.cycle}:
\begin{subeqnarray}
\stirlingcycle{n}{n-1}^{(r)}
        &=&
\dfrac{(r(n-1)+1)!}{(r+1) r^{n-2} (n-2)!}  \\[2mm]
\stirlingcycle{n}{n-2}^{(r)}
        &=&
\dfrac{(r(n-2)+2)!\;\bigl(r(r+2)(n-1)\,+\,2\bigr)}{2\,(r+1)^2 \, (r+2) \, r^{n-3}\, (n-3)! }
\end{subeqnarray}

With this information, we prove the following proposition:

\begin{proposition}
Let $D_{r}(n)$ be the discriminant of the 
quadratic polynomial\\
$\dfrac{d^{n-3}}{dx^{n-3}}(c_{r,n}(x)/x)$.
For all $r,n\geq 3$,
we have $D_r(n)<0$.
\label{prop.cycle.discriminant}
\end{proposition}

\proof
From~\reff{eq.complex.quadratic}, we write out the exact expression for $D_r(n)$,
which is
\begin{subeqnarray}
D_{n}(r) 
	&=&
\left((n-2)!  \stirlingcycle{n}{n-1}^{(r)}\right)^2
\,-\,
	2(n-1)!(n-3)! \stirlingcycle{n}{n}^{(r)} \stirlingcycle{n}{n-2}^{(r)}\\
	&=&
	\dfrac{ (r(n-1)+1)!^2}{(r+1)^2 r^{2n-4}}
	\nonumber\\
	&&
	\,-\,
	\dfrac{  (rn)! (r(n-2)+2)! ((r+2)r(n-1)+2) }{ n (r+2)(r+1)^2 r^{2n-3}}
	\\[2mm]
	&=&
	\dfrac{(r(n-1)+1)!(r(n-2)+2)!}{(r+1)^2 r^{2n-4}}\cdot
	\left( \dfrac{(r(n-1)+1)!}{(r(n-2)+2)!}\right.
	\nonumber\\
	&&
	\,-\,
	\left.\dfrac{ ((r+2)r(n-1)+2)}{n (r+2)r}\cdot \dfrac{(rn)!}{(r(n-1)+1)!}
	\right)\\[2mm]
	&=&
	\dfrac{(r(n-1)+1)!(r(n-2)+2)!}{(r+1)^2 r^{2n-4}}\cdot\nonumber\\
	&&
	 \Biggl( 
	 (rn - (r-1))(rn-(r-1)-1)\cdots (rn - 2(r-1)+1)
	\nonumber\\
	&&
	\,-\,
	\dfrac{ ((r+2)r(n-1)+2)}{n (r+2)r}\cdot
	(rn)(rn-1)\cdots (rn - (r-1)+1)
	\Biggr)
\end{subeqnarray}
Notice that when $r,n\geq 3$,
we have 
\begin{subeqnarray}
(rn - (r-1))(rn - (r-1)-1) - \dfrac{ ((r+2)r(n-1)+2)}{n (r+2)r}\cdot (rn)(rn-1)
\nonumber\\
=\;\dfrac{-n(r^3-2r)+(r^3-4r+2)}{r+2}<0\;.
\end{subeqnarray}
It is easy to compare the remaining terms 
and it follows that $D_r(n)<0$ for $r,n\geq 3$.
\qed

\begin{corollary}
   \label{cor.cycle.nonreal}
For $r,n \ge 3$, the polynomial $c_{r,n}(x)$ has some nonreal zeros.
\end{corollary}

One can use a  similar argument to show that the discriminant of the
quadratic polynomial $\dfrac{d^{n-3}}{dx^{n-3}}(s_{r,n}(x)/x)$ is negative
when $r\geq 4$ and $n\geq 3$;
we leave the details to the reader.
This implies:

\begin{corollary}
   \label{cor.subset.nonreal}
For $r \ge 4$ and $n \ge 3$, the polynomial $s_{r,n}(x)$ has some nonreal zeros.
\end{corollary}

For $r=3$, computational evidence suggests that the zeros of
$\dfrac{d^{n-3}}{dx^{n-3}}(s_{3,n}(x)/x)$ are real,
but that the polynomials $s_{3,n}(x)$ nevertheless have
some non-real complex zeros for $n \ge 4$.

\bigskip

Let us first establish some notation before we provide the plots 
for the zeros of the polynomials
$c_{r,n}(x)$ and $s_{r,n}(x)$.
Let $\scra \subset \mathbb{C}$ be a subset of complex numbers,
let $\alpha\in \mathbb{R}$ be a real number. 
We then define the set $\alpha \scra$ as follows:
\be
\alpha \, \scra \;\eqdef\; \{\alpha x \: : x\in \scra \}\;.
\ee
For a fixed $r\geq 1$, we  let $\RC_{r,n}$ and $\RS_{r,n}$ 
denote the set of zeros of the
$r$th-order Stirling cycle polynomial $c_{r,n}(x)$,
and
of the $r$th-order Stirling subset polynomial $s_{r,n}(x)$,
respectively:
\begin{subeqnarray}
	\RC_{r,n} &\eqdef& \{z\in \mathbb{C} \:\colon c_{r,n}(z) = 0\}\\
	\RS_{r,n} &\eqdef& \{z\in \mathbb{C} \:\colon s_{r,n}(z) = 0\}\;.
\end{subeqnarray}
In Sections~\ref{sec.cycle.plots} and \ref{sec.subset.plots},
we provide normalised plots of $\RC_{r,n}$ and $\RS_{r,n}$, respectively,
for $3 \le r \le 10$.
These plots will immediately suggest some problems 
and conjectures which we mention.

\subsection[Stirling cycle polynomials of order $r \ge 3$]{Stirling cycle polynomials of order $\bm{r \ge 3}$}
\label{sec.cycle.plots}

\begin{figure}[t!]
	\subfloat[\centering $r=3$-rd order Stirling cycle ]{
	\scalebox{0.25}{
		\includegraphics{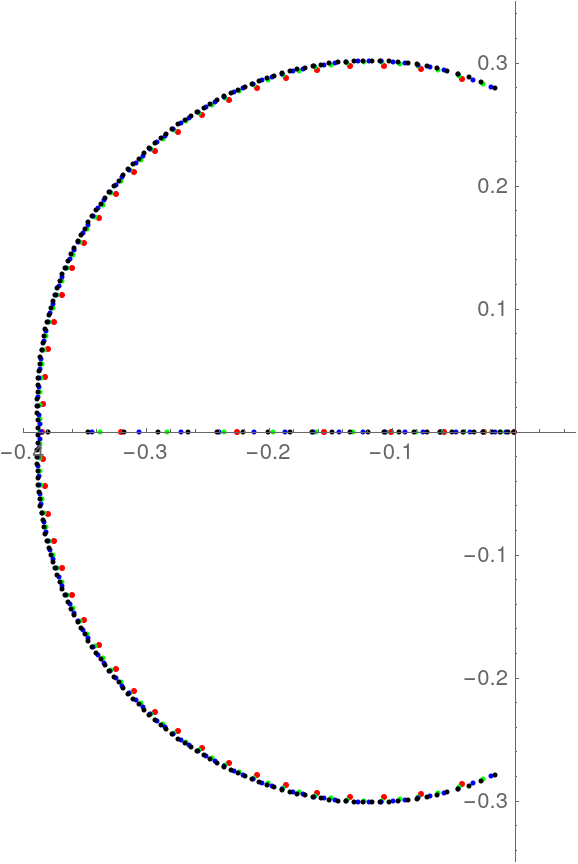}}}
	\qquad
        \subfloat[\centering $r=4$-th order Stirling cycle ]{
		 \scalebox{0.25}{
		 \includegraphics{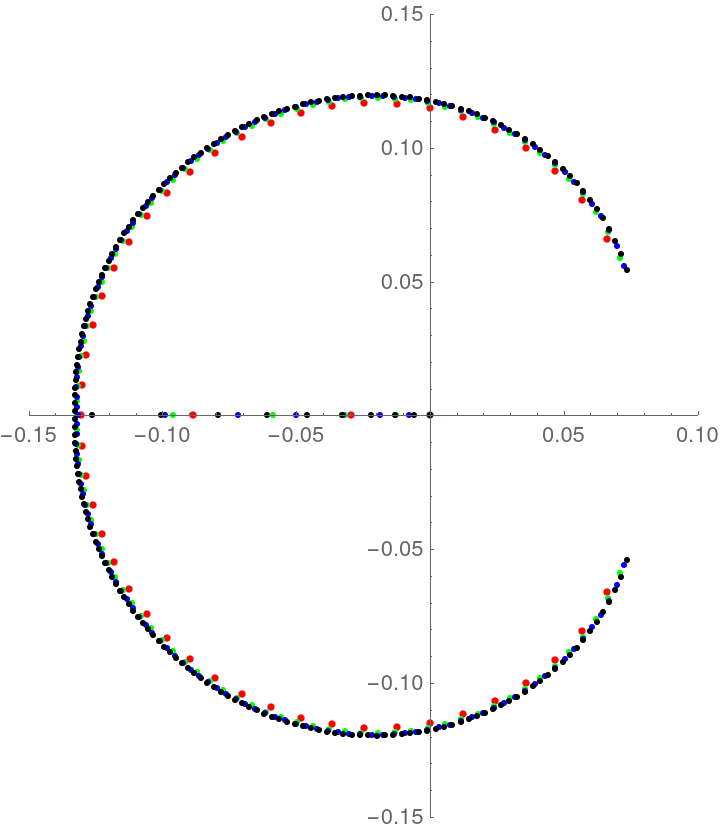}}}\\
        \subfloat[\centering $r=5$-th order Stirling cycle]{
 	\scalebox{0.25}{
		\includegraphics{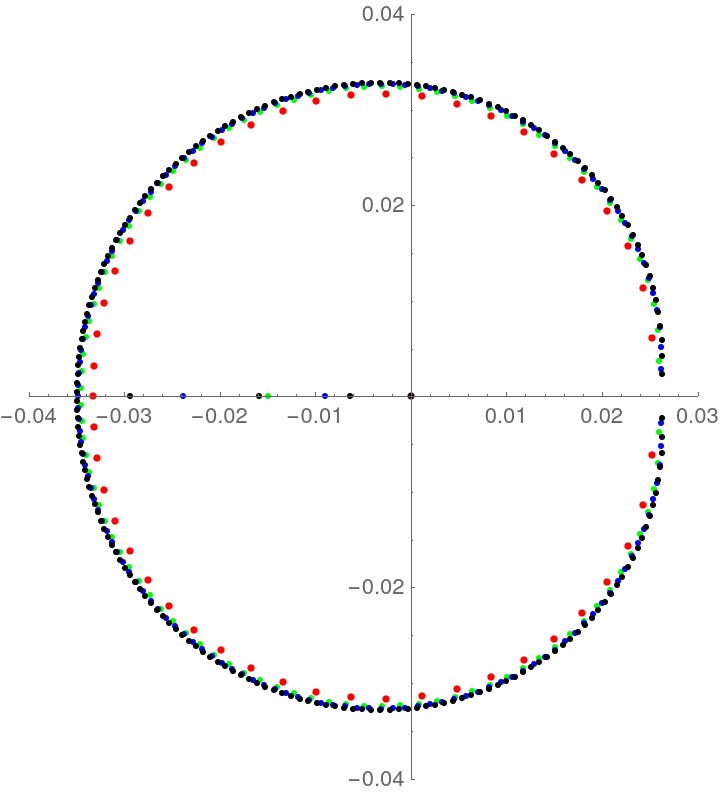}}}
	\qquad
	\subfloat[\centering $r=6$-th order Stirling cycle ]{
        \scalebox{0.25}{
		\includegraphics{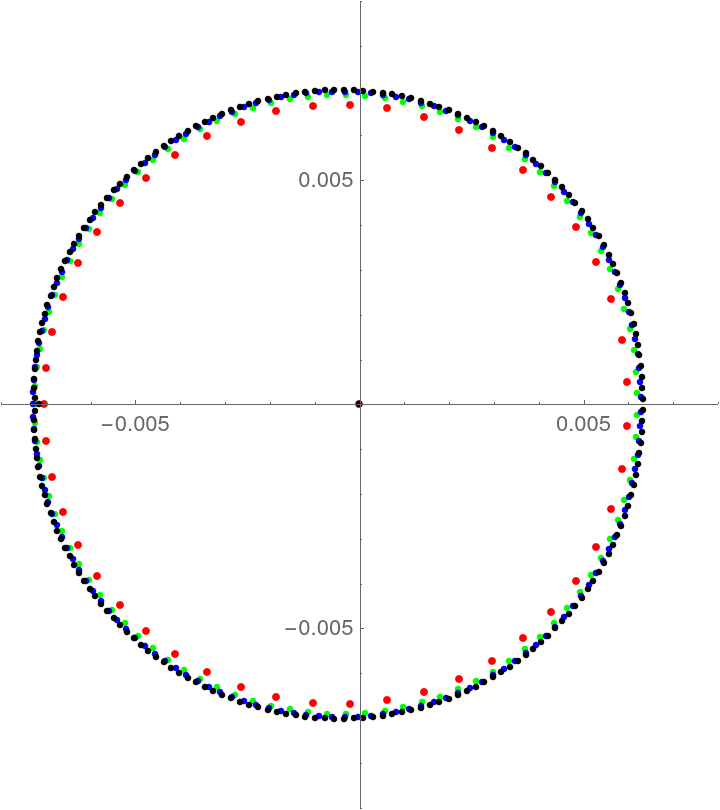}}}
\caption{
We plot $n^{r-2}\RC_{r,n}$ with $r=3,4,5,6$ in subfigures (a), (b), (c)
and (d), respectively,
and with $n=50,100,150,200$ using the different colors
red, green, blue and black, respectively.
}
\label{fig.cycle.a}
\end{figure}

\begin{figure}[t!]
        \subfloat[\centering $r=7$-th order Stirling cycle ]{
        \scalebox{0.25}{
                \includegraphics{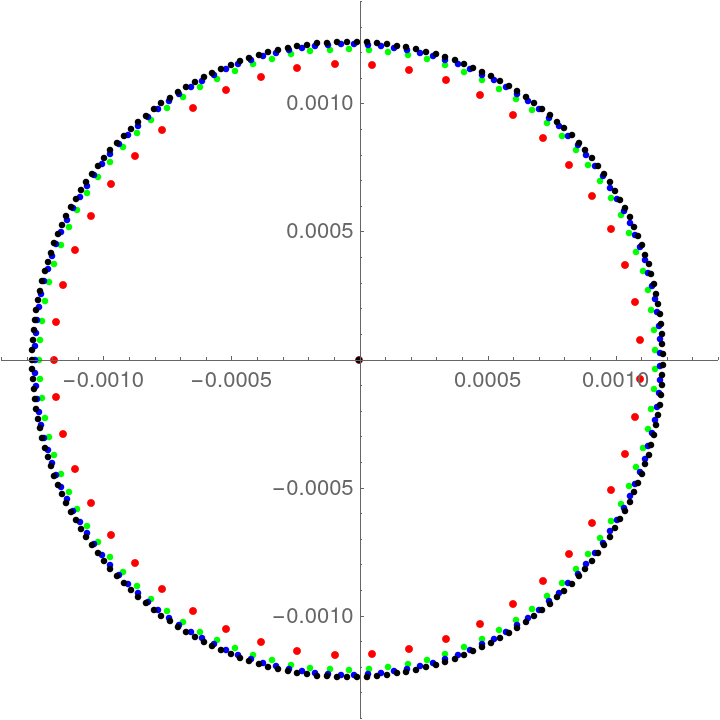}}}
        \qquad
        \subfloat[\centering $r=8$-th order Stirling cycle ]{
                 \scalebox{0.25}{
                 \includegraphics{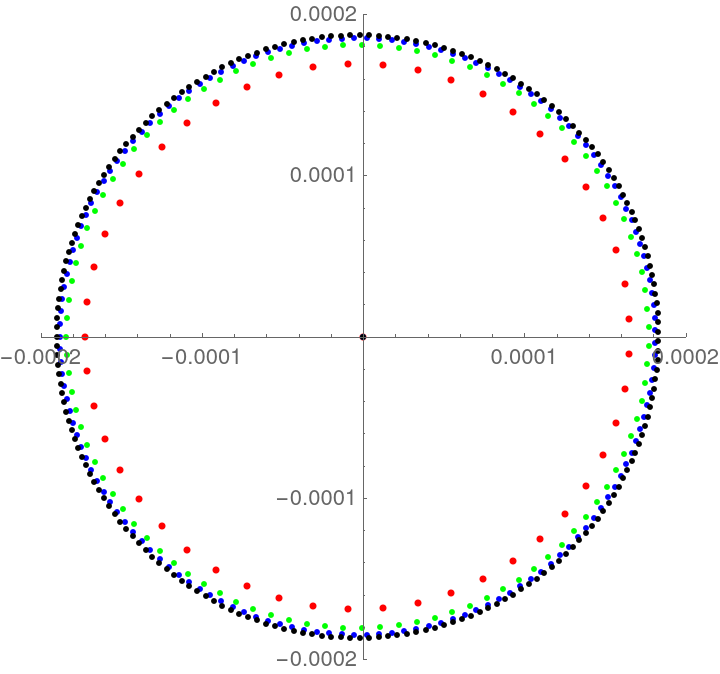}}}\\
        \subfloat[\centering $r=9$-th order Stirling cycle]{
        \scalebox{0.25}{
                \includegraphics{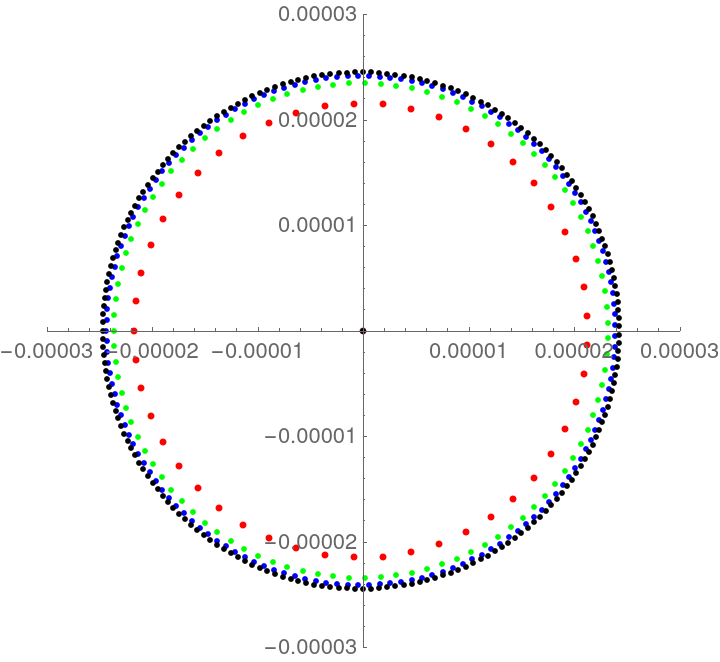}}}
        \qquad
        \subfloat[\centering $r=10$-th order Stirling cycle ]{
        \scalebox{0.25}{
                \includegraphics{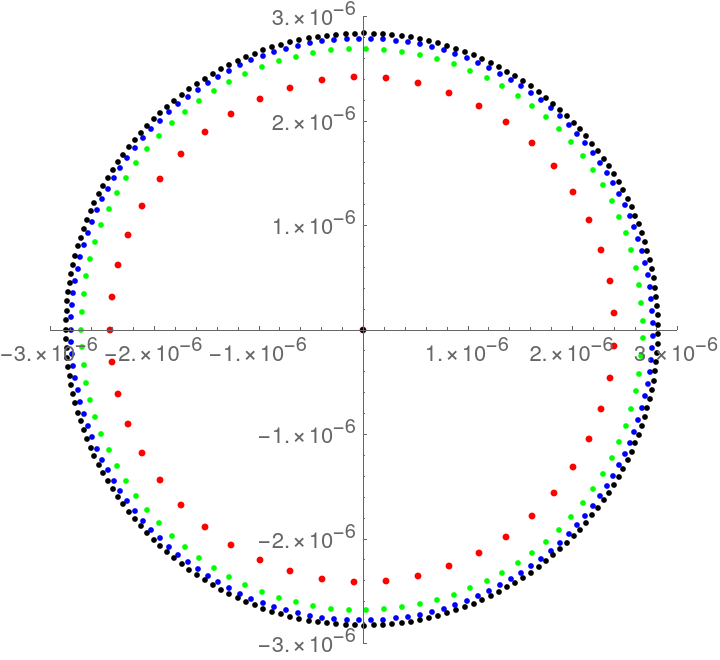}}}
\caption{
We plot $n^{r-2}\RC_{r,n}$ with $r=7,8,9,10$ in subfigures (a), (b), (c)
and (d), respectively,
and with $n=50,100,150,200$ using the different colors
red, green, blue and black, respectively.
}
	\label{fig.cycle.b}
\end{figure}

We fix $r$ and plot the set $n^{r-2}\RC_{r,n}$ 
for $n=50,100,150,200$ 
using the different colors red, green, blue and black, respectively.
The overlapping plots suggest that the zeros have some asymptotic 
distribution supported on some limiting curves.
We plot the cases $r=3,4,5,6$ in Figure~\ref{fig.cycle.a} and
$r=7,8,9,10$ in Figure~\ref{fig.cycle.b}.

Our plots suggest the following problem:

\begin{problem} For a fixed $r\geq 1$, study the asymptotic distribution of
	the set of normalised zeros $n^{r-2}\,\RC_{r,n}$.
	Do they accumulate on a curve?  If so, what curve?
\end{problem}

In the case $r=1$, we clearly get that 
the normalised zeros accumulate uniformly 
on the unit interval $[-1,0]$.

When $r=3$, we conjecture the following:
\begin{conjecture} The zeros of the third-order Stirling cycle polynomials
        $c_{3,n}(x)$ lie on the left-half plane.
\end{conjecture}

\clearpage

\subsection[Stirling subset polynomials of order $r \ge 3$]{Stirling subset polynomials of order $\bm{r \ge 3}$}
\label{sec.subset.plots}

\begin{figure}[t!]
        \subfloat[\centering $r=3$-rd order Stirling subset ]{
        \scalebox{0.25}{
                \includegraphics{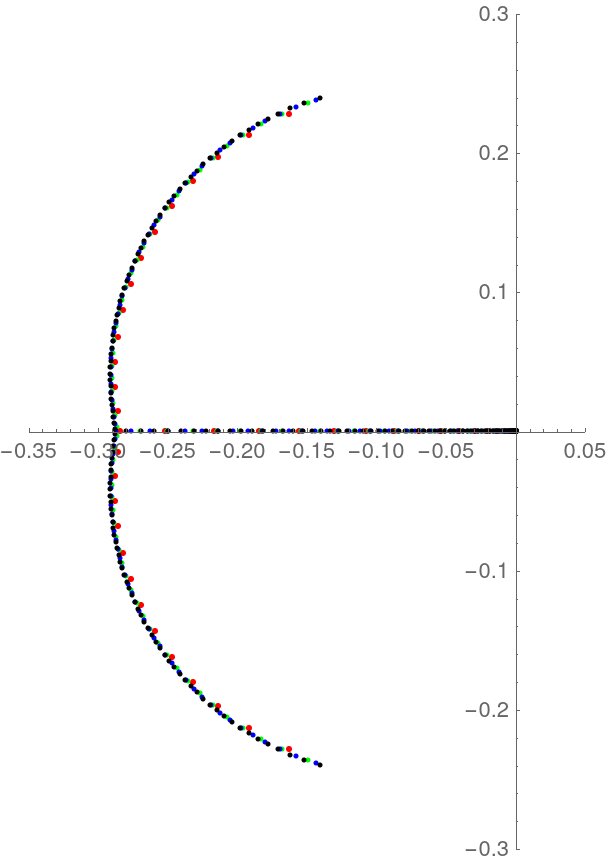}}}
        \qquad
        \subfloat[\centering $r=4$-th order Stirling subset ]{
                 \scalebox{0.25}{
                 \includegraphics{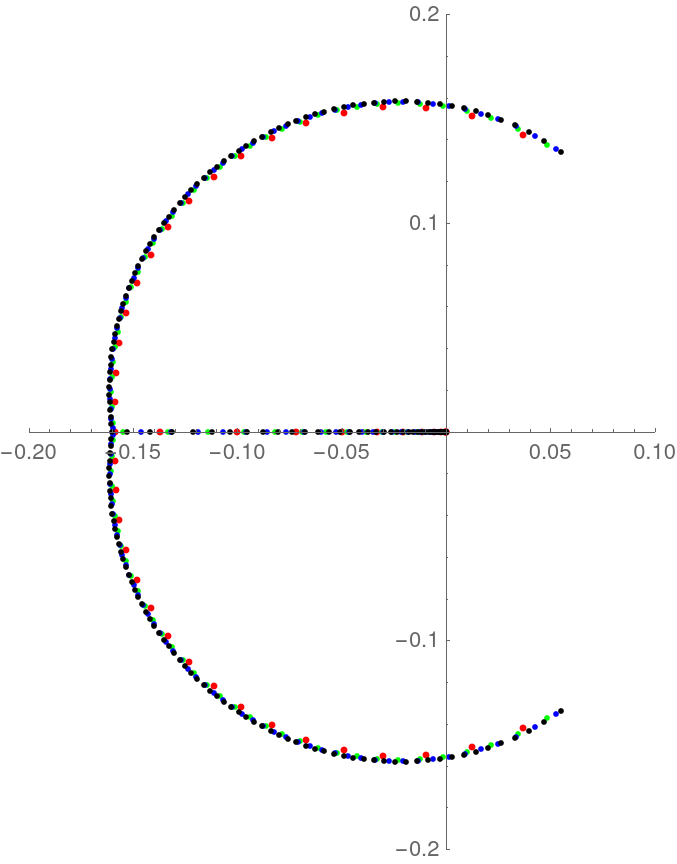}}}\\
        \subfloat[\centering $r=5$-th order Stirling subset]{
        \scalebox{0.25}{
                \includegraphics{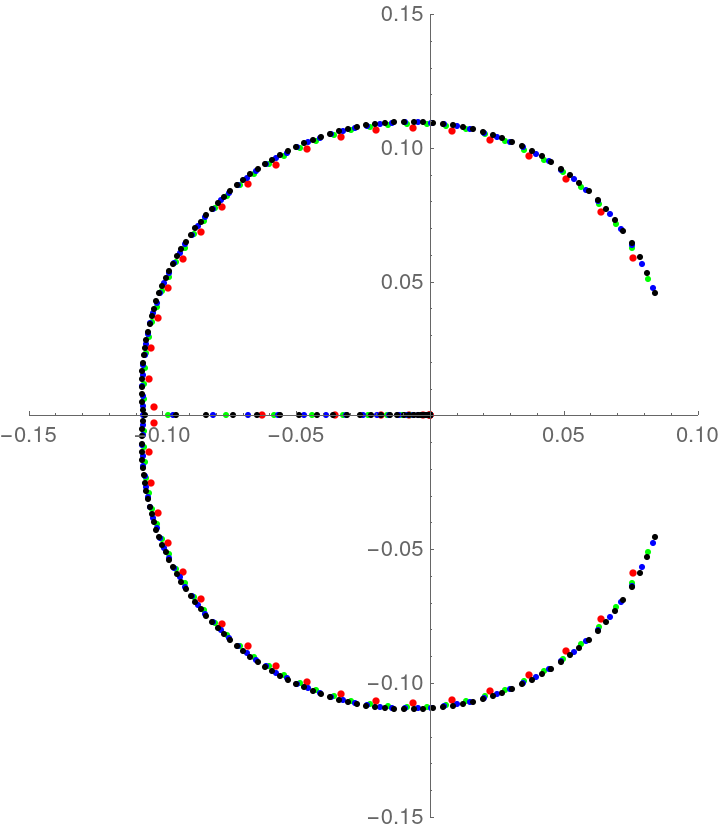}}}
        \qquad
        \subfloat[\centering $r=6$-th order Stirling subset ]{
        \scalebox{0.25}{
                \includegraphics{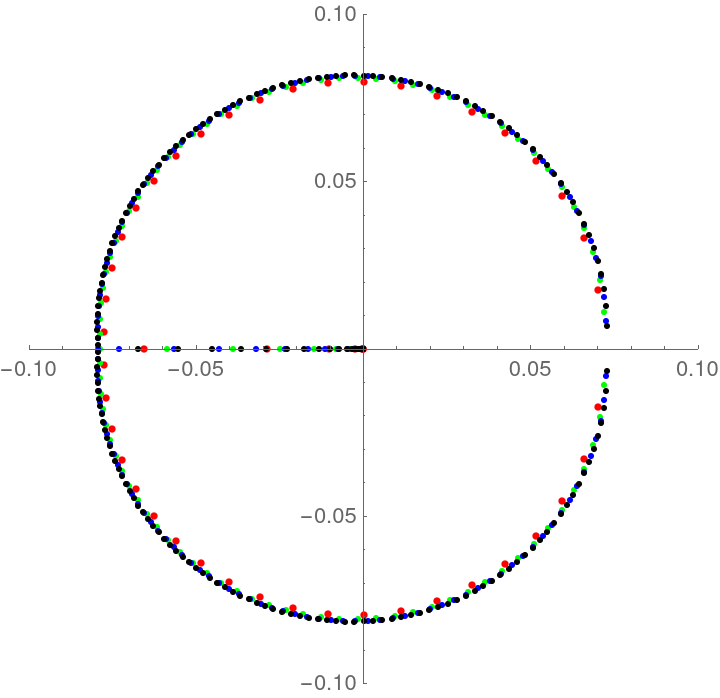}}}
\caption{
We plot $n^{r-2}\RS_{r,n}$ with $r=3,4,5,6$ in subfigures (a), (b), (c)
and (d), respectively,
and with $n=50,100,150,200$ using the different colors
red, green, blue and black, respectively.
}
	\label{fig.subset.a}
\end{figure}

\begin{figure}[t!]
        \subfloat[\centering $r=7$-th order Stirling subset ]{
        \scalebox{0.25}{
                \includegraphics{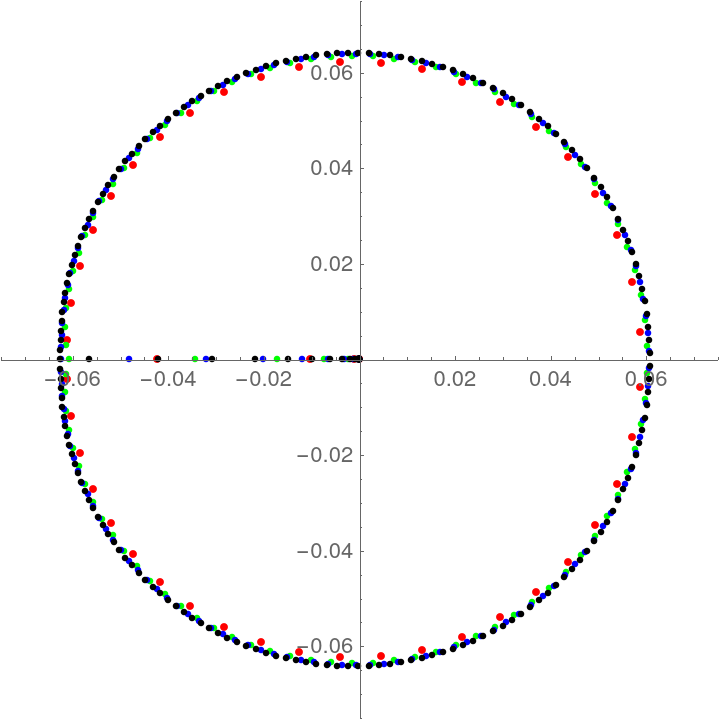}}}
        \qquad
        \subfloat[\centering $r=8$-th order Stirling subset ]{
                 \scalebox{0.25}{
                 \includegraphics{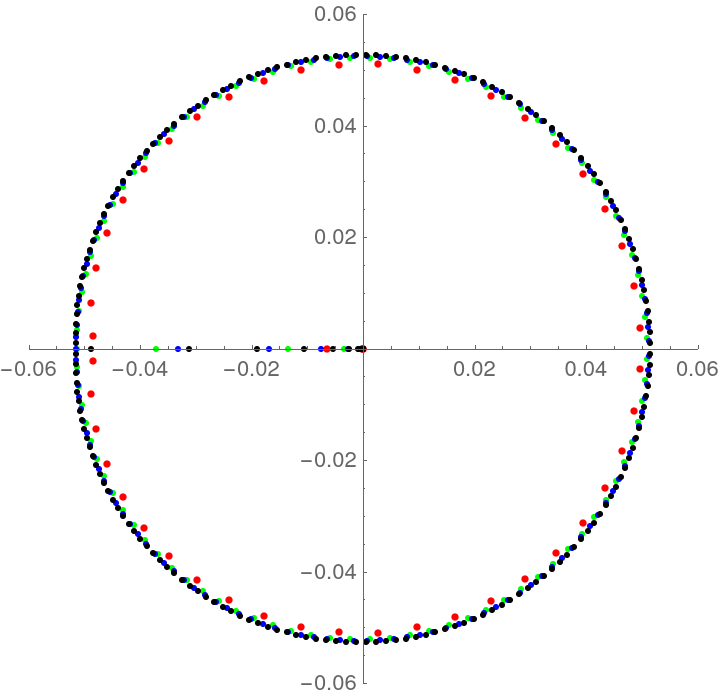}}}\\
        \subfloat[\centering $r=9$-th order Stirling subset]{
        \scalebox{0.25}{
                \includegraphics{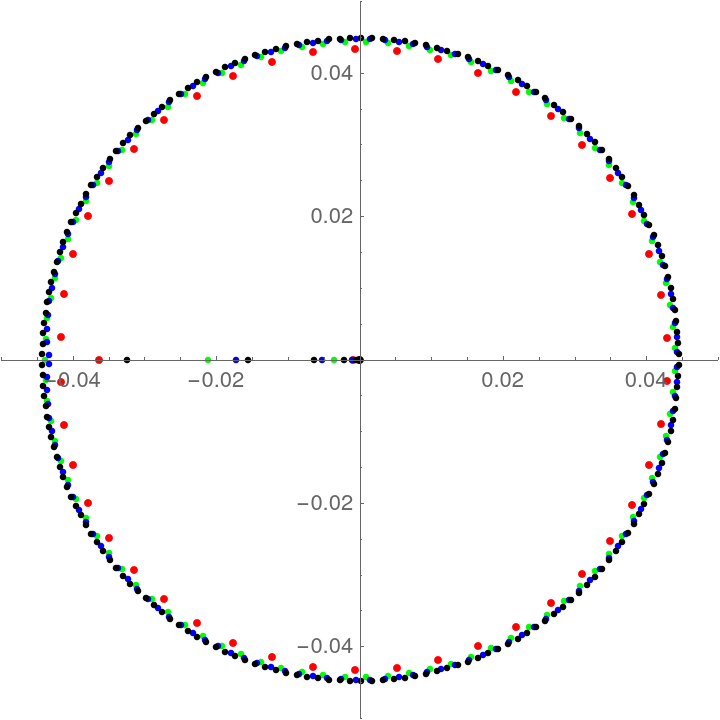}}}
        \qquad
        \subfloat[\centering $r=10$-th order Stirling subset ]{
        \scalebox{0.25}{
                \includegraphics{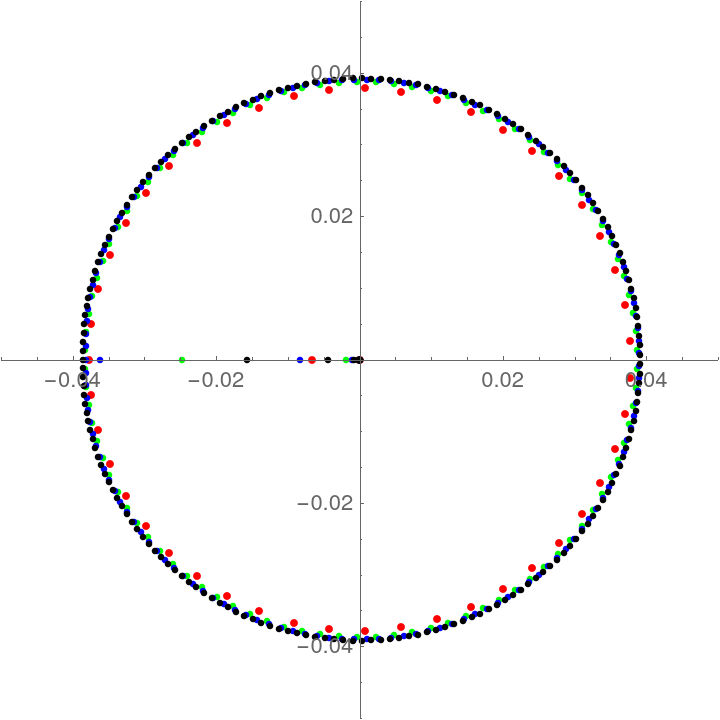}}}
\caption{
We plot $n^{r-2}\RS_{r,n}$ with $r=7,8,9,10$ in subfigures (a), (b), (c)
and (d), respectively,
and with $n=50,100,150,200$ using the different colors
red, green, blue and black, respectively.
}
	\label{fig.subset.b}
\end{figure}

We fix $r$ and plot the set $n^{r-2}\RS_{r,n}$
for $n=50,100,150,200$
using the different colors red, green, blue and black, respectively.
The overlapping plots suggest that the zeros have some asymptotic 
distribution supported on some limiting curves.
We plot the cases $r=3,4,5,6$ in Figure~\ref{fig.subset.a} and
$r=7,8,9,10$ in Figure~\ref{fig.subset.b}.

%

Our plots suggest the following problem:

\begin{problem} For a fixed $r\geq 1$, study the asymptotic distribution of
        the set of normalised zeros $n^{r-2}\,\RS_{r,n}$.
	Do they accumulate on a curve?  If so, what curve?
\end{problem}

In the case $r=1$, 
the asymptotic distribution of the normalised zeros of the Bell polynomials
--- which lie on the negative real axis ---
was studied by Elbert \cite{Elbert_2001a, Elbert_2001b}.

When $r=3$, we conjecture the following:
\begin{conjecture} The zeros of the third-order Stirling subset polynomials
        $s_{3,n}(x)$ lie on the left-half plane.
\end{conjecture}


\subsection*{Acknowledgements}
We thank David Callan, Alex Dyachenko and Markus Kuba for helpful correspondence.
The content in 
Section~\ref{sec.interpretation.stirperm} is due to David Callan.
Alex Dyachenko carried out the numerical tests of Conjecture~\ref{conj.secondeuler}.

The first author was supported
by DIMERS project ANR-18-CE40-0033 funded by Agence Nationale de la Recherche
(ANR, France).
He is currently supported by the
Tsinghua University Shuimu scholarship.

\appendix
\section{Tables of higher-order Stirling cycle and subset numbers}
\label{app.sec.tables}

In this Appendix we provide tables of the $r$th-order
Stirling cycle and subset numbers for $1 \le r \le 4$
and the first few values of $n$.

\subsection{Higher-order Stirling cycle numbers}

\bigskip

%
%
\begin{table}[H]
\centering
\small
\begin{tabular}{c|rrrrrrrrr|r|}
\cline{2-10}
   & \multicolumn{9}{|c|}{$\bm{r=1}$ \bf (ordinary Stirling cycle numbers)} &  \multicolumn{1}{c}{\quad} \\
\hline
$n \setminus k$ & 0 & 1 & 2 & 3 & 4 & 5 & 6 & 7 & 8 & Row sums \\
\hline
0 & 1 &  &  &  &  &  &  &  &  & 1  \\
1 & 0 & 1 &  &  &  &  &  &  &  & 1  \\
2 & 0 & 1 & 1 &  &  &  &  &  &  & 2  \\
3 & 0 & 2 & 3 & 1 &  &  &  &  &  & 6  \\
4 & 0 & 6 & 11 & 6 & 1 &  &  &  &  & 24  \\
5 & 0 & 24 & 50 & 35 & 10 & 1 &  &  &  & 120  \\
6 & 0 & 120 & 274 & 225 & 85 & 15 & 1 &  &  & 720  \\
7 & 0 & 720 & 1764 & 1624 & 735 & 175 & 21 & 1 &  & 5040  \\
8 & 0 & 5040 & 13068 & 13132 & 6769 & 1960 & 322 & 28 & 1 & 40320  \\
\hline
\end{tabular}
\end{table}

\bigskip
\bigskip

%
%
\begin{table}[H]
\centering
\scriptsize
\begin{tabular}{c|rrrrrrrrr|r|}
\cline{2-10}
   & \multicolumn{9}{|c|}{\small $\bm{r=2}$} &  \multicolumn{1}{c}{\quad} \\
\hline
$n \setminus k$ & 0 & 1 & 2 & 3 & 4 & 5 & 6 & 7 & 8 & Row sums \\
\hline
0 & 1 &  &  &  &  &  &  &  &  & 1  \\
1 & 0 & 1 &  &  &  &  &  &  &  & 1  \\
2 & 0 & 2 & 3 &  &  &  &  &  &  & 5  \\
3 & 0 & 6 & 20 & 15 &  &  &  &  &  & 41  \\
4 & 0 & 24 & 130 & 210 & 105 &  &  &  &  & 469  \\
5 & 0 & 120 & 924 & 2380 & 2520 & 945 &  &  &  & 6889  \\
6 & 0 & 720 & 7308 & 26432 & 44100 & 34650 & 10395 &  &  & 123605  \\
7 & 0 & 5040 & 64224 & 303660 & 705320 & 866250 & 540540 & 135135 &  & 2620169  \\
8 & 0 & 40320 & 623376 & 3678840 & 11098780 & 18858840 & 18288270 & 9459450 & 2027025 & 64074901  \\
\hline
\end{tabular}
\end{table}

\bigskip
\bigskip

%
%
\begin{table}[H]
\hspace*{-1.5cm}
\scriptsize
\begin{tabular}{c|rrrrrrrr|r|}
\cline{2-9}
   & \multicolumn{8}{|c|}{\small $\bm{r=3}$} &  \multicolumn{1}{c}{\quad} \\
\hline
$n \setminus k$ & 0 & 1 & 2 & 3 & 4 & 5 & 6 & 7 & Row sums \\
\hline
0 & 1 &  &  &  &  &  &  &  & 1  \\
1 & 0 & 2 &  &  &  &  &  &  & 2  \\
2 & 0 & 6 & 40 &  &  &  &  &  & 46  \\
3 & 0 & 24 & 420 & 2240 &  &  &  &  & 2684  \\
4 & 0 & 120 & 3948 & 50400 & 246400 &  &  &  & 300868  \\
5 & 0 & 720 & 38304 & 859320 & 9609600 & 44844800 &  &  & 55352744  \\
6 & 0 & 5040 & 396576 & 13665960 & 258978720 & 2690688000 & 12197785600 &  & 15161519896  \\
7 & 0 & 40320 & 4419360 & 216339552 & 6112906800 & 105205900800 & 1042910668800 & 4635158528000 & 5789608803632  \\
\hline
\end{tabular}
\end{table}

\bigskip
\bigskip

%
%
\begin{table}[H]
\hspace*{-1.8cm}
\scriptsize
\begin{tabular}{c|rrrrrrr|r|}
\cline{2-8}
   & \multicolumn{7}{|c|}{\small $\bm{r=4}$} &  \multicolumn{1}{c}{\quad} \\
\hline
$n \setminus k$ & 0 & 1 & 2 & 3 & 4 & 5 & 6 & Row sums \\
\hline
0 & 1 &  &  &  &  &  &  & 1  \\
1 & 0 & 6 &  &  &  &  &  & 6  \\
2 & 0 & 24 & 1260 &  &  &  &  & 1284  \\
3 & 0 & 120 & 18144 & 1247400 &  &  &  & 1265664  \\
4 & 0 & 720 & 223776 & 38918880 & 3405402000 &  &  & 3444545376  \\
5 & 0 & 5040 & 2756160 & 889945056 & 185253868800 & 19799007228000 &  & 19985153803056  \\
6 & 0 & 40320 & 35307360 & 18478684224 & 6780291598080 & 1663116607152000 & 210384250804728000 & 212054166217509984  \\
\hline
\end{tabular}
\end{table}

\bigskip
\bigskip

\subsection{Higher-order Stirling subset numbers}

\bigskip

%
%
\begin{table}[H]
\centering
\small
\begin{tabular}{c|rrrrrrrrr|r|}
\cline{2-10}
   & \multicolumn{9}{|c|}{$\bm{r=1}$ \bf (ordinary Stirling subset numbers)} &  \multicolumn{1}{c}{\quad} \\
\hline
$n \setminus k$ & 0 & 1 & 2 & 3 & 4 & 5 & 6 & 7 & 8 & Row sums \\
\hline
0 & 1 &  &  &  &  &  &  &  &  & 1  \\
1 & 0 & 1 &  &  &  &  &  &  &  & 1  \\
2 & 0 & 1 & 1 &  &  &  &  &  &  & 2  \\
3 & 0 & 1 & 3 & 1 &  &  &  &  &  & 5  \\
4 & 0 & 1 & 7 & 6 & 1 &  &  &  &  & 15  \\
5 & 0 & 1 & 15 & 25 & 10 & 1 &  &  &  & 52  \\
6 & 0 & 1 & 31 & 90 & 65 & 15 & 1 &  &  & 203  \\
7 & 0 & 1 & 63 & 301 & 350 & 140 & 21 & 1 &  & 877  \\
8 & 0 & 1 & 127 & 966 & 1701 & 1050 & 266 & 28 & 1 & 4140  \\
\hline
\end{tabular}
\end{table}

\bigskip
\bigskip

%
%
\begin{table}[H]
\centering
\small
\begin{tabular}{c|rrrrrrrrr|r|}
\cline{2-10}
   & \multicolumn{9}{|c|}{\small $\bm{r=2}$ \bf (Ward numbers)} &  \multicolumn{1}{c}{\quad} \\
\hline
$n \setminus k$ & 0 & 1 & 2 & 3 & 4 & 5 & 6 & 7 & 8 & Row sums \\
\hline
0 & 1 &  &  &  &  &  &  &  &  & 1  \\
1 & 0 & 1 &  &  &  &  &  &  &  & 1  \\
2 & 0 & 1 & 3 &  &  &  &  &  &  & 4  \\
3 & 0 & 1 & 10 & 15 &  &  &  &  &  & 26  \\
4 & 0 & 1 & 25 & 105 & 105 &  &  &  &  & 236  \\
5 & 0 & 1 & 56 & 490 & 1260 & 945 &  &  &  & 2752  \\
6 & 0 & 1 & 119 & 1918 & 9450 & 17325 & 10395 &  &  & 39208  \\
7 & 0 & 1 & 246 & 6825 & 56980 & 190575 & 270270 & 135135 &  & 660032  \\
8 & 0 & 1 & 501 & 22935 & 302995 & 1636635 & 4099095 & 4729725 & 2027025 & 12818912  \\
\hline
\end{tabular}
\end{table}

\bigskip
\bigskip


%
%
\begin{table}[H]
\hspace*{-1.5cm}
\scriptsize
\begin{tabular}{c|rrrrrrrrr|r|}
\cline{2-10}
   & \multicolumn{9}{|c|}{\small $\bm{r=3}$} &  \multicolumn{1}{c}{\quad} \\
\hline
$n \setminus k$ & 0 & 1 & 2 & 3 & 4 & 5 & 6 & 7 & 8 & Row sums \\
\hline
0 & 1 &  &  &  &  &  &  &  &  & 1  \\
1 & 0 & 1 &  &  &  &  &  &  &  & 1  \\
2 & 0 & 1 & 10 &  &  &  &  &  &  & 11  \\
3 & 0 & 1 & 35 & 280 &  &  &  &  &  & 316  \\
4 & 0 & 1 & 91 & 2100 & 15400 &  &  &  &  & 17592  \\
5 & 0 & 1 & 210 & 10395 & 200200 & 1401400 &  &  &  & 1612206  \\
6 & 0 & 1 & 456 & 42735 & 1611610 & 28028000 & 190590400 &  &  & 220273202  \\
7 & 0 & 1 & 957 & 158301 & 10335325 & 333533200 & 5431826400 & 36212176000 &  & 41988030184  \\
8 & 0 & 1 & 1969 & 549549 & 57962905 & 3073270200 & 89625135600 & 1394168776000 & 9161680528000 & 10648606224224  \\
\hline
\end{tabular}
\end{table}

\bigskip
\bigskip


%
%
\begin{table}[H]
\hspace*{-1.5cm}
\scriptsize
\begin{tabular}{c|rrrrrrrr|r|}
\cline{2-9}
   & \multicolumn{8}{|c|}{\small $\bm{r=4}$} &  \multicolumn{1}{c}{\quad} \\
\hline
$n \setminus k$ & 0 & 1 & 2 & 3 & 4 & 5 & 6 & 7 & Row sums \\
\hline
0 & 1 &  &  &  &  &  &  &  & 1  \\
1 & 0 & 1 &  &  &  &  &  &  & 1  \\
2 & 0 & 1 & 35 &  &  &  &  &  & 36  \\
3 & 0 & 1 & 126 & 5775 &  &  &  &  & 5902  \\
4 & 0 & 1 & 336 & 45045 & 2627625 &  &  &  & 2673007  \\
5 & 0 & 1 & 792 & 231231 & 35735700 & 2546168625 &  &  & 2582136349  \\
6 & 0 & 1 & 1749 & 981981 & 300179880 & 53469541125 & 4509264634875 &  & 4563035339611  \\
7 & 0 & 1 & 3718 & 3741738 & 2002016016 & 666586946025 & 135277939046250 & 13189599057009375 & 13325545588763123  \\
\hline
\end{tabular}
\end{table}

\bigskip
\bigskip

\section{Tables of higher-order quasi-Eulerian cycle and subset numbers}
\label{app.sec.tables.qe}

In this Appendix we provide tables of the $r$th-order
quasi-Eulerian cycle and subset numbers for $1 \le r \le 5$
and the first few values of $n$.

\subsection{Higher-order quasi-Eulerian cycle numbers}
   \label{app.eq.cycle}

%
%
\begin{table}[H]
\centering
\small
\begin{tabular}{c|rrrrrrrrr|r|}
\cline{2-10}
   & \multicolumn{9}{|c|}{\small $\bm{r=1}$} &  \multicolumn{1}{c}{\quad} \\
\hline
$n \setminus k$ & 0 & 1 & 2 & 3 & 4 & 5 & 6 & 7 & 8 & Row sums \\
\hline
0 & 1 &  &  &  &  &  &  &  &  & 1  \\
1 & 1 & 0 &  &  &  &  &  &  &  & 1  \\
2 & 0 & 1 & 0 &  &  &  &  &  &  & 1  \\
3 & 0 & $-1$ & 2 & 0 &  &  &  &  &  & 1  \\
4 & 0 & 2 & $-7$ & 6 & 0 &  &  &  &  & 1  \\
5 & 0 & $-6$ & 29 & $-46$ & 24 & 0 &  &  &  & 1  \\
6 & 0 & 24 & $-146$ & 329 & $-326$ & 120 & 0 &  &  & 1  \\
7 & 0 & $-120$ & 874 & $-2521$ & 3604 & $-2556$ & 720 & 0 &  & 1  \\
8 & 0 & 720 & $-6084$ & 21244 & $-39271$ & 40564 & $-22212$ & 5040 & 0 & 1  \\
\hline
\end{tabular}
\end{table}

\bigskip
\bigskip

\begin{table}[H]
\centering
\small
\begin{tabular}{c|rrrrrrrrr|r|}
\cline{2-10}
   & \multicolumn{9}{|c|}{\small $\bm{r=2}$ \bf (second-order Eulerian numbers) } &  \multicolumn{1}{c}{\quad} \\
\hline
$n \setminus k$ & 0 & 1 & 2 & 3 & 4 & 5 & 6 & 7 & 8 & Row sums \\
\hline
 0 & 1 &  &  &  &  &  &  &  &  & 1 \\
 1 & 1 &  &  &  &  &  &  &  &  & 1 \\
 2 & 1 & 2 &  &  &  &  &  &  &  & 3 \\
 3 & 1 & 8 & 6 &  &  &  &  &  &  & 15 \\
 4 & 1 & 22 & 58 & 24 &  &  &  &  &  & 105 \\
 5 & 1 & 52 & 328 & 444 & 120 &  &  &  &  & 945 \\
 6 & 1 & 114 & 1452 & 4400 & 3708 & 720 &  &  &  & 10395 \\
 7 & 1 & 240 & 5610 & 32120 & 58140 & 33984 & 5040 &  &  & 135135 \\
 8 & 1 & 494 & 19950 & 195800 & 644020 & 785304 & 341136 & 40320 &  & 2027025 \\
\hline
\end{tabular}
\end{table}


\bigskip
\bigskip

\begin{table}[H]
\hspace*{-1.5cm}
\scriptsize
\begin{tabular}{c|rrrrrrrr|r|}
\cline{2-9}
   & \multicolumn{8}{|c|}{\small $\bm{r=3}$} &  \multicolumn{1}{c}{\quad} \\
\hline
$n \setminus k$ & 0 & 1 & 2 & 3 & 4 & 5 & 6 & 7 & Row sums \\
\hline
 0 & 1 &  &  &  &  &  &  &  & 1 \\
 1 & 2 &  &  &  &  &  &  &  & 2 \\
 2 & 34 & 6 &  &  &  &  &  &  & 40 \\
 3 & 1844 & 372 & 24 &  &  &  &  &  & 2240 \\
 4 & 199828 & 42864 & 3588 & 120 &  &  &  &  & 246400 \\
 5 & 36056936 & 8002992 & 748728 & 35424 & 720 &  &  &  & 44844800 \\
 6 & 9752801896 & 2212167336 & 220309896 & 12130056 & 371376 & 5040 &  &  & 12197785600 \\
 7 & 3691552813712 & 849994084272 & 88121628912 & 5290935792 & 194847552 & 4177440 & 40320 &  & 4635158528000 \\
\hline
\end{tabular}
\end{table}
%

\bigskip
\bigskip

\begin{table}[H]
\hspace*{-1.8cm}
\scriptsize
\begin{tabular}{c|rrrrrrr|r|}
\cline{2-8}
   & \multicolumn{7}{|c|}{\small $\bm{r=4}$} &  \multicolumn{1}{c}{\quad} \\
\hline
$n \setminus k$ & 0 & 1 & 2 & 3 & 4 & 5 & 6 & Row sums \\
\hline
 0 & 1 &  &  &  &  &  &  & 1 \\
 1 & 6 &  &  &  &  &  &  & 6 \\
 2 & 1236 & 24 &  &  &  &  &  & 1260 \\
 3 & 1229376 & 17904 & 120 &  &  &  &  & 1247400 \\
 4 & 3366706176 & 38473488 & 221616 & 720 &  &  &  & 3405402000 \\
 5 & 19614640553136 & 183482227008 & 881706816 & 2736000 & 5040 &  &  & 19799007228000 \\
 6 & 208727896045756896 & 1649611318980672 & 6725066986368 & 18337857984 & 35105760 & 40320 &  & 210384250804728000 \\
\hline
\end{tabular}
\end{table}
%

\bigskip
\bigskip

\begin{table}[H]
\hspace*{-1.8cm}
\scriptsize
\begin{tabular}{c|rrrrrr|r|}
\cline{2-7}
   & \multicolumn{6}{|c|}{\small $\bm{r=5}$} &  \multicolumn{1}{c}{\quad} \\
\hline
$n \setminus k$ & 0 & 1 & 2 & 3 & 4 & 5 & Row sums \\
\hline
 0 & 1 &  &  &  &  &  & 1 \\
 1 & 24 &  &  &  &  &  & 24 \\
 2 & 72456 & 120 &  &  &  &  & 72576 \\
 3 & 1742235984 & 1329120 & 720 &  &  &  & 1743565824 \\
 4 & 162123744912336 & 69701970960 & 20323440 & 5040 &  &  & 162193467211776 \\
 5 & 41351875243477668864 & 11349535075620480 & 2003358856320 & 303776640 & 40320 &  & 41363226782215962624 \\
\hline
\end{tabular}
\end{table}

%

\bigskip
\bigskip

\subsection{Higher-order quasi-Eulerian subset numbers}
   \label{app.eq.subset}

%
%
\begin{table}[H]
\centering
\small
\begin{tabular}{c|rrrrrrrrr|r|}
\cline{2-10}
   & \multicolumn{9}{|c|}{\small $\bm{r=1}$} &  \multicolumn{1}{c}{\quad} \\
\hline
$n \setminus k$ & 0 & 1 & 2 & 3 & 4 & 5 & 6 & 7 & 8 & Row sums \\
\hline
0 & 1 &  &  &  &  &  &  &  &  & 1  \\
1 & 1 & 0 &  &  &  &  &  &  &  & 1  \\
2 & 0 & 1 & 0 &  &  &  &  &  &  & 1  \\
3 & $-1$ & 1 & 1 & 0 &  &  &  &  &  & 1  \\
4 & 1 & $-5$ & 4 & 1 & 0 &  &  &  &  & 1  \\
5 & 2 & 1 & $-14$ & 11 & 1 & 0 &  &  &  & 1  \\
6 & $-9$ & 36 & $-29$ & $-24$ & 26 & 1 & 0 &  &  & 1  \\
7 & 9 & $-104$ & 281 & $-244$ & 1 & 57 & 1 & 0 &  & 1  \\
8 & 50 & $-83$ & $-454$ & 1401 & $-1259$ & 225 & 120 & 1 & 0 & 1  \\
\hline
\end{tabular}
\end{table}

\bigskip
\bigskip

\begin{table}[H]
\centering
\small
\begin{tabular}{c|rrrrrrrrr|r|}
\cline{2-10}
   & \multicolumn{9}{|c|}{\small $\bm{r=2}$ \bf (reversed second-order Eulerian numbers) } &  \multicolumn{1}{c}{\quad} \\
\hline
$n \setminus k$ & 0 & 1 & 2 & 3 & 4 & 5 & 6 & 7 & 8 & Row sums \\
\hline
 0 & 1 &  &  &  &  &  &  &  &  & 1 \\
 1 & 1 &  &  &  &  &  &  &  &  & 1 \\
 2 & 2 & 1 &  &  &  &  &  &  &  & 3 \\
 3 & 6 & 8 & 1 &  &  &  &  &  &  & 15 \\
 4 & 24 & 58 & 22 & 1 &  &  &  &  &  & 105 \\
 5 & 120 & 444 & 328 & 52 & 1 &  &  &  &  & 945 \\
 6 & 720 & 3708 & 4400 & 1452 & 114 & 1 &  &  &  & 10395 \\
 7 & 5040 & 33984 & 58140 & 32120 & 5610 & 240 & 1 &  &  & 135135 \\
 8 & 40320 & 341136 & 785304 & 644020 & 195800 & 19950 & 494 & 1 &  & 2027025 \\
\hline
\end{tabular}
\end{table}


\bigskip
\bigskip

\begin{table}[H]
\hspace*{-1.5cm}
\scriptsize
\begin{tabular}{c|rrrrrrrrr|r|}
\cline{2-10}
   & \multicolumn{9}{|c|}{\small $\bm{r=3}$} &  \multicolumn{1}{c}{\quad} \\
\hline
$n \setminus k$ & 0 & 1 & 2 & 3 & 4 & 5 & 6 & 7 & 8 & Row sums \\
\hline
 0 & 1 &  &  &  &  &  &  &  &  & 1 \\
 1 & 1 &  &  &  &  &  &  &  &  & 1 \\
 2 & 9 & 1 &  &  &  &  &  &  &  & 10 \\
 3 & 246 & 33 & 1 &  &  &  &  &  &  & 280 \\
 4 & 13390 & 1921 & 88 & 1 &  &  &  &  &  & 15400 \\
 5 & 1211386 & 180036 & 9771 & 206 & 1 &  &  &  &  & 1401400 \\
 6 & 164131730 & 24931166 & 1486131 & 40921 & 451 & 1 &  &  &  & 190590400 \\
 7 & 31103704820 & 4795137550 & 303467476 & 9711671 & 153531 & 951 & 1 &  &  & 36212176000 \\
 8 & 7854121032724 & 1223909199718 & 80747636454 & 2846874725 & 55244660 & 537756 & 1962 & 1 &  & 9161680528000 \\
\hline
\end{tabular}
\end{table}


\bigskip
\bigskip

\begin{table}[H]
\hspace*{-1.8cm}
\scriptsize
\begin{tabular}{c|rrrrrrrr|r|}
\cline{2-9}
   & \multicolumn{8}{|c|}{\small $\bm{r=4}$} &  \multicolumn{1}{c}{\quad} \\
\hline
$n \setminus k$ & 0 & 1 & 2 & 3 & 4 & 5 & 6 & 7 & Row sums \\
\hline
 0 & 1 &  &  &  &  &  &  &  & 1 \\
 1 & 1 &  &  &  &  &  &  &  & 1 \\
 2 & 34 & 1 &  &  &  &  &  &  & 35 \\
 3 & 5650 & 124 & 1 &  &  &  &  &  & 5775 \\
 4 & 2582915 & 44376 & 333 & 1 &  &  &  &  & 2627625 \\
 5 & 2510663365 & 35275610 & 228861 & 788 & 1 &  &  &  & 2546168625 \\
 6 & 4456094293397 & 52872120317 & 297244421 & 974995 & 1744 & 1 &  &  & 4509264634875 \\
 7 & 13054985706631155 & 133950756253880 & 660603311240 & 1987086224 & 3723163 & 3712 & 1 &  & 13189599057009375 \\
\hline
\end{tabular}
\end{table}
%

\bigskip
\bigskip

\begin{table}[H]
\hspace*{-1.5cm}
\scriptsize
\begin{tabular}{c|rrrrrrr|r|}
\cline{2-8}
   & \multicolumn{7}{|c|}{\small $\bm{r=5}$} &  \multicolumn{1}{c}{\quad} \\
\hline
$n \setminus k$ & 0 & 1 & 2 & 3 & 4 & 5 & 6 & Row sums \\
\hline
 0 & 1 &  &  &  &  &  &  & 1 \\
 1 & 1 &  &  &  &  &  &  & 1 \\
 2 & 125 & 1 &  &  &  &  &  & 126 \\
 3 & 125665 & 460 & 1 &  &  &  &  & 126126 \\
 4 & 487856621 & 1006503 & 1251 & 1 &  &  &  & 488864376 \\
 5 & 5187834064414 & 6833491661 & 5300301 & 2999 & 1 &  &  & 5194672859376 \\
 6 & 123266182967274060 & 112433009305765 & 59083404656 & 23048178 & 6716 & 1 &  & 123378675083039376 \\
\hline
\end{tabular}
\end{table}


\bigskip
\bigskip

\section{Report on systematic verification of conjectures}
\label{app.sec.report}


We will now report the results of our systematic experiments, 
i.e., we will write out up to which order we have verified each of our conjectures.

Our total positivity tests were performed by using the following methods:
\begin{itemize}
	\item For a finite matrix $A$ with real entries, 
		an efficient method to test total positivity is 
		by computing the bidiagonal factorization
		using Neville elimination \cite{Gasca_92}.
		We used this technique to test the total positivity
                of lower-triangular matrices.

	\item For checking negative-real-rootedness of polynomials, 
		we used {\sc Mathematica}'s built-in {\tt NSolve} command
                to list out approximate values of the zeros. 
		This method is very efficient.

	\item For checking coefficientwise Hankel-total posivity of a sequence of polynomials,
		we unfortunately have to compute and check all minors of the Hankel matrix.
		No other efficient method is known.
		Compared to the other two types of tests, 
		this method is very slow and we can check this only for matrices of small sizes.

\end{itemize}

We now present the reports for
Conjectures~\ref{conj.cycle}, \ref{conj.subset},
\ref{conj.quasi-Eulerian.TP}, \ref{conj.quasi-Eulerian.hankel},
and \ref{conj.phylogenetic.ordered}.

\bigskip

\noindent{\bf Conjecture~\ref{conj.cycle}:}

The order up to which the various conjectures in Conjecture~\ref{conj.cycle} 
have been verified is recorded in Table~\ref{tab.conj.cyc.ver}.
We have checked that $\widecheck{C}^{(r)}$ is not totally positive for $r\geq 3$.
For $r=3$, the $6\times 6$ top left submatrix fails to be totally positive.
For $r\geq 4$, we in fact have a $2\times 2$ minor which is negative;
it is 
\be
        \stirlingcycle{2}{2}^{(r)}\stirlingcycle{3}{2}^{(r)}
                -
        \stirlingcycle{2}{1}^{(r)}\stirlingcycle{3}{3}^{(r)}\;.
\ee

\bigskip

\noindent{\bf Conjecture~\ref{conj.subset}:}

The order up to which the various conjectures in Conjecture~\ref{conj.subset}
have been verified is recorded in Table~\ref{tab.conj.subset.ver}.
We have checked that $\widecheck{S}^{(r)}$ is not totally positive for $r\geq 3$.
For $r=3$, the $13\times 13$ top left submatrix fails to be totally positive.
For $r\geq 4$, we in fact have a $2\times 2$ minor which is negative;
it is
\be
        \stirlingsubset{2}{2}^{(r)}\stirlingsubset{3}{2}^{(r)}
                -
        \stirlingsubset{2}{1}^{(r)}\stirlingsubset{3}{3}^{(r)}\;.
\ee

For $r\geq 3$, the Hankel matrices $(s_{r,n+k}(x))_{n,k\geq 0}$ fail to be 
coefficientwise totally positive.
In fact, it seems that for $r\geq 3$,
the $3\times 3$ minor $|(s_{r,n+k}(x))_{1\leq n \leq 3, 0\leq k \leq 2}|$
is has some negative coefficients.

\setcounter{topnumber}{3}
%
\begin{table}[t]
\begin{tabular}{|c|c|c|}
\hline
Conjecture & $r$ & Order up to which verified\\
\hline
(a) Total positivity of $C^{(r)}$ & $r=2,\ldots,10$ & $70\times 70$ \\
\hline
(b) Total positivity of $\widecheck{C}^{(r)}$ & $r=2$ & $90\times 90$ \\
\hline
(c) Log-concavity of coefficients	&&\\
 of $c_{r,n}(x)$ & $r=3,4,5$ & $1000$ \\
\hline
(d) Coefficientwise Hankel-total	&&\\
positivity of $(c_{r,n}(x))_{n\geq 0}$
	& $r=3,\ldots,7$
	& $11\times 11$\\
	& $r=8,9,10$&$10\times 10$\\
\hline
\end{tabular}
\caption{Systematic report for verification of Conjecture~\ref{conj.cycle}.}
\label{tab.conj.cyc.ver}
\end{table}
\begin{table}[t]
\begin{tabular}{|c|c|c|}
\hline
Conjecture & $r$ & Order up to which verified\\
\hline
(a) Total positivity of $S^{(r)}$ & $r=2,\ldots,10$ & $70\times 70$ \\
\hline
(b) Total positivity of $\widecheck{S}^{(r)}$ & $r=2$ & $70\times 70$ \\
\hline
(c) Log-concavity of coefficients       &&\\
 of $s_{r,n}(x)$ & $r=3,4,5$ & $1000$ \\
\hline
\end{tabular}
\caption{Systematic report for verification of Conjecture~\ref{conj.subset}.}
\label{tab.conj.subset.ver}
\end{table}
\begin{table}[t]
\centering
\begin{tabular}{|c|c|}
\hline
Conjecture &  Order up to which verified\\
\hline
(a) Total positivity of $D$  & $100\times 100$ \\
\hline
(b) Total positivity of $D^{\rm rev}$ & $80\times 80$ \\
\hline
(c) Coefficientwise Hankel-total &\\
positivity of $(d_n(x))_{n\geq 0}$
        & $10\times 10$\\
\hline
\end{tabular}
\caption{Systematic report for verification of Conjecture~\ref{conj.phylogenetic.ordered}.}
\label{tab.conj.phylogenetic.ordered}
\end{table}

\bigskip

\noindent {\bf Conjectures~\ref{conj.quasi-Eulerian.TP} and
                           \ref{conj.quasi-Eulerian.hankel}:}

Conjectures~\ref{conj.quasi-Eulerian.TP}(b) and (b${}'$)
have been tested for $r=2,\ldots, 10$
for the $50\times 50$ leading principal submatrix.
Conjecture~\ref{conj.quasi-Eulerian.hankel}
has been tested for $r=3,\ldots, 10$
for the $9\times 9$ leading principal submatrix.

%

\bigskip

\noindent {\bf Conjecture~\ref{conj.phylogenetic.ordered}:}

The order up to which the various conjectures in Conjecture~\ref{conj.phylogenetic.ordered}
have been verified is recorded in Table~\ref{tab.conj.phylogenetic.ordered}.

\end{document}